\documentclass[11pt]{amsart}
\usepackage{amscd,amssymb}
\usepackage{amsthm,amsmath,amssymb}
\usepackage{youngtab}
\usepackage{tikz}
\usetikzlibrary{calc}
\usepackage{graphicx}
\usepackage{color}
\usepackage{pdflscape}
\usepackage{subfigure}
\usepackage{bm}
\usepackage{mathrsfs}
\usepackage{verbatim}
\usepackage{braids}
\usepackage{adjustbox, collectbox, xkeyval, graphicx, ifpdf, adjcalc}
\newcommand{\ms}{\medskip}

\usepackage{amsmath,amsthm,verbatim, color, amssymb,amsfonts,amscd, graphicx, hyperref, enumerate}
\usepackage{graphics}
\usepackage{tikz, tikz-cd}
\usetikzlibrary{matrix}
\usepackage[all]{xy}
\usepackage{listings}
\theoremstyle{plain}
\newtheorem{thm}{Theorem}[section]
\newtheorem{lem}[thm]{Lemma}        
\newtheorem{prop}[thm]{Proposition}

\newtheorem*{thm*}{Theorem}
\newtheorem{abcthm}{Theorem}
\newtheorem{abccor}[abcthm]{Corollary}

\theoremstyle{definition}
\newtheorem{defn}[thm]{Definition} 
\newtheorem*{defn*}{Definition}   
   
\newtheorem*{exam*}{Example}  
\newtheorem{exam}[thm]{Example}  
\newtheorem{rem}[thm]{Remark}

\newcommand{\bb}{\mathbin{\backslash \mkern-6mu \backslash}}
\newcommand{\ff}{\mathbin{/\mkern-6mu/}}
\newcommand{\cd}{\partial}
\newcommand{\sm}{\sigma}
\newcommand{\Dl}{\Delta}
\newcommand{\balpha}{\overline{\alpha}}
\newcommand{\btau}{\overline{\tau}}

\newcommand{\lng}{\ell}

\newcommand{\Z}{\mathbb{Z}}
\newcommand{\D}{\mathcal{D}}
\newcommand{\A}{\mathcal{A}}
\newcommand{\C}{\mathcal{C}}
\newcommand{\B}{\mathcal{B}}

\newcommand{\bt}{\bullet}
\newcommand{\lb}{\lbrack \! \lbrack}
\newcommand{\rb}{\rbrack \! \rbrack}
\newcommand{\ba}{\bar{a}}

\newcommand{\altprod}{\pi}
\newcommand{\Aforheader}{\A^n_\bt}
\newcommand{\EndMon}{\operatorname{EndMon}}
\newcommand{\EndGen}{\operatorname{EndGen}}
\newcommand{\id}{\operatorname{id}}
\newcommand{\colim}{\operatorname{colim}}

\definecolor{amethyst}{rgb}{0.6, 0.4, 0.8}
\definecolor{amber}{rgb}{1.0, 0.49, 0.0}
\definecolor{brightmaroon}{rgb}{0.76, 0.13, 0.28}
\definecolor{darkpastelgreen}{rgb}{0.01, 0.75, 0.24}

\usepackage[
paper=a4paper,
headsep=15pt,text={138mm,213mm},centering,includehead
]{geometry}

\usepackage{lipsum}

\begin{document}

\title{Homological stability for Artin monoids}
\author{Rachael Boyd}
\address{Max Planck Institute for Mathematics, Bonn}
\email{rachaelboyd@mpim-bonn.mpg.de}

\def\subjclassname{\textup{2010} Mathematics Subject Classification}
\expandafter\let\csname subjclassname@1991\endcsname=\subjclassname
\expandafter\let\csname subjclassname@2000\endcsname=\subjclassname
\subjclass{
57M07, % Topological methods in group theory
20F36, %Artin and braid groups
20F55, %Reflection and Coxeter groups
55U10, %Simplicial sets and complexes
20M32. %Algebraic monoids
}
\keywords{Artin monoids, homological stability, Artin groups, $K(\pi,1)$ conjecture.}

\maketitle

\begin{abstract}
	We prove that certain sequences of Artin monoids containing the braid monoid as a submonoid satisfy homological stability. When the $K(\pi,1)$ conjecture holds for the associated family of Artin groups this establishes homological stability for these groups. In particular, this recovers and extends Arnol'd's proof of stability for the Artin groups of type $\textbf{A}$, $\textbf{B}$ and $\textbf{D}$.
\end{abstract}

\section{Introduction}

A sequence of groups or monoids with maps between them
	\begin{equation*}
	G_1 \rightarrow G_2 \rightarrow \cdots \rightarrow G_n \rightarrow \cdots	
	\end{equation*}
is said to satisfy {\it homological stability} if the induced maps on homology
	$$H_i(G_n) \to H_i(G_{n+1})$$
are isomorphisms for $n$ sufficiently large compared to $i$.

This paper concerns homological stability for sequences of \emph{Artin monoids and groups}, and in this paper the associated maps will always be inclusions. In particular we consider sequences of Artin groups that have the braid group as a subgroup, and the corresponding sequences of monoids.

We recall the definition of Artin groups. Given a finite set $\Sigma$, to every unordered pair $\{\sigma_s,\sigma_t\} \in \Sigma \times \Sigma$ associate either a natural number greater than 2 or the symbol $\infty$, and denote this by $m(s,t)$. An Artin group $A$ with generating set $\Sigma$ has the following presentation
$$
	A=\langle \,\Sigma \, \vert \,  \altprod(\sm_s,\sm_t; m(s,t))=\altprod(\sm_t,\sm_s; m(s,t))\rangle,
$$
where $\altprod(\sm_s,\sm_t; m(s,t))$ is the alternating product of $\sm_s$ and $\sm_t$ starting with $\sm_s$ and of length $m(s,t)$. The braid group with its standard presentation is the archetypal example of an Artin group, with presentation
	$$
	\B_n=\Big\langle \sm_i \,\text{ for } 1\leq i\leq n-1 \,\mid\,\begin{array}{lr}\sm_i\sm_j=\sm_j\sm_i & |i-j|\geq 2\\ \sm_i\sm_{i+1}\sm_i=\sm_{i+1}\sm_i\sm_{i+1} &1\leq i \leq n-2\end{array}\Big\rangle.
	$$
Every Artin group has an associated \emph{Coxeter group} (discussed in Section \ref{chap: BG Artin groups}), and in fact Artin groups were first introduced by Brieskorn \cite{Brieskorn2} as the fundamental groups of hyperplane complements built from Coxeter groups. The information of the presentation can be packaged into a \emph{Coxeter diagram}. This diagram has vertex set $\Sigma$ and edges corresponding to $m(s,t)$ for each pair of vertices: no edge when $m(s,t)=2$, an unlabelled edge when $m(s,t)=3$ and an edge labelled with $m(s,t)$ otherwise. For example the braid group $\B_n$ has diagram
	\begin{center}
		\begin{tikzpicture}[scale=0.2, baseline=0]
		\draw[fill= black] (5,0) circle (0.4) node[below] {$\sm_1$};	
		\draw[line width=1] (5,0) -- (15,0);
		\draw[fill= black] (10,0) circle (0.4) node[below] {$\sm_2$};
		\draw[fill= black] (15,0) circle (0.4) node[below] {$\sm_3$};		
		\draw (17.5,0) node {$\ldots$};
		\draw[fill= black] (20,0) circle (0.4) node[below] 	{$\sm_{n-2}$};
		\draw[line width=1] (20,0) -- (25,0);				
		\draw[fill= black] (25,0) circle (0.4) node[below] {$\sm_{n-1}$};
		\end{tikzpicture}
	\end{center}
which is known as the Coxeter diagram of \emph{type} $\textbf{A}$ and has corresponding Coxeter group the symmetric group $S_n$.

The sequences of Artin groups studied in this paper correspond to the following sequence of diagrams:

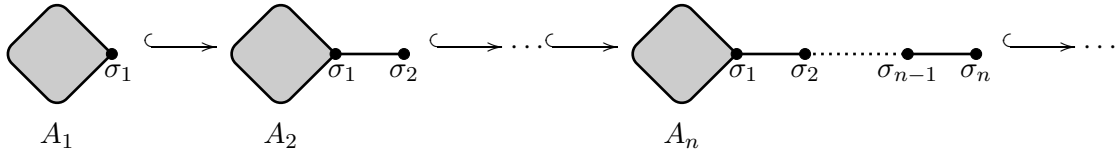
\centerline{	\xymatrix@R=3mm {
		\begin{tikzpicture}[scale=0.18, baseline=0]
		\draw[line width=1, fill = white!80!black, rounded corners=5 pt]
		(5,0) --  (1,4) -- (-3,0) -- (1,-4) -- (5,0);
		\draw (1,-6) node {$A_1$};
		\draw[fill= black] (5,0) circle (0.4);
		\draw (5.5,0) node[below] {$\sm_1$};
		\end{tikzpicture} \ar@{^{(}->}[r]
		&
		\begin{tikzpicture}[scale=0.18, baseline=0]
		\draw[line width=1, fill = white!80!black, rounded corners=5 pt]
		(5,0) --  (1,4) -- (-3,0) -- (1,-4) -- (5,0);
		\draw[line width=1] (5,0) -- (10,0);
		\draw (1,-6) node {$A_2$};
		\draw[fill= black] (5,0) circle (0.4);		
		\draw (5.5,0) node[below] {$\sm_1$};
		\draw[fill= black] (10,0) circle (0.4) node[below] {$\sm_2$};
		\end{tikzpicture}  \ar@{^{(}->}[r]
		&
		\cdots  \ar@{^{(}->}[r]
		&
		\begin{tikzpicture}[scale=0.18, baseline=0]
		\draw[line width=1, fill = white!80!black, rounded corners=5 pt]
		(5,0) --  (1,4) -- (-3,0) -- (1,-4) -- (5,0);
		\draw[fill= black] (5,0) circle (0.4);
		\draw (5.5,0) node[below] {$\sm_1$};
		\draw[line width=1] (5,0) -- (10,0);
		\draw[fill= black] (10,0) circle (0.4) node[below] {$\sm_2$};
		\draw[line width=1, dotted] (10,0) -- (17.5,0);
		\draw[fill= black] (17.5,0) circle (0.4)node[below] {$\sm_{n-1}$};
		\draw[line width=1] (17.5,0) -- (22.5,0);
		\draw[fill= black] (22.5,0) circle (0.4)node[below] {$\sm_{n}$};
		\draw (1,-6) node {$A_n$};
		\end{tikzpicture}  \ar@{^{(}->}[r]
		&\cdots
}}

\noindent where the grey box indicates that the sequence begins with an arbitrary diagram: an arbitrary Artin group with finite generating set. The type $\textbf{A}$ subdiagram corresponds to a subgroup of $A_n$ being the braid group $\B_{n+1}$, with an increasing number of generators as~$n$ increases. This gives rise to a sequence of groups and inclusions
$$
A_1 \hookrightarrow A_2 \hookrightarrow \cdots \hookrightarrow A_n \hookrightarrow \cdots
$$
\noindent and the goal of this paper is to discuss stability for sequences of Artin groups of this form. This was motivated by work of Hepworth \cite{Hepworth}, who proved homological stability for the associated sequence of Coxeter groups.

While the argument used for the proof of homological stability is similar to that used by several authors, the novel part of this paper comes from dealing with Artin groups and monoids. Very little is known for Artin groups in general, for instance the centre of a generic Artin group is unknown and it is not known whether all Artin groups are torsion free. In particular there are no tools to date for working with Artin cosets (for example there is no canonical way to choose a coset representative), something that is usually desirable when proving homological stability for a family of groups. Therefore the results of this paper are stated and proved for the corresponding \emph{Artin monoids}, for which a technical `coset' theory is developed in Section \ref{chap:BG Artin monoids} (the notion of coset of a submonoid is not defined in general). Key properties of Artin monoids, such as the existence of a well defined length function and the existence of lowest common multiples under certain conditions, allow us to define a canonical choice of `coset representative' for these monoids. From the monoid result we then deduce homological stability for  Artin groups that satisfy the $K(\pi, 1)$ conjecture (discussed in more detail below).

We denote the Artin monoid corresponding to $A_n$ by $A_n^+$. The inclusion map between the monoids is denoted $s$ and called the \emph{stabilisation map}. The main result of this paper is the following, which to my knowledge is the first instance where homological stability is proved for monoids.

\begin{abcthm} \label{thm:A}
	The sequence of Artin monoids
	\begin{equation*}
	A^+_1 \hookrightarrow A^+_2 \hookrightarrow \cdots \hookrightarrow A^+_n \hookrightarrow \cdots	
	\end{equation*}
	satisfies homological stability. More precisely, the induced map on homology
	$$
	H_*(BA^+_{n-1})\overset{s_*}{\longrightarrow} H_*(BA^+_n)
	$$
	is an isomorphism when $*<\frac{n}{2}$ and a surjection when $*=\frac{n}{2}$. Here, homology is taken with arbitrary constant coefficients i.e.~coefficients in an abelian group.
\end{abcthm}

The classifying space of an Artin monoid $BA^+$ is homotopy equivalent to some interesting spaces that arise naturally in mathematics. One manifestation of this is that in the study of Artin groups, there is a well known conjecture by Arnol'd, Brieskorn, Pham and Thom called the $K(\pi,1)$ conjecture (discussed in Section \ref{sec:K(pi,1)conjecture}). For this introduction it suffices to know the following fact, due to Dobrinskaya.

\begin{thm*}[{\cite[Theorem 6.3]{Dobrinskaya2006}}]\label{thm:monoid relation to kpi1 dobrinskaya}
	Given an Artin group $A$ and its associated monoid $A^+$, the $K(\pi,1)$ conjecture holds if and only if the induced map between their classifying spaces $BA^+ \to BA$ is a homotopy equivalence.
\end{thm*}

\noindent Thus if the $K(\pi,1)$ conjecture holds for a family of Artin groups, Theorem \ref{thm:A} establishes homological stability as below.

\begin{abccor} \label{corollary: result for groups}
	When the $K(\pi, 1)$ conjecture holds for all $A_n$, the sequence of Artin groups
	\begin{equation*}
	A_1 \hookrightarrow A_2 \hookrightarrow \cdots \hookrightarrow A_n \hookrightarrow \cdots	
	\end{equation*}
	satisfies homological stability. More precisely, the induced map on homology
	$$
	H_*(BA_{n-1})\to H_*(BA_n)
	$$is an isomorphism when $*<\frac{n}{2}$ and a surjection when $*=\frac{n}{2}$. Here, homology is taken with arbitrary constant coefficients i.e.~coefficients in an abelian group.
\end{abccor}

The $K(\pi,1)$ conjecture has been proven for large classes of Artin groups \cite{Paris}: the conjecture holds for Artin groups for which the corresponding Coxeter groups are \emph{finite} (this is \emph{Deligne's Theorem} \cite{Deligne}, Theorem \ref{thm:delignes}), of \emph{large type}, of \emph{dimension two} and of \emph{FC type}. Proving the~$K(\pi,1)$ conjecture for families of Artin groups continues to be an active area of research to this day.

\begin{abccor}\label{cor:C}
	Homological stability holds for the sequences of Artin groups~$(A_n)$ for which the corresponding Coxeter groups are either finite, of large type, of dimension two or of FC type.
\end{abccor}

\begin{rem}\label{rem:finite egs}
	In the case of finite Coxeter groups, Corollary \ref{cor:C} therefore recovers the few known cases of stability for families of Artin groups of the form studied in this paper, i.e.~homological stability holds for the sequences of Artin groups~$\{A_n\}_{n\geq 1}$ of type $\textbf{A}$, $\textbf{B}$ and $\textbf{D}$, given by the following diagrams:
	
		
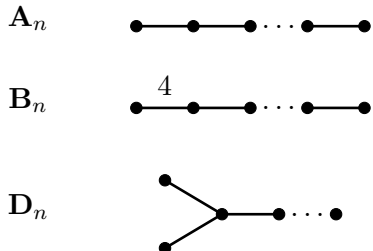
\centerline{	\xymatrix@R=3mm {	
			\textbf{A}_n & \begin{tikzpicture}[scale=0.15, baseline=0]
			\draw[fill= black] (5,0) circle (0.5);
			\draw[line width=1] (5,0) -- (15,0);
			\draw[fill= black] (10,0) circle (0.5);
			\draw[fill= black] (15,0) circle (0.5);		
			\draw (17.5,0) node {$\ldots$};
			\draw[fill= black] (20,0) circle (0.5);
			\draw[line width=1] (20,0) -- (25,0);				
			\draw[fill= black] (25,0) circle (0.5);
			\end{tikzpicture}
			  \\
			\textbf{B}_n & \begin{tikzpicture}[scale=0.15, baseline=0]
			\draw[fill= black] (5,0) circle (0.5);
			\draw[line width=1] (5,0) -- (15,0);
			\draw[fill= black] (10,0) circle (0.5);
			\draw[fill= black] (15,0) circle (0.5);		
			\draw (17.5,0) node {$\ldots$};
			\draw[fill= black] (20,0) circle (0.5);
			\draw[line width=1] (20,0) -- (25,0);				
			\draw[fill= black] (25,0) circle (0.5);
			\draw (7.5,2) node {4};
			\end{tikzpicture}
			\\
			\textbf{D}_n & \begin{tikzpicture}[scale=0.15, baseline=0]
			\draw[fill= black] (5,-3) circle (0.5);
			\draw[line width=1] (5,-3) -- (10,0);
			\draw[line width=1] (5,3) -- (10,0);
			\draw[line width=1] (15,0) -- (10,0);
			\draw[fill= black] (5,3) circle (0.5);
			\draw[fill= black] (10,0) circle (0.5);
			\draw[fill= black] (15,0) circle (0.5);
			\draw (17.5,0) node {$\ldots$};				
			\draw[fill= black] (20,0) circle (0.5);
			\end{tikzpicture}
		}}
		These three sequences consist of Artin groups which relate to finite Coxeter groups Hence by Corollary \ref{cor:C}, the sequences of Artin groups satisfy homological stability.
\end{rem}

The three examples in Remark~\ref{rem:finite egs} were proved by Arnol'd, who computed the full (co)homology of the groups in question, using the associated hyperplane complement. The results and proofs are in Brieskorn's Bourbaki seminar~\cite{Brieskorn}. Despite the theorem in this paper generalising these results, the method of proof is not a straightforward generalisation of the proof of stability for the braid group.

\begin{rem}
	Krannich \cite{krannich} introduced a framework to study homological stability phenomena in the context of~$E_2$-algebras. This generalised a categorical framework of Randal-Williams and Wahl \cite{randalWilliamsWahl}, which ‘automated’ parts of the homological stability proof for sequences of discrete groups. The sequence of monoids studied in this paper does not fit into the categorical set up of \cite{randalWilliamsWahl}, however the classifying spaces of the sequence can be shown to assemble into an~$E_1$-module over an~$E_2$ algebra, as in \cite{krannich}. Combining the results of \cite{krannich} with the high connectivity results established in this paper, our homological stability result with constant coefficients (Theorem \ref{thm:A}) can most likely be enhanced to one with abelian coefficients and coefficient systems of finite degree in the sense of \cite[Section 4]{krannich}.
\end{rem}

\subsection{Outline of proof}\label{section:setup}

The proof follows the outline of a standard homological stability argument, which we describe below for the benefit of the reader. We indicate where the new ingredients are used. 

The proof of Theorem \ref{thm:A} requires the introduction of a semi-simplicial space $\A^n_\bt$ for each monoid in the sequence $A^+_n$ such that:

\begin{enumerate}
	\item there exist homotopy equivalences $\A^n_p \simeq BA^+_{n-p-1}$ for every $p \geq 0$; and \label{point 2}
	\item there is a highly connected map from the geometric realisation of $\A^n_\bt$ to the classifying space $BA^+_n$, which we denote $\|\phi_\bt\|$
	$$
	\| \A^n_\bt \| \overset{\|\phi_\bt \|}{\longrightarrow} BA^+_n
	$$
	i.e.~ $\|\phi_\bt \|$ induces an isomorphism on a large range of homotopy groups. \label{point 4}
\end{enumerate}

The skeletal filtration of $\|\A^n_\bt\|$ gives rise to a spectral sequence

$$
E^1_{p,q}=H_q(\A^n_p)\Rightarrow H_{p+q}(\|\A_\bt^n\|).
$$
From Point (\ref{point 2}) it follows:
$$
E^1_{p,q}=H_q(\A^n_p)=H_q(BA^+_{n-p-1}).
$$

\noindent We prove that on the $E^1$ page under the above equality the differentials are given by either the zero map or the stabilisation map $s_* \colon H_q(BA^+_{n-p-1}) \to H_q(BA^+_{n-p})$. Following this, applying the inductive hypotheses that previous monoids in the sequence satisfy stability gives that in a range (when $q$ is small compared to $n$) the spectral sequence converges to $H_q(BA^+_{n-1})$:

$$
E^1_{p,q}=H_q(BA^+_{n-p-1})\Rightarrow H_{p+q}(\|\A_\bt^n\|)=H_q(BA^+_{n-1}) \text{ in a range}.
$$
The highly connected map of Point (\ref{point 4}) above now gives that in a range the spectral sequence also converges to the homology of $BA_n^+$, which completes the proof.

For sequences of discrete groups, a usual candidate for $\A^n_\bt$ would be built out of cosets of previous groups in the sequence. However the fact that no tools exist for manipulating Artin cosets means that this approach cannot be taken. The coset theory developed in this paper for the corresponding sequence of Artin monoids is used to build $\A^n_\bt$. The main obstacle in the proof is the high connectivity argument for Point \ref{point 4} which follows a `union of chambers' argument inspired by, but more involved than, work of Paris \cite{Paris} and Davis \cite{Davis}.

\subsection{Organisation of the paper}
Sections \ref{chap: BG Artin groups}, and \ref{sec:K(pi,1)conjecture} provide background on Coxeter groups and Artin groups, and the $K(\pi, 1)$ conjecture respectively. Section \ref{chap:BG Artin monoids} then introduces Artin monoids and develops a novel theory of `cosets' and corresponding technical results. Following this Section \ref{sec: semi simp constructions with monoids} details the required semi-simplicial background and particular monoid constructions used in the proof, some of which are new. Section \ref{section:artin} applies the theory of Section \ref{chap:BG Artin monoids}, and introduces notation used throughout the proof. Section \ref{sec:introduce tri ssspace} introduces the semi-simplicial space $\A^n_
\bt$ and the map $\|\phi_\bt \|$ described above. High connectivity of $\|\phi_\bt \|$ is then the topic of Section \ref{section:highcon}, in which the general method of proof for the high connectivity argument is introduced, before the proof is split into several cases, due to the complexity of using Artin monoids. Finally the spectral sequence argument and homological stability result are given in Section \ref{sec:HS proof}.

\subsection*{Acknowledgments}
I would like to thank my PhD advisor Richard Hepworth for his wonderful teaching, our many maths conversations and his unconditional support. I would also like to thank my thesis examiners Mark Grant and Ulrike Tillman for their careful reading and related comments and suggestions. Finally I would like to thank the anonymous referee, whose helpful comments significantly improved the exposition of this work.

\section{Coxeter groups and Artin groups}\label{chap: BG Artin groups}

\subsection{Coxeter groups}
This section follows \emph{The Geometry and Topology of Coxeter Groups} by Davis~\cite{Davis}.

\begin{defn}
	A \emph{Coxeter matrix} on a finite set of generators $S$ is a symmetric matrix $M$ indexed by elements of $S$, i.e.~with entries $m(s,t)$ in $\mathbb{N}\cup \infty$ for $\{s, t\}$ in~$S\times S$. This matrix must satisfy
	\begin{itemize}
		\item $m(s,s)=1$ for all $s$ in $S$
		\item $m(s,t)=m(t,s)$ must be either greater than 1, or $\infty$, when $s\neq t$.
	\end{itemize}
\end{defn}

\begin{defn}\label{def:Coxeter system}
	A Coxeter matrix $M$ with generating set $S$ has an associated \emph{Coxeter group} $W$, with presentation
	$$
	W=\langle S \, \vert \,  (st)^{m(s,t)}=e\rangle.
	$$
	\noindent Here $m(s,t)=\infty$ means there is no relation between $s$ and $t$. We call $(W,S)$ a \emph{Coxeter system}. We adopt the convention that $(W,\emptyset)$ is the trivial group.
\end{defn}

\begin{rem}
	Note that the condition $m(s,s)=1$ on the Coxeter matrix implies that the generators of the group are involutions i.e.~$s^2=e$ for all $s$ in $S$. %Therefore a Coxeter group is a group which is generated by reflections.
\end{rem}

\begin{defn}
	Define the \emph{length function }on a Coxeter system $(W,S)$
	$$\lng: W \to \mathbb{N}$$
	to be the function which maps $w$ in $W$ to the minimum word length required to express $w$ in terms of the generators.
\end{defn}

\begin{defn}\label{def: pi(s,t;m)}
	Define $\altprod(a,b ; k)$ to be the word of length $k$, given by the alternating product of $a$ and $b$ i.e.
	$$
	\altprod(a,b ; k)=\overbrace{abab\ldots }^{\text{length }k}
	$$
\end{defn}

\begin{rem}\label{rem:second pres for W}
	The relations $(st)^{m(s,t)}=e$ can be rewritten as
	$$
	\altprod(s,t;m(s,t))=\altprod(t,s;m(s,t))
	$$
	when $m(s,t)\neq \infty$. Therefore the presentation of a Coxeter group $W$ can also be given as
	$$
	W=\Big\langle S \, \Big\vert \, \begin{array}{cr}  (s)^2=e &  s\in S\\	\altprod(s,t;m(s,t))=\altprod(t,s;m(s,t)) & s, t \in S\end{array}\Big\rangle.
	$$
\end{rem}

\begin{defn}
	Given a Coxeter matrix corresponding to a Coxeter system $(W,S)$, there is an associated graph called the \emph{Coxeter diagram}, denoted $\D_W$. It is the graph with vertices indexed by the elements of the generating set $S$. Edges are drawn between the vertices corresponding to $s$ and $t$ in $S$ when $m(s,t)\geq 3$ and labelled with $m(s,t)$ when $m(s,t)\geq 4$, as shown below:
	
	\centerline{\xymatrix@R=2mm{
			\begin{tikzpicture}[scale=0.2, baseline=0]
			\draw[fill= black] (5,0) circle (0.4) node[below] {$s$};			
			\draw[fill= black] (10,0) circle (0.4) node[below] {$t$};
			\end{tikzpicture} &
			\begin{tikzpicture}[scale=0.2, baseline=0]
			\draw[fill= black] (5,0) circle (0.4) node[below] {$s$};
			\draw (7.5,1.2) node {};
			\draw[line width=1] (5,0) -- (10,0);			
			\draw[fill= black] (10,0) circle (0.4) node[below] {$t$};
			\end{tikzpicture} &
			\begin{tikzpicture}[scale=0.2, baseline=0]
			\draw[fill= black] (5,0) circle (0.4) node[below] {$s$};
			\draw (8.5,1.4) node {$m(s,t)$};
			\draw[line width=1] (5,0) -- (12,0);			
			\draw[fill= black] (12,0) circle (0.4) node[below] {$t$};
			\end{tikzpicture} \\
			m(s,t)=2 &m(s,t)=3 & m(s,t)\geq 4 \text{ or }\infty.			
	}}	
	\ms
	\noindent When the diagram $\D_W$ is connected, $W$ is called an \emph{irreducible} Coxeter group. The disjoint union of two diagrams gives the product of their corresponding Coxeter groups.
\end{defn}

\begin{thm}[{Classification of finite Coxeter groups, \cite{Coxeter}}]\label{prop:classification of finite Coxeter}
	A Coxeter group is finite if and only if it is a (direct) product of finitely many finite irreducible Coxeter groups.\\
	\noindent The following is a complete list of the diagrams corresponding to finite irreducible Coxeter groups.
	
	\ms
	
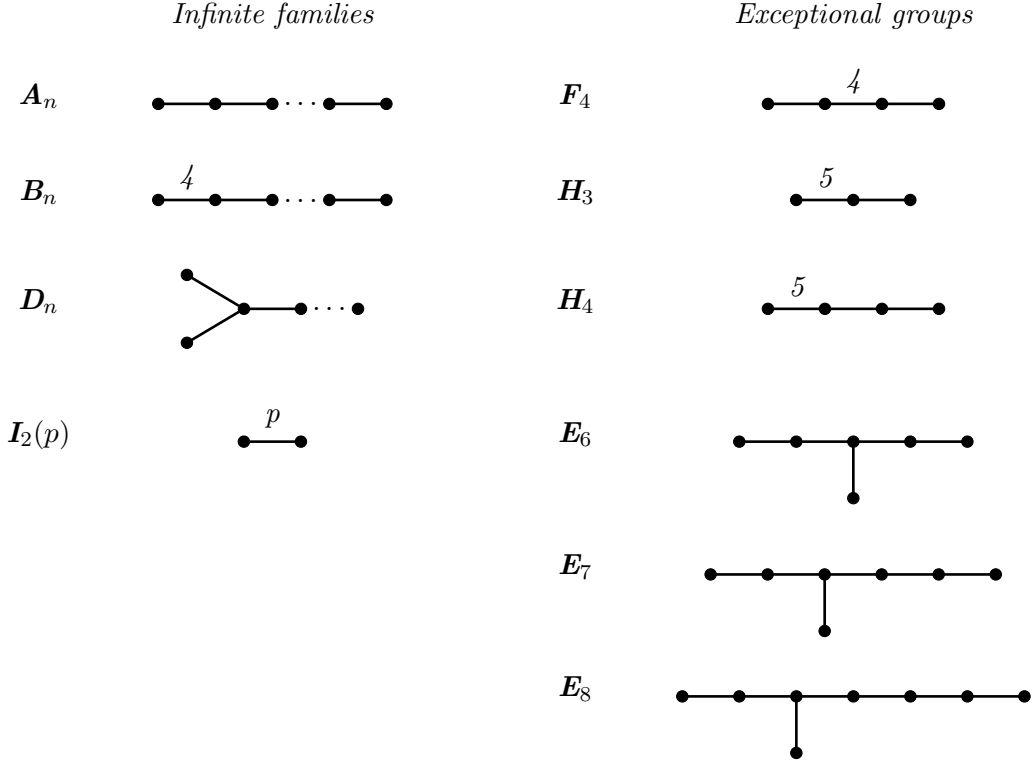
\centerline{	\xymatrix@R=3mm {	
			& \textrm{Infinite families} & & & \textrm{Exceptional groups}\\	
			\textbf{A}_n & \begin{tikzpicture}[scale=0.15, baseline=0]
			\draw[fill= black] (5,0) circle (0.5);
			\draw[line width=1] (5,0) -- (15,0);
			\draw[fill= black] (10,0) circle (0.5);
			\draw[fill= black] (15,0) circle (0.5);		
			\draw (17.5,0) node {$\ldots$};
			\draw[fill= black] (20,0) circle (0.5);
			\draw[line width=1] (20,0) -- (25,0);				
			\draw[fill= black] (25,0) circle (0.5);
			\end{tikzpicture}
			& & \textbf{F}_4 & \begin{tikzpicture}[scale=0.15, baseline=0]
			\draw[fill= black] (5,0) circle (0.5);
			\draw[line width=1] (5,0) -- (20,0);
			\draw[fill= black] (10,0) circle (0.5);
			\draw[fill= black] (15,0) circle (0.5);
			\draw[fill= black] (20,0) circle (0.5);
			\draw (12.5,2) node {4};
			\end{tikzpicture}  \\
			\textbf{B}_n & \begin{tikzpicture}[scale=0.15, baseline=0]
			\draw[fill= black] (5,0) circle (0.5);
			\draw[line width=1] (5,0) -- (15,0);
			\draw[fill= black] (10,0) circle (0.5);
			\draw[fill= black] (15,0) circle (0.5);		
			\draw (17.5,0) node {$\ldots$};
			\draw[fill= black] (20,0) circle (0.5);
			\draw[line width=1] (20,0) -- (25,0);				
			\draw[fill= black] (25,0) circle (0.5);
			\draw (7.5,2) node {4};
			\end{tikzpicture}
			& & \textbf{H}_3 & \begin{tikzpicture}[scale=0.15, baseline=0]
			\draw[fill= black] (5,0) circle (0.5);
			\draw[line width=1] (5,0) -- (15,0);
			\draw[fill= black] (10,0) circle (0.5);
			\draw[fill= black] (15,0) circle (0.5);
			\draw (7.5,2) node {5};
			\end{tikzpicture}  \\
			\textbf{D}_n & \begin{tikzpicture}[scale=0.15, baseline=0]
			\draw[fill= black] (5,-3) circle (0.5);
			\draw[line width=1] (5,-3) -- (10,0);
			\draw[line width=1] (5,3) -- (10,0);
			\draw[line width=1] (15,0) -- (10,0);
			\draw[fill= black] (5,3) circle (0.5);
			\draw[fill= black] (10,0) circle (0.5);
			\draw[fill= black] (15,0) circle (0.5);
			\draw (17.5,0) node {$\ldots$};				
			\draw[fill= black] (20,0) circle (0.5);
			\end{tikzpicture}
			& & \textbf{H}_4 & \begin{tikzpicture}[scale=0.15, baseline=0]
			\draw[fill= black] (5,0) circle (0.5);
			\draw[line width=1] (5,0) -- (20,0);
			\draw[fill= black] (10,0) circle (0.5);
			\draw[fill= black] (15,0) circle (0.5);
			\draw[fill= black] (20,0) circle (0.5);
			\draw (7.5,2) node {5};
			\end{tikzpicture}  \\
			\textbf{I}_2(p) & \begin{tikzpicture}[scale=0.15, baseline=0]
			\draw[fill= black] (5,0) circle (0.5);
			\draw (7.5,2) node {p};
			\draw[line width=1] (5,0) -- (10,0);			
			\draw[fill= black] (10,0) circle (0.5);
			\end{tikzpicture}
			& &\textbf{E}_6 & \begin{tikzpicture}[scale=0.15, baseline=0]
			\draw[fill= black] (5,0) circle (0.5);
			\draw[line width=1] (5,0) -- (25,0);
			\draw[fill= black] (10,0) circle (0.5);
			\draw[fill= black] (15,0) circle (0.5);
			\draw[fill= black] (20,0) circle (0.5);		
			\draw[fill= black] (25,0) circle (0.5);			
			\draw[fill= black] (15,-5) circle (0.5);
			\draw[line width=1] (15,0) -- (15,-5);
			\end{tikzpicture}  \\
			& & & \textbf{E}_7 & \begin{tikzpicture}[scale=0.15, baseline=0]
			\draw[fill= black] (5,0) circle (0.5);
			\draw[line width=1] (5,0) -- (30,0);
			\draw[fill= black] (10,0) circle (0.5);
			\draw[fill= black] (15,0) circle (0.5);
			\draw[fill= black] (20,0) circle (0.5);		
			\draw[fill= black] (25,0) circle (0.5);		
			\draw[fill= black] (30,0) circle (0.5);		
			\draw[fill= black] (15,-5) circle (0.5);
			\draw[line width=1] (15,0) -- (15,-5);
			\end{tikzpicture}   \\
			& & & \textbf{E}_8 & \begin{tikzpicture}[scale=0.15, baseline=0]
			\draw[fill= black] (5,0) circle (0.5);
			\draw[line width=1] (5,0) -- (35,0);
			\draw[fill= black] (10,0) circle (0.5);
			\draw[fill= black] (15,0) circle (0.5);
			\draw[fill= black] (20,0) circle (0.5);		
			\draw[fill= black] (25,0) circle (0.5);		
			\draw[fill= black] (30,0) circle (0.5);		
			\draw[fill= black] (35,0) circle (0.5);		
			\draw[fill= black] (15,-5) circle (0.5);
			\draw[line width=1] (15,0) -- (15,-5);
			\end{tikzpicture}  \\
	}}
\end{thm}

\begin{defn}
	We say that a finite irreducible Coxeter group $W$ is \emph{of type} $\D$ if its corresponding diagram is given by $\D$, and we denote this Coxeter group $W(\D)$.
\end{defn}

\begin{rem}
	The Coxeter group $W(\textbf{A}_n)$ is isomorphic to~$S_{n+1}$, the symmetric group, which is the reflection group of the regular $(n+1)$-simplex.
\end{rem}

\begin{defn}
	Let $(W,S)$ be a Coxeter system. For each $T\subseteq S$, $T$ generates a subgroup $W_T$ such that $(W_T,T)$ is a Coxeter system in its own right. We call subgroups that arise in this way \emph{parabolic subgroups}. If the subgroup is finite we call it a \emph{spherical subgroup}.
\end{defn}

\subsection{Artin groups}
This section follows Charney \cite[Section 1.]{Charney} and notes by Paris \cite{Paris}.

Given a Coxeter system $(W,S)$, the corresponding Artin group is given by forgetting the involution relations i.e.~setting $m(s,s)=\infty$.

\begin{defn}
	For every Coxeter system $(W,S)$ there is a corresponding \emph{Artin system} $(A_W,\Sigma)$ comprising of an {\emph Artin group} $A_W$ with generating set $$\Sigma:=\{ \sm_s \text{ for } s \in S\}$$
	\noindent and presentation
	$$
	A_W=\langle \,\Sigma \, \vert \,  s, t \in S, \altprod(\sm_s,\sm_t; m(s,t))=\altprod(\sm_t,\sm_s; m(s,t))\rangle.
	$$
\end{defn}
\noindent We note that the Coxeter diagram $\D_W$ contains all the information about the Artin group presentation.

\begin{exam}
	The Artin group $A_W$ corresponding to the Coxeter group $W\cong S_n$ is the \emph{braid group} $\B_n$. The corresponding diagram $\D_W$ is 		
	\begin{center}
		\begin{tikzpicture}[scale=0.2, baseline=0]
		\draw[fill= black] (5,0) circle (0.4) node[below] {$\sm_1$};	
		\draw[line width=1] (5,0) -- (15,0);
		\draw[fill= black] (10,0) circle (0.4) node[below] {$\sm_2$};
		\draw[fill= black] (15,0) circle (0.4) node[below] {$\sm_3$};		
		\draw (17.5,0) node {$\ldots$};
		\draw[fill= black] (20,0) circle (0.4) node[below] 	{$\sm_{n-2}$};
		\draw[line width=1] (20,0) -- (25,0);				
		\draw[fill= black] (25,0) circle (0.4) node[below] {$\sm_{n-1}$};
		\end{tikzpicture}
	\end{center}
	where we relabel $\sm_{s_i}$ to $\sm_i$ for ease of notation. The presentation is therefore given by
	$$
	\B_n=\Big\langle \sm_i \,\text{ for } s_i \in S \,\mid\,\begin{array}{lr}\sm_i\sm_j=\sm_j\sm_i & |i-j|\geq 2\\ \sm_i\sm_{i+1}\sm_i=\sm_{i+1}\sm_i\sm_{i+1} &1\leq i \leq (n-2)\end{array}\Big\rangle,
	$$
	the standard presentation for the braid group on $n$ strands.
\end{exam}

\begin{exam}
	When all possible edges in the Coxeter diagram $\D_W$ are present and labelled with $\infty$ the corresponding Artin group is the free group on $|S|$ generators. The group has presentation
	$$
	A_W=\langle \sm_s \text{ for } s \in S\rangle.
	$$
\end{exam}

\begin{exam}
	When there are no edges in the Coxeter diagram $\D_W$ the corresponding Artin group is the free abelian group on $|S|$ generators. The group has presentation
	$$
	A_W=\langle \sm_s \text{ for } s \in S \,\vert\, \sm_s\sm_t=\sm_t\sm_s \, s\neq t \in S \rangle.
	$$
\end{exam}

\begin{exam}
	When all of the edges in the Coxeter diagram are labelled with~$\infty$, but not necessarily all possible edges are present (some $m(s,t)$ may be equal to $2$) then the corresponding Artin group is called a \emph{right angled} Artin group, or RAAG.
\end{exam}

\begin{defn}
	When the Coxeter group $W$ is finite, i.e.~when its diagram $\D_W$ is a disjoint union of diagrams from Proposition \ref{prop:classification of finite Coxeter}, then the corresponding Artin group $A_W$ is called a \emph{finite type} Artin group, or a \emph{spherical} Artin group.
\end{defn}

Much of the known theory of Artin groups is concentrated around RAAGs and finite type Artin groups, though we do not restrict ourselves to either of these families in our results. In general little is known about Artin groups. For instance the following properties hold for finite type Artin groups \cite{Charney}:

\begin{itemize}
	\item there exists a finite model for the classifying space $K(A_W,1)$,
	\item $A_W$ is torsion free,
	\item the centre of $A_W$ is $\Z$, for $A$ irreducible,
	\item $A_W$ has solvable word and conjugacy problem.
\end{itemize}

\noindent To date these properties are not known for general Artin groups. In the next section we consider the first point in detail.

\section{The $K(\pi,1)$ conjecture}\label{sec:K(pi,1)conjecture}

This section introduces the $K(\pi,1)$ conjecture, following \cite{Paris}.

In general, one can associate a hyperplane arrangement $\mathcal{A}$ and associated complement $M(\mathcal{A})$ to each Coxeter group $W$, such that there is a free action of $W$ on~$M(\mathcal{A})$. When we consider this hyperplane complement modulo this $W$ action, the corresponding quotient $M(\mathcal{A})/W$ has as its fundamental group the Artin group~$A_W$. In some cases this quotient space is known to be a $K(A_W,1)$, in particular we recall Deligne's theorem for finite type Artin groups.

\begin{thm}[{Deligne's Theorem, \cite{Deligne}}]\label{thm:delignes}
	For $W$ a finite Coxeter group and~$A_W$ the associated Artin group, $M(\mathcal{A})/W$ is aspherical with fundamental group~$A_W$, that is $M(\mathcal{A})/W$ is a $K(A_W,1)$.
\end{thm}

For arbitrary Artin groups, the \emph{$K(\pi,1)$ conjecture} was formulated by Arnol'd, Brieskorn, Pham and Thom, and states than an analogue of Deligne's theorem holds for all Artin groups. The analogue of the hyperplane complement was formulated by Vinberg. For a more detailed description see Davis \cite{Davis}, notes by Paris \cite{Paris} and the introduction to a paper on RAAGs by Charney \cite{Charney}.

\begin{rem}
	It is worth noting here a reformulation of the conjecture in terms of a finite dimensional CW-complex called the \emph{Salvetti complex}, denoted by $Sal(\mathcal{A})$ and introduced by Salvetti in \cite{Salvetti}, for a hyperplane arrangement $\mathcal{A}$ in a finite dimensional real vector space $V$. The Salvetti complex is defined in terms of cosets of finite subgroups of the Coxeter group~\cite{Salvetti94}. Paris extended this definition to any infinite hyperplane arrangement in a nonempty convex cone $I$ \cite{Paris} and proved that $Sal(\mathcal{A})$ and $M(\mathcal{A})$ have the same homotopy type. The $K(\pi,1)$ conjecture can therefore be restated as a conjecture about the Salvetti complex.
\end{rem}

The $K(\pi,1)$ conjecture has been proven for large classes of Artin groups \cite{Paris}. However the conjecture has not been proven for general Artin groups. We will apply a reformulation of the $K(\pi,1)$ conjecture to our results, involving the Artin monoid~$A^+$ and discussed in Section \ref{chap:BG Artin monoids}.

\section{Artin monoids}\label{chap:BG Artin monoids}

The start of this section follows Jean Michel \emph{A note on words in braid monoids} \cite{Michel} and Brieskorn and Saito \emph{Artin-Gruppen und Coxeter-Gruppen} \cite{BreiskornSaito}. Much of the material in Section~\ref{sec:generic artin monoid divisors and cosets 2} is new.

\subsection{Definition and examples}

\begin{defn}
	The \emph{Artin monoid system} $(A_W^+, \Sigma)$ associated to a Coxeter system~$(W,S)$ is given by the generating set $\Sigma$ for the corresponding Artin system~$(A_W,\Sigma)$, and the monoid with the same presentation as the Artin group $A_W$:
	$$
	A_W^+=\langle \Sigma \, \vert \, \altprod(\sm_s,\sm_t; m(s,t))=\altprod(\sm_t,\sm_s; m(s,t))\rangle^+.
	$$
	\noindent Words in $A^+_W$ are therefore strings of letters for which the alphabet consists of $\sm_s$ in $\Sigma$.
	
\end{defn}

\begin{rem}
	The group completion of $A_W^+$ is $A_W$.
\end{rem}

\begin{exam}
	The braid monoid $\mathcal{B}_n^+$ is the monoid associated to the Coxeter group~$S_n$, the symmetric group, with group completion the braid group~$\mathcal{B}_n$. Given the standard generating set for the symmetric group, the braid monoid consists of words in the braid group made from the positive generators $\sm_i$. In terms of the braid diagrams these can be viewed as braids consisting of only positive twists.
\end{exam}

\begin{defn}
	We call a submonoid $M^+$ of an Artin monoid $A^+$ a \emph{parabolic submonoid} if the monoid $M^+$ is generated by the set $M^+\cap \Sigma$. We call this generating set for the monoid $\Sigma_M$, giving a system~$(M^+,\Sigma_M)$.
\end{defn}

In this paper, by convention every submonoid of an Artin monoid considered will be a parabolic submonoid.

\subsection{Divisors in Artin monoids}\label{sec:generic artin monoid divisors and cosets}

\begin{defn}	
	Define the \emph{length function }on an Artin monoid $A^+$ with system~$(A^+,\Sigma)$
	$$\lng: A^+ \to \mathbb{N}$$
	to be the function which maps an element $\alpha$ in $A^+$ to the unique word length required to express $\alpha$ in terms of the generators in $\Sigma$.
\end{defn}

\begin{rem}
	Note here that since there are no inverses in Artin monoids, multiplication corresponds to addition of lengths, i.e.~$\lng(ab)=\lng(a)+\lng(b)$ ($\lng$ is a monoid homomorphism).
\end{rem}

\begin{defn}
	For elements $\alpha$ and $\beta$ in an Artin monoid $A^+$ with system $(A^+,\Sigma)$, we say that $\alpha \preceq_R \beta$ if for some $\gamma$ in $A^+$ we have $\beta=\gamma \alpha$, that is a word representing~$\alpha$ appears on the right of some word representing  $\beta$, in terms of the generating set~$\Sigma$. We say that $\beta$ is \emph{right-divisible} by $\alpha$, or alternatively that $\alpha$ \emph{right divides} $\beta$.
\end{defn}

\begin{prop}[{Michel \cite[Prop 2.4]{Michel}}] \label{prop:cancellation}
	Artin monoids satisfy left and right cancellation, i.e. for $a$, $b$ and $c$ in $A^+$,
	\begin{eqnarray*}
		ab=ac \Rightarrow b=c\\
		ba=ca \Rightarrow b=c.
	\end{eqnarray*}
\end{prop}

We now consider work by Brieskorn and Saito in their 1972 paper \textit{Artin-Gruppen und Coxeter-Gruppen} \cite{BreiskornSaito}. They consider notions of least common multiples and greatest common divisors of sets of elements in the Artin monoid. We are interested in the notion of least common multiple.

\begin{defn}\label{def:lcm}
	Given a set of elements $\{g_j\}_{j\in J}$ in an Artin monoid $A^+$ with system~$(A_W,\Sigma)$, a \emph{common multiple} $\beta$ is an element in $A^+$ which is right divisible by all $g_j$. That is $g_j \preceq_R \beta$ for all $j$ in $J$. A \emph{least common multiple} of $\{g_j\}$ is a common multiple that right-divides all other common multiples.
		Let $E$ be a set of elements in the Artin monoid $A^+$. Denote the least common multiple (if it exists) of $E$ by $\Dl(E)$. For $\alpha$ and $\beta$ two elements in $A^+$ denote the least common multiple of $\alpha$ and $\beta$ (if it exists) by $\Dl(\alpha, \beta)$.
\end{defn}

\begin{rem}
	Should a least common multiple exist, it will be unique.
\end{rem}

\begin{prop}[{Brieskorn and Saito \cite[4.1]{BreiskornSaito} }]\label{prop:lcm}
	A finite set of elements in an Artin monoid either has a least common multiple or no common multiple at all.
\end{prop}

\begin{rem}\label{rem:set of letters in element well defined}
	Since the relations in an Artin monoid have the same letters appearing on each side, the set of letters present in any word representing an element of an Artin monoid is fixed. Therefore the notion of `letters appearing in an element' is well defined.
\end{rem}

\begin{lem}[{Brieskorn and Saito \cite{BreiskornSaito}}] \label{rem:lcm}
	For a fixed generating set, the letters arising in a least common multiple of a set of elements in an Artin monoid are only those letters which appear in the elements themselves.
\end{lem}

\begin{defn}\label{def:end sets}
	Consider a submonoid $M^+$ of an Artin monoid $A^+$, with system~$(M^+,\Sigma_M)$. Given an element $\alpha$ in $A^+$ we define two \emph{end sets}
	\begin{eqnarray*}
		\EndGen_M (\alpha)&=&\{\sm_s \mid s \in \Sigma_M,  \sm_s \preceq_R \alpha\}\\
		\EndMon_{M} (\alpha)&=&\{\beta \in M^+ \mid \beta \preceq_R \alpha\}.
	\end{eqnarray*}
\end{defn}

\begin{rem}
	$\EndGen_M(\alpha)$ is exactly the letters $\sm_s$ in~$\Sigma_M$ that a word representing $\alpha$ can end with, and $\EndMon_M(\alpha)$ is exactly the elements in $M^+$ that a word representing $\alpha$ can end with. Note that $\EndGen_M(\alpha)$ is a subset of $\EndMon_M(\alpha)$, consisting of words of length $1$ and $\EndMon_M(\alpha)=\emptyset$ if and only if $\alpha$ has no right-divisors in $M^+$.
\end{rem}

\subsection{Required theory}\label{sec:generic artin monoid divisors and cosets 2}

Much of the proof of Theorem \ref{thm:A} is concerned with algebraic manipulation of words in the Artin monoid. Here we introduce some technical definitions and lemmas used in the proof. We build up a theory of cosets in the case of Artin monoids, which is new unless cited.

\begin{lem}
	Given $\alpha$ in $A^+$, and $M^+$ a submonoid of $A^+$, the set $\EndMon_M(\alpha)$ has a least common multiple $\Dl(\EndMon_M(\alpha))=\beta$ which lies in the submonoid~$M^+$. That is, there exists $\beta$ in $M^+$ and $\gamma$ in $A^+$ such that $\alpha=\gamma \beta$ for some words representing~$\alpha, \beta, \gamma$, and if $\beta'$ in $A^+$ and $\gamma'$ in $A^+$ satisfy $\alpha=\gamma'\beta'$, it follows that~$\beta \preceq_R \beta'$.
	\begin{proof}
		From Proposition \ref{prop:lcm} if a common multiple exists, then~$\Dl(\EndMon_M(\alpha))$ exists. We have that $\alpha$ itself is a common multiple of all elements in $\EndMon_M(\alpha)$, by definition of $\EndMon_M(\alpha)$. Furthermore Lemma \ref{rem:lcm} notes that only letters appearing in $\EndMon_M(\alpha)$ will appear in $\Dl(\EndMon_M(\alpha))$. By definition these are letters in $M^+$ and so $\Dl(\EndMon_M(\alpha))$ lies in $M^+$.
	\end{proof}
\end{lem}

\begin{rem}\label{rem:balpha def}
	For n element~$\alpha$ in $A^+$ let $\Dl(\EndMon_M(\alpha)) = \beta$. We write $\balpha$ \emph{with respect to $M^+$} for the element $\balpha$ in $A^+$ such that $\alpha=\balpha \beta$. It will always be clear in the text for which submonoid $M^+$ we are taking the reduction with respect to.
\end{rem}

\begin{defn} \label{def:cosets}
	For $A^+$ an Artin monoid and $M^+$ a submonoid, let $A^+(M)$ be the following set
	$$
	A^+(M)=\{\balpha \text{ with respect to }M^+ \mid \alpha \in A^+\}.
	$$
	That is, $ A^+(M)$ is the set of elements in $A^+$ which are not right divisible by any element of~$M$.
\end{defn}

\begin{lem} \label{lemma:barmult}
	For all $\alpha$ in $A^+$ and all $\beta$ in $M^+$, $\balpha=\overline{\alpha\beta}$ where the reduction is taken with respect to $M^+$.
	\begin{proof}
		Let $\balpha=\gamma$, so $\alpha=\gamma \eta $ for some $\eta$ in $M^+$, and $\EndMon_M(\gamma)=\emptyset$ i.e.~$\gamma$ has no right divisors in $M^+$.
		Then $\alpha \beta = \gamma \eta \beta$ and since $\eta$ and $\beta$ are both in~$M^+$, it follows that $\eta \beta \in \EndMon_M(\alpha \beta)$. If $\eta \beta$ is  the least common multiple of $\EndMon_M(\alpha \beta)$ then~$\overline{\alpha\beta}=\gamma =\balpha$ so we are done. Suppose for a contradiction that~$\eta \beta$ is not the least common multiple of $\EndMon_M(\alpha \beta)$, and note that~$\eta \beta$ is a right divisor of the least common multiple. Then there exists some $\zeta$ in $M^+$ of length at least $1$ such that $\zeta \eta \beta$ is the least common multiple of $\EndMon_M(\alpha \beta)$. It follows that there exists a $\gamma'=\overline{\alpha\beta}$ with $\EndMon_M(\gamma')=\emptyset$ and~$\alpha \beta = \gamma' \zeta \eta \beta$. But~$\alpha \beta = \gamma \eta \beta$ and it follows from cancellation that $\gamma=\gamma' \zeta$. Since $\zeta$ is in $M^+$ with length at least $1$ it follows that $\zeta \in \EndMon_M(\gamma)$ which contradicts $\EndMon_M(\gamma)=\emptyset$. Therefore $\eta \beta$ is  the least common multiple of $\EndMon_M(\alpha \beta)$ and it follows that $\overline{\alpha\beta}=\gamma =\balpha$.
	\end{proof}
\end{lem}

\begin{defn}\label{def:monoid coset}
	Consider the relation $\sim$ on $A^+$ given by
	$$
	\alpha_1 \sim \alpha_2 \iff \alpha_1 \beta_1=\alpha_2 \beta_2 \text{ for some }\beta_1\text{ and }\beta_2\text{ in }M^+
	$$
	\noindent where $M^+$ is a submonoid of $A^+$.~When we use this relation it will be made clear which submonoid $M^+$ is being considered. The relation~$\sim$ is symmetric and reflexive. Let $\approx$ be the transitive closure of $\sim$. That is, $\alpha_1 \approx \alpha_2$ if there is a chain of elements in $A^+$:
	$$
	\alpha_1 \sim \tau_1 \sim \tau_2 \sim \cdots \sim \tau_k \sim \alpha_2
	$$
	\noindent for some $k$. Denote the equivalence class of $\alpha$ in $A^+$ under the relation $\approx$ with respect to the submonoid $M^+$ as $[\alpha]_M$.
\end{defn}

\begin{defn}\label{def:quotient map}
	Let $q:A^+(M)\to A^+/\approx$ be the quotient map, taken with respect to the equivalence relation $\approx$.	
\end{defn}

\begin{lem} \label{lem:barcosetequality}
	The map $q$ of Definition \ref{def:quotient map} is a bijection. That is for all $\alpha_1$ and~$\alpha_2$ in $A^+$:
	\begin{equation*}
	[\alpha_1]_M=[\alpha_2]_M \iff \overline{\alpha_1}=\overline{\alpha_2}
	\end{equation*}
	\begin{proof}
		($\Leftarrow$) If $\overline{\alpha_1} =\overline{\alpha_2}=\gamma$ with respect to $M^+$, then $\alpha_1 \sim \gamma \sim \alpha_2$ so it follows~$\alpha_1~\approx~\alpha_2$.\\
		($\Rightarrow$) We want to show that if $\alpha_1 \approx \alpha_2$ then $\overline{\alpha_1} =\overline{\alpha_2}$. Since $\alpha_1 \approx \alpha_2$ there is a chain $$\alpha_1 \sim \tau_1 \sim \tau_2 \sim \cdots \sim \tau_k \sim \alpha_2$$ so if we show that $\overline{\eta} =\overline{\zeta}$ whenever $\eta \sim \zeta$ for $\eta$ and $\zeta$ in $A^+$ it will follow that $$\balpha_1 = \btau_1 = \btau_2 = \cdots = \btau_k = \balpha_2.$$ Since $\eta \sim \zeta$ it follows that for some $\beta_1$ and $\beta_2$ in $M^+$, $\eta\beta_1=\zeta \beta_2$. Applying Lemma~\ref{lemma:barmult} gives
		$$
		\overline{\eta}=\overline{\eta\beta_1}=\overline{\zeta \beta_2}=\overline{\zeta}
		$$
		as required.
	\end{proof}
\end{lem}

\begin{prop}\label{prop:decomp} For $M^+$ a submonoid of $A^+$,
	$A^+ \cong A^+(M) \times M^+$ as sets, via the bijection
	\begin{eqnarray*}
		p:A^+ &\to& A^+(M) \times M^+ \\
		\alpha &\mapsto& (\balpha, \beta) \text{ where } \alpha=\balpha \beta
	\end{eqnarray*}
	\noindent where~$\beta=\Delta(EndMon_M(\alpha))$. This decomposition respects the right action of $M^+$ on $A^+$, i.e.~$M^+$ acts trivially on the first factor and as right multiplication on the second.
	\begin{proof}
		To show $p$ is surjective: consider $(\gamma, \beta) \in A^+(M) \times M^+$. Due to Lemma~\ref{lemma:barmult} for~$\alpha \in A^+$ ans any $\beta \in M^+$ we have $\overline{\alpha\beta}=\balpha$. Therefore $\gamma\beta$ satisfies $p(\gamma\beta)= (\gamma, \beta)$ since $\overline{\gamma\beta}=\overline{\gamma}=\gamma$ (we have $\gamma \in A^+(M)$ so $\EndMon_p(\gamma)=\emptyset$).		
		To show injectivity, suppose $p(\alpha_1)=p(\alpha_2)$, that is $(\overline{\alpha_1}, \beta_1)=(\overline{\alpha_2}, \beta_2)$. This translates to
		$$\alpha_1=\overline{\alpha_1}\beta_1=\overline{\alpha_2}\beta_2=\alpha_2$$ 
		therefore $p$ is injective.
		Under this decomposition, the action of $m$ in $M^+$ satisfies~$p(\alpha\cdot m)=(\balpha, \beta\cdot m)$ where $\alpha=\balpha \beta$, again due to Lemma \ref{lemma:barmult}.
	\end{proof}
\end{prop}

\begin{prop}[{\cite[1.5]{Michel}}] \label{prop:deltagens}
	If generators $s$ and $t$ in $S_M$ are in $\EndGen_M(\alpha)$ for some $\alpha$ in $A^+$ then $\Dl(s,t)$ lies in $\EndMon_M(\alpha)$.
\end{prop}

\begin{lem}\label{lem:deltaendsets}
	Consider a subset $F$ of $\EndMon_M(\alpha)$ for some submonoid $M^+$ of~$A^+$ and some $\alpha$ in~$A^+$. Then $\Dl(F)$ is in $\EndMon_M(\alpha)$.
	\begin{proof}
		Since $F$ is a subset of $\EndMon_M(\alpha)$, which has a least common multiple, $F$ has a common multiple so $\Dl(F)$ exists by Proposition \ref{prop:lcm}. Since $\Dl(\EndMon_M(\alpha))$ is a common multiple for $\EndMon_M(\alpha)$, it is a common multiple for $F$. The element $\Dl(F)$ right-divides all other common multiples of $F$ by definition. Therefore~$\Dl(F)\preceq_R\Dl(\EndMon_M(\alpha))$ and it follows that $\Dl(F)$ is in $\EndMon_M(\alpha)$.
	\end{proof}
\end{lem}

\begin{defn} \label{def:LW}
	Recall from Remark~\ref{rem:set of letters in element well defined} that the set of letters present in an element of an Artin monoid is well defined. We say elements $\alpha$ and $\beta$ in an Artin monoid with system $(A^+,\Sigma)$ \emph{letterwise commute} if:
	\begin{itemize}
		\item the set of letters in $\Sigma$ that $\alpha$ contains is disjoint from the set of letters that~$\beta$ contains, and
		\item each letter of $\Sigma$ that $\alpha$ contains commutes with each letter of $\Sigma$ that~$\beta$ contains.
	\end{itemize}
\end{defn}

\begin{lem}\label{lem:LWCends}
	If $\beta$ and $\gamma$ are in $\EndMon_M(\alpha)$  and $\beta$ and $\gamma$ letterwise commute, it follows that:
	\begin{itemize}
		\item $\Dl(\beta,\gamma)=\beta\gamma=\gamma\beta$
		\item $\Dl(\beta, \gamma) \in \EndMon_M(\alpha)$.
	\end{itemize}
	\begin{proof}
		Since $\beta$ and $\gamma$ letterwise commute, they contain distinct generators. From Remark~\ref{rem:set of letters in element well defined} every word representative for $\beta$ and $\gamma$ contains the same set of letters. It follows each of these letters must appear in $\Dl(\beta, \gamma)$. If both~$\beta$ and $\gamma$ have length~$1$, say $\beta=\sigma$ and $\gamma=\tau$ for generators $\sm$ and $\tau$ then since the words letterwise commute it follows that $\sigma$ commutes with $\tau$. Therefore since~$\sm\tau=\tau\sm$ and both generators must appear in $\Dl(\beta, \gamma)$ it follows that
		$$
		\Dl(\beta, \gamma)=\sm\tau=\tau\sm=\beta\gamma=\gamma\beta
		$$
		\noindent as required. Similarly, if $\beta=\sm_1\ldots \sm_k$ has length $k$, and $\gamma=\tau$ has length $1$ then since the words contain distinct generators (which all must appear in the lowest common multiple) and~$\tau\sm_i=\sm_i\tau\, \forall \,i$ it follows that
		$$
		\Dl(\beta, \tau)=\Dl(\sm_1\ldots \sm_{k},\tau)=(\sm_1\ldots \sm_{k})\tau=\tau(\sm_1\ldots \sm_{k})=\beta \tau=\tau\beta.
		$$
		Suppose now that $\beta=\sm_1\ldots \sm_k$ has length $k$ and $\gamma=\tau_1\ldots \tau_l$ has length $l$. It is true that $\beta \preceq_R \beta \gamma $ and $\gamma \preceq_R \beta \gamma$. \\
		{\bf Claim:}
		If $x$ in $A^+$ is a common multiple of $\beta$ and $\gamma$ then $\beta\gamma=\gamma\beta$ is in $\EndMon_M(x)$.\\
		{\bf Proof of claim:} we proceed by induction on~$\lg(\gamma)$. The base case~$\lg(\gamma)=1$ is covered above. As our inductive hypothesis, we suppose if~$\lg(\gamma')<l$ and~$\gamma'$ letterwise commutes with~$\beta$, then~$\Dl(\gamma',\beta)=\gamma'\beta=\beta\gamma'$, hence the claim holds. We prove the claim for 
		$$
		\beta=\sm_1\ldots\sm_k \text{  and  } \gamma=\tau_1\ldots \tau_l.
		$$
		Since~$x$ is a common multiple of~$\beta$ and~$\gamma$, there exist~$y$ and~$z$ in~$A^+$ such that 
		$$
		x=y\beta \text{  and  } x=z\gamma=z\tau_1\ldots \tau_l.
		$$
		In particular,~$\beta$ and~$\gamma'=\tau_2\ldots \tau_l$ are in~$\EndMon_{M}(x)$, so $$\Dl(\gamma',\beta)=\gamma'\beta=\beta\gamma' \in \EndMon_{M}(x)$$ 
		by the inductive hypothesis. Therefore, there exists~$w\in A^+$ such that
		$$ x=z\gamma=z\tau_1\ldots \tau_l\text{  and  }x=w\beta\gamma'=w\beta\tau_2\ldots \tau_l$$
		which by cancellation of $\tau_2\ldots \tau_l$ gives
		$ z\tau_1=w\beta.$ By the base case, since~$\tau_1$ and~$\beta$ are both in~$\EndMon_{M}(z\tau_1)$ so is~$\Dl(\tau_1,\beta)=\tau_1\beta=\beta\tau_1$. Therefore there exists~$v\in A^+$ such that
		$$
		z\tau_1=v\beta\tau_1
		$$
		and so by cancellation of~$\tau_1$, it follows that~$z=v\beta$. Reinserting this in the previous equation gives
		$$
		x=z\gamma=v\beta \gamma
		$$
		and so~$\beta \gamma=\gamma\beta$ is in~$\EndMon_{M}(x)$ as required.
	\end{proof}
\end{lem}

\begin{lem}\label{lem:mini}
	If words $\alpha$, $a$ and $b$ in $A^+$ are such that $b \preceq_R \alpha a$ and $a$ and $b$ letterwise commute then it follows that $b \preceq_R\alpha$.
	\begin{proof}
		An equivalent way of writing $m\preceq_R n$ for $m$, $n$ in $A^+$ is $m \in \EndMon_A(n)$ where the end set is taken with respect to the full monoid $A^+$. Since $a$ and $b$ are both in $\EndMon_A(\alpha a)$ it follows that $\Dl(a,b)$ is in $\EndMon_A(\alpha a)$, from Lemma~\ref{lem:deltaendsets}. Since $a$ and $b$ letterwise commute, $\Dl(a,b)=ab=ba$ from Lemma~\ref{lem:LWCends}. Therefore~$ba$ is in~$\EndMon_A(\alpha a)$ and, by cancellation of $a$, $b$ is in $\EndMon_A(\alpha)$ as required.
	\end{proof}
\end{lem}

\subsection{Relation to the $K(\pi,1)$ conjecture}

In 2002 Dobrinskaya published a paper relating the classifying space of the Artin monoid $B A_W^+$ to the $K(\pi, 1)$ conjecture. This was later translated into English as \emph{Configuration Spaces of Labelled Particles and Finite Eilenberg - MacLane Complexes} \cite{Dobrinskaya2006}. The main result of the paper was the following:

\begin{thm}[{Dobrinskaya \cite[Theorem 6.3]{Dobrinskaya2006}}]\label{thm:monoid relation to kpi1 dobrinskaya}
	Given an Artin group $A_W$ and its associated monoid $A_W^+$, the $K(\pi,1)$ conjecture holds if and only if the natural map between their classifying spaces, $BA_W^+ \to BA_W$ is a homotopy equivalence.
\end{thm}

This theorem has been reproved using a different Morse-theoretic approach by Ozornova \cite{ozornova2017discrete} and her result has in turn been strengthened by Paolini \cite{paolini2017classifying}.

\section{Semi-simplicial constructions with monoids}\label{sec: semi simp constructions with monoids}

This section is split into three subsections. The first introduces background semi-simplicial theory before the second introduces theory for generic monoids and submonoids, including some new results. The third subsection focuses on Artin monoids and contains results required later in the proof.

\subsection{Semi-simplicial objects}

This subsection consists of the required background and follows Ebert and Randal-Williams \cite{EberyRandalWilliams}.

\begin{defn}[{\cite[1.1]{EberyRandalWilliams}}]
	Let $\Delta$ denote the category which has as objects the non-empty finite ordered sets  $[n]=\{0,1,\ldots,n\}$, and as morphisms monotone increasing functions. These functions are generated by the basic functions which act on the ordered sets as follows:
	\begin{eqnarray*}
		D^i:[n]&\to& [n+1] \text{   for } 0\leq i \leq n\\
		\{0,1,\ldots,n\}&\mapsto& \{0,1,\ldots, \widehat{i},\ldots,n+1\}\\
		S^i:[n+1]&\to& [n]\text{   for } 0\leq i \leq n\\
		\{0,1,\ldots,n+1\} &\mapsto& \{0,1,\ldots,i,i, \ldots n\}
	\end{eqnarray*}
	The opposite category $\Delta^{op}$ is known as the simplicial category. We denote the opposite of the maps $D^i$ by $\cd_i$ and the opposite of the maps $S^i$ by $s_i$. We call these the \emph{face maps} and the \emph{degeneracy maps} respectively.
	
	Let $\Delta_{inj}\subset \Delta$ be the subcategory of $\Delta$ which has the same objects but only the injective monotone maps as morphisms, generated by the $D_i$. The opposite category $\Delta_{inj}^{op}$ is known as the semi-simplicial category and its morphisms are therefore generated by the face maps $\cd_i$.
\end{defn}

\begin{defn} [{\cite[1.1]{EberyRandalWilliams}}]
	A \emph{simplicial object} in a category $\mathcal{C}$ is a covariant functor $X_{\bullet}: \Delta^{op} \rightarrow \mathcal{C}$. A \emph{semi-simplicial object} is a functor $X_{\bullet}: \Delta_{inj}^{op} \rightarrow \mathcal{C}$. We denote $X_\bt([n])$ by $X_n$. A \emph{(semi-)simplicial map} $f: X_{\bullet} \rightarrow Y_{\bullet}$ is a natural transformation of functors, and in particular has components $f_n: X_n \rightarrow Y_n$. Simplicial objects in $\mathcal{C}$ form a category denoted $\bf{s}\mathcal{C}$, and semi-simplicial objects a category denoted $\bf{ss}\mathcal{C}$. When $\mathcal{C}$ equals \textit{\bf{Set}} a (semi-)simplicial object is called a \emph{(semi-)simplicial set} and when $\mathcal{C}$ equals \textit{\bf{Top}} it is called a \emph{(semi-)simplicial space}.
\end{defn}

\begin{rem}
	A semi-simplicial object in a category $\mathcal{C}$ is equivalent to the following data:
	\begin{enumerate}[(a)]
		\item An object $X_p$ in $\mathcal{C}$, for $p\geq 0$
		\item Morphisms in $\mathcal{C}$ $\cd^p_i:X_p \to X_{p-1}$ for $0 \leq i \leq p$ and all $p\geq 0$ called \emph{face maps}, which satisfy the following \emph{simplicial identities}
		$$
		\cd^{p-1}_i\cd^p_j=\cd^{p-1}_{j-1}\cd^p_i \text{ if } i <j.
		$$
	\end{enumerate}
\end{rem}

\begin{defn}[{\cite[1.3]{EberyRandalWilliams}}]\label{def:augmentation}
	An augmented semi-simplicial object in $\mathcal{C}$ is a triple~($X_\bt, X_{-1}, \epsilon_\bt)$ such that $X_\bt$ is a semi-simplicial object in $\mathcal{C}$, $X_{-1}$ is an object of $\mathcal{C}$ and $\epsilon_\bt$ is a family of morphisms such that $\epsilon_p:X_p\to X_{-1}$ and $\epsilon_{p-1}\circ \cd_i=\epsilon_{p}$ for all~$p\geq 1$ and $0\leq i \leq p$.
\end{defn}

\begin{exam}[{\cite[1.2]{EberyRandalWilliams}}]
	The standard $n$-simplex has two equivalent manifestations: as a simplicial object in $\bf{Set}$ and as an object in $\bf{Top}$. When viewed as a simplicial set the standard $n$-simplex is denoted $\Dl^n_\bt$ and is defined via the functor $\Dl^n_m=\Dl^n_\bt([m])=\hom_{\Dl}([m],[n])$ for all $[m]$ in $\Dl^{op}$. When viewed as an object in $\bf{Top}$ the standard $n$-simplex is denoted $|\Dl^n|$ and defined to be
	$$
	|\Dl^n|:=\Big\{ (t_0, \ldots , t_n) \in \mathbb{R}^{n+1} \, \vert \, \sum_{i=0}^n t_i =1 \text{ and } t_i \geq 0 \forall i \Big\}.
	$$
	\noindent One can associate to a morphism $\phi: [n]\to [m]$ in $\Dl$ a continuous map
	\begin{eqnarray*}
		\phi_*: |\Dl^n| &\to& |\Dl^m|\\
		(t_0,\ldots, t_n) &\mapsto& (s_0,\ldots, s_m) \text{   where  } s_j=\sum_{\phi(i)=j}t_i.
	\end{eqnarray*}
	That is, morphisms send the $j$th vertex of the simplex~$|\Dl^n|$ to the $\phi(j)$th vertex of~$|\Dl^m|$ and extend linearly. Under this viewpoint the map $D^i_*$ sends~$|\Dl^n|$ to the $i$th face of~$|\Dl^{n+1}|$ and the map $S^i_*$ collapses together the $i$th and $(i+1)$st vertices of~$|\Dl^{n+1}|$ to give a map to $|\Dl^{n}|$.
\end{exam}	

A tuple $(\cd^{p-1}_{i_1}, \cd^{p-2}_{i_2},\ldots, \cd^{p-k}_{i_k})$ denotes the application of several face maps in a row, where $\cd^{p-1}_{i_1}$ is the first face map to be applied, followed by $\cd^{p-2}_{i_2}$, etc. For ease of notation we dispense with superscripts, writing the tuple as $(\cd_{i_1}, \cd_{i_2},\ldots, \cd_{i_k})$ and assuming the domain and targets are such that the composite map is defined.

\begin{lem}\label{lem:ijs background}
	With the above notation, the tuple of face maps can be written such that $i_{j+1}\geq i_j$ for all $j$.
	\begin{proof}
		Suppose  $i_{j+1}< i_j$ in the tuple $(\cd_{i_1}, \cd_{i_2},\ldots, \cd_{i_k})$. The simplicial identities show 
		$
		\cd_{i_{j+1}}\cd_{i_{j}}=\cd_{i_{j}-1}\cd_{i_{j+1}} \text{ since } i_{j+1} <i_j
		$.
		Therefore $$(\cd_{i_1}, \cd_{i_2}\ldots ,\cd_{i_j},\cd_{i_{j+1}},\ldots,\cd_{i_{k}})=(\cd_{i_1}, \cd_{i_2}\ldots ,\cd_{i_{j+1}},\cd_{i_j-1},\ldots,\cd_{i_{k}}).$$ Since~$i_{j+1}< i_j$, it follows that  $i_{j}-1\geq i_{j+1}$. Relabelling $i_j:=i_{j+1}$ and $i_{j+1}:=i_j-1$ gives $(\cd_{i_1}, \cd_{i_2}\ldots ,\cd_{i_j},\cd_{i_{j+1}},\ldots,\cd_{i_{k}})$ such that $i_{j+1}\geq i_j$. This procedure reduces the sum~$\sum_{j=1}^{k} i_j$ by one, and therefore upon iteration must terminate. Applying this process enough times gives $i_{j+1}\geq i_j$ for all $j$.
	\end{proof}
\end{lem}

\begin{defn}[{\cite[1.2]{EberyRandalWilliams}}]
	The \emph{geometric realisation} of a semi-simplicial set or space is the topological space denoted by $\|X_{\bt}\|$ and defined to be
	\begin{equation*}
	\|X_\bt \|:=\coprod_{n\geq 0} X_n \times |\Delta^n|/\sim
	\end{equation*}
	\noindent where $\sim$ is generated by $(x,t) \sim (y,u)$ whenever $\cd_i(x)=y$ and $D^i(u)=t$.
\end{defn}

The geometric realisation is an example of a coequaliser or colimit, see Dugger~\cite{dugger2008primer}.

\begin{defn}
	Given a semi-simplicial map $f_\bt:X_\bt\to Y_\bt$ there is an induced map~$\|f_\bt\|:\|X_\bt\|\to \|Y_\bt\|$ which we call the \emph{geometric realisation of the semi-simplicial map $f_\bt$}.
\end{defn}

\begin{defn}[{\cite[1.4]{EberyRandalWilliams}}]\label{def:bisemi}
	A \emph{bi-semi-simplicial} object in a category $\mathcal{C}$ is a functor $X_{\bullet \bt}: (\Delta_{inj} \times \Delta_{inj})^{op} \rightarrow \mathcal{C}$. We write $X_{p,q}=X_{\bullet \bt}([p],[q])$. We write the image of the standard face maps in each simplicial direction $(\cd_i\times \id)$ and $(\id\times \cd_j)$, as~$\cd_{i, \bt}$ and $\cd_{\bt, j}$. We note that $$(\cd_i\times \cd_j)= (\cd_{i, \bt}\circ\cd_{\bt, j})=(\cd_{\bt,j}\circ\cd_{i,\bt}):X_{p,q}\to X_{(p-1),(q-1)}$$ and we denote this map $\cd_{i,j}$.
	When $\mathcal{C}$ is equal to \textit{\bf{Top}} the bi-semi-simplicial object is called a \emph{bi-semi-simplicial space}.
\end{defn}

\begin{rem}
	A bi-semi-simplicial space can be viewed as a semi-simplicial object in $\bf{ssTop}$ in two ways:
	\begin{enumerate}[1.]
		\item $X_{\bt,q}:[p]\mapsto (X_\bt:[q]\mapsto X_{p,q} )$ with face maps $\cd_{i, \bt}$.
		\item $X_{p,\bt}:[q]\mapsto (X_\bt:[p]\mapsto X_{p,q} )$ with face maps $\cd_{\bt, j}$.
	\end{enumerate}
\end{rem}

\begin{defn}[{\cite[1.2]{EberyRandalWilliams}}]\label{def:bisemi geom rel}
	Given a bi-semi-simplicial space $X_{\bt, \bt}$ we define its \emph{geometric realisation} to be
	\begin{equation*}
	\|X_{\bt,\bt} \|=\coprod_{p,q\geq 0} X_{p,q} \times |\Delta^p| \times |\Dl^q|/\sim
	\end{equation*}
	where $\sim$ is generated by $(x,t_1,t_2) \sim (y,u_1,u_2)$ whenever $(\cd_{i,j})(x)=y$, $D^i(u_1)=t_1$ and $D^j(u_2)=t_2$. This is equivalent to taking the geometric realisation of the semi-simplicial space first in the $p$ direction, followed by the $q$ direction, or in the $q$ direction followed by the $p$ direction. This is due to the following homeomorphisms~\cite[1.9,1.10]{EberyRandalWilliams}
	$$
	\|X_{\bt,\bt} \|\cong \|X_{\bt,q}:[p]\mapsto \|X_\bt:[q]\mapsto X_{p,q} \|\|\cong \| X_{p,\bt}:[q]\mapsto \|X_\bt:[p]\mapsto X_{p,q} \|\|.
	$$
\end{defn}

\subsection{Semi-simplicial constructions using monoids and submonoids}

The following description of the geometric bar construction and related definitions loosely follows Chapter 7 of May's \emph{Classifying spaces and fibrations} \cite{MayCSandF}. In this section we view monoids and groups as discrete spaces.

\begin{defn}
	Let $M$ be a monoid and let $X$ and $Y$ be  spaces with a left and right action of $M$ respectively. Then the \emph{bar construction} denoted $B(Y,M,X)$ is the geometric realisation of the semi-simplicial space $B_\bt(Y,M,X)$ given by
	$$
	B_n(Y,M,X)=Y\times M^n \times X.
	$$
	Elements in $B_n(Y,M,X)$ are written as $y[g_1,\ldots, g_n]x$ for $y \in Y$, $g_i \in M$ for $1\leq i \leq n$ and $x \in X$. Face maps are given by
	$$
	\cd_i(y[g_1,\ldots, g_n]x)=\begin{cases}
	yg_1[g_2,\ldots, g_n]x & \text{ if } i=0\\
	y[g_1,\ldots,g_ig_{i+1}, \ldots, g_n]x & \text{ if } 1\leq i\leq n-1 \\
	y[g_1,\ldots, g_{n-1}]g_nx & \text{ if } i=n.
	\end{cases}
	$$
\end{defn}

\begin{defn}\label{def:homotopy quotient}
	Consider the bar construction $B(*,M,Y)$ for $Y$ a space with an action of the monoid $M$ on the left and $*$ a point on which $M$ acts trivially. Define this to be the \emph{homotopy quotient} of $Y$ over $M$ (or $M$ under $Y$) and denote it $B(*,M,Y)=:M\bb Y$. This is the geometric realisation of the semi-simplicial space~$B_\bt(*,M,Y)$ given by
	$$
	B_j(*,M,Y)=*\times M^j \times Y.
	$$
	Elements are written as $[m_1,\ldots, m_j]y$ for $m_i$ in $M$ for $1\leq i \leq j$ and $y$ in $Y$. Face maps are given by
	$$
	\cd_i([m_1,\ldots, m_j]y)=\begin{cases}
	[m_2,\ldots, m_j]y & \text{ if } i=0\\
	[m_1,\ldots, m_im_{i+1}, \ldots, m_j]y & \text{ if } 1\leq i\leq j-1 \\
	[m_1,\ldots, m_{j-1}]m_jy & \text{ if } i=j.
	\end{cases}
	$$
	In the situation of a monoid $M$ acting on a space $Y$ on the right we define the homotopy quotient to be $B(Y,M,*)=:Y\ff M$.
\end{defn}

\begin{exam}
	Consider the bar construction $B(*,N,M)$, for $N$ a submonoid of $M$ acting on $M$ on the left, by left multiplication, and $*$ a point, on which $N$ necessarily acts trivially. Then the homotopy quotient of $M$ over $N$ is $$B(*,N,M)=N\bb M.$$ This is the geometric realisation of the semi-simplicial space $B_\bt(*,N,M)$ given by
	$$
	B_j(*,N,M)=*\times N^j \times M.
	$$
	Elements are written as $[n_1,\ldots, n_j]m$ for $n_i$ in $N$ for $1\leq i \leq j$ and $m$ in $M$. Face maps are given by
	$$
	\cd_i([n_1,\ldots, n_j]m)=\begin{cases}
	[n_2,\ldots, n_j]m & \text{ if } i=0\\
	[n_1,\ldots, n_in_{i+1}, \ldots, n_j]m & \text{ if } 1\leq i\leq j-1 \\
	[n_1,\ldots, n_{j-1}]n_jm & \text{ if } i=j.
	\end{cases}
	$$
	We can build a similar homotopy quotient for a submonoid $N$ acting on $M$ on the right by right multiplication. Then the associated homotopy quotient is the geometric realisation $B(M,N,*)=M\ff N$.
\end{exam}

\begin{lem}\label{lem:homotopy quotient as classifying space}
	The homotopy quotient of a group $G$ or monoid $M$ under a point $*$ is a model for the classifying space of the group or monoid, i.e.~$BG\simeq G\bb * \simeq *\ff G$ and $BM\simeq M\bb * \simeq *\ff M$.
	\begin{proof}
		Writing down the simplices and face maps for the homotopy quotients $G\bb *$ and $G\ff *$ gives exactly the simplices and face maps for the \emph{standard resolution} or \emph{bar resolution} of $G$, which is a model for $BG$ (see e.g.~\cite{Brown}). This holds similarly for the monoid $M$ (see e.g.~\cite[p.~31]{MayCSandF}).
	\end{proof}
\end{lem}

\begin{lem}\label{lem:m bb m htpy pt}
	For a monoid $M$, $M\bb M \simeq *$.
	\begin{proof}
	This is a consequence of \cite[Lemma 1.12]{EberyRandalWilliams} using the augmentation to a point.
	\end{proof}
\end{lem}

\begin{lem}\label{lem:htpy quotient of product}
	Let $N$ be a monoid and~$S$ be a space with right~$N$ action. Suppose~$S$ can be decomposed as $S\cong X\times Y$ and, under this decomposition, the action of~$N$ restricts to a right action on the $Y$ component and trivial action on the $X$ component. Then the map given by the geometric realisation of the levelwise map on the bar construction
	\begin{eqnarray*}
		B_p((X \times Y),N, *)&\to&X \times B_p(Y,N, *)\\
		(x,y)[n_1,\ldots, n_p] &\mapsto& (x,y[n_1,\ldots, n_p])
	\end{eqnarray*}
	\noindent for $x \in X$, $y \in Y$ and $n_i \in N$ for all $i$ is a homotopy equivalence. That is the homotopy quotient satisfies
	$$
	S\ff N \cong (X \times Y)\ff N \simeq X \times (Y \ff N)
	$$
	\noindent where the homotopy equivalence is given by the geometric realisation of the levelwise map on the bar construction
	\begin{eqnarray*}
		B_p((X \times Y),N, *)&\to&X \times B_p(Y,N, *)\\
		(x,y)[n_1,\ldots, n_p] &\mapsto& (x,y[n_1,\ldots, n_p])
	\end{eqnarray*}
	\noindent for $x \in X$, $y \in Y$ and $n_i \in N$ for all $i$.
	\begin{proof}
		The homotopy quotient $S\ff N$ is the geometric realisation of the simplicial space $B_\bt(S,N,*)$ with $j$-simplices given by
		$$
		B_j(S,N,*)=S\times N^j
		$$
		\noindent and face maps given by Definition \ref{def:homotopy quotient}, the first face map $\cd_1$ encoding the right action of $N$ on $S$. Under the decomposition $S\cong X\times Y$ the $j$-simplices are given by
		$$
		B_j(S,N,*)\cong (X\times Y)\times N^j\cong X\times(Y\times N^j)
		$$
		\noindent where the second isomorphism highlights that the action of $N$ on $S$ can be restricted to a right action on $Y$, since the action is trivial on the $X$ component. Note that the second factor is precisely the $j$-simplices in $B_j(Y,N,*)$, and since the face maps act trivially on the $X$ factor, the face maps in $B_j(S,N,*)$ induce face maps in $B_j(Y,N,*)$ under the decomposition. The proof is concluded by taking the geometric realisation of $B_\bt(S,N,*)$ and the geometric realisation of $X \times B_\bt(Y,N,*)$, noting that $\|X\times B_\bt(Y,N,*)\|\simeq X\times \|B_\bt(Y,N,*)\|$.	
	\end{proof}
\end{lem}

\subsection{Semi-simplicial constructions for Artin monoids}

Given an Artin monoid $A^+$ and a parabolic submonoid $M^+$, recall from Section \ref{chap:BG Artin monoids} that~$A^+(M)$ is the set of elements in $A^+$ which do not end in elements in $M^+$ and there is a decomposition as sets (Proposition \ref{prop:decomp}),~$A^+ \cong A^+(M) \times M^+$. This decomposition maps $\alpha$ in $A^+$ to $(\balpha, \beta)$ where $\alpha=\balpha \beta$ (as defined in Remark \ref{rem:balpha def}) and the right action of $M^+$ on~$A^+$ descends to a trivial action on~$A^+(M)$ and a right action on~$M^+$.

In this section we view monoids and sets as discrete spaces.

\begin{prop}\label{prop:homotopy quotient to set}
	The map $$A^+\ff M^+ \to A^+(M)$$ defined levelwise on the bar construction $B_\bt(A^+,M^+, *)$ by
	\begin{eqnarray*}
		B_p(A^+,M^+, *)&\to& A^+(M)\\
		\alpha[m_1,\ldots, m_p] &\mapsto&\balpha
	\end{eqnarray*}
	\noindent is a homotopy equivalence.
	\begin{proof}
		From Proposition \ref{prop:decomp}, $A^+ \cong A^+(M) \times M^+$ as sets, hence as discrete spaces, and this decomposition respects the right action of $M^+$ on $A^+$. Then
		\begin{eqnarray*}
			A^+\ff M^+
			&=&(A^+(M) \times M^+)\ff M^+\\
			&\simeq& A^+(M) \times (M^+\ff M^+) \\
			&\simeq& A^+(M) \times *\\
			&=& A^+(M)
		\end{eqnarray*}
		where the first homotopy equivalence uses Lemma \ref{lem:htpy quotient of product} and the second homotopy equivalence uses Lemma \ref{lem:m bb m htpy pt}. The levelwise map given by the composition of the maps in these two lemmas is precisely the map in the statement.
	\end{proof}
\end{prop}

\begin{prop}\label{prop: equivariant htpy between maps in double homotopy quotient}
	Let $A^+$ be a monoid and $M^+$ be a submonoid. Consider two maps $f$ and $g:A^+ \bb A^+ \to A^+ \bb A^+$ which are both equivariant with respect to the action of $M^+$ on the right of $A^+\bb A^+$. Then there exists an $M^+$ equivariant homotopy between the two maps.
	
	\begin{proof}
		Denote the set of $k$-cells in $A^+\bb A^+$ as  $(A^+ \bb A^+)_k$.
		Let the $k$-cell of $A^+ \bb A^+$ corresponding to geometric realisation of the $k$-simplex~$[p_1,\ldots , p_k ]a$ of~$B_k(*,A^+,A^+)$ (as in Definition~\ref{def:homotopy quotient}) be denoted by the tuple~$(p_1,\ldots , p_k, a)$, with~$p_i$ and~$a$ in~$A^+$. There is a right action of $A^+$ on the~$k$-cells given by  $$(p_1,\ldots , p_k, a)\cdot \mu= (p_1,\ldots , p_k, a\mu).$$ Define the set of {\it elementary $k$-cells} to be those with tuple~$(p_1,\ldots , p_k, e)$ where~$e$ is the identity element in the monoid, and denote this cell~$D(p_1,\ldots , p_k)$. Then every $k$-cell is uniquely determined by an elementary~$k$-cell and an element~$a$ in~$A^+$, since~$(p_1,\ldots , p_k, a)=D(p_1,\ldots , p_k)\cdot a$.  The isomorphism of Proposition \ref{prop:decomp} shows that $A^+ = A^+(M) \times M^+$ and we let $a=\bar{a}m$ under this decomposition. Then we get the following description for $k$-cells:
		
		\ms
		\centerline{\xymatrix@R=3mm@C=-2.5mm {
		(A^+ \bb A^+)_k&\cong& \underset{(p_1,\ldots, p_k)}{\bigcup} D(p_1,\ldots, p_k)\times A^+&\cong&\underset{(p_1,\ldots, p_k)}{\bigcup} D(p_1,\ldots, p_k)\times (A^+(M) \times M^+)\\
		(p_1,\ldots , p_k, a)\ar@{|->}[rr]&& (D(p_1,\ldots , p_k), a)\ar@{|->}[rr]& & (D(p_1,\ldots , p_k), (\bar{a}, m))
		}}
		\ms
		
		Let $f_k$ be the restriction of the map~$f$ to the $k$-cells of~$A^+ \bb A^+$ and similarly for $g_k$. We first define an equivariant homotopy between $f_0$ and~$g_0$. Under the above decomposition,~$(A^+ \bb A^+)_0\cong (A^+(M) \times M^+)$. Consider~$f_0(\balpha)$ and~$g_0(\balpha)$ in~$A^+ \bb A^+$ for $\balpha$ in~$A^+(M)$. Then since~$A^+ \bb A^+ \simeq *$ by Lemma \ref{lem:m bb m htpy pt} it follows that there exists a path between~$f_0(\balpha)$ and~$g_0(\balpha)$:~call this~$h_0(\balpha, t)$ for~$t \in [0,1]$.
		Extend this homotopy to all~$0$-cells by setting~$h_0(\balpha m, t)=h_0(\balpha, t)\cdot m$ for all~$m$ in~$M^+$. Then, since~$f_0$ and~$g_0$ are~$M^+$ equivariant,~$$h_0(\balpha m, 0)=h_0(\balpha, 0)\cdot m=f_0(\balpha)\cdot m=f_0(\balpha m)$$ and similarly~$$h_0(\balpha m, 1)=h_0(\balpha,1)\cdot m=g_0(\balpha)\cdot m=g_0(\balpha m).$$ The homotopy~$h_0(x,t)$ is~$M^+$ equivariant, since~$h_0(x,t)\cdot \mu=h_0(x\mu,t)$ for~$\mu$ in $M^+$.
		
		Now assume that we have built an equivariant homotopy~$h_{k-1}(x,t)$ on the~$(k-1)$-skeleton and we show how to extend it to the~$k$-cells. The homotopy~$h_{k-1}(x,t)$ satisfies~$h_{k-1}(x,0)=f_{k-1}(x)$ and~$h_{k-1}(x,1)=g_{k-1}(x)$. For some $\balpha$ in $A^+(M)$, consider the~$k$-cell~$D(p_1,\ldots, p_k)\cdot \balpha$. Then its boundary consists of $(k-1)$-cells and it follows that $h_{k-1}$ defines a homotopy
		$$
		(\cd (D(p_1,\ldots, p_k))\cdot \balpha) \times I \to  A^+ \bb A^+
		$$
		and the maps $f_k$ and $g_k$ also define maps
		\begin{eqnarray*}
			f_k: ((D(p_1,\ldots, p_k))\cdot \balpha) \times \{0\} &\to& A^+ \bb A^+\\
			g_k: ((D(p_1,\ldots, p_k))\cdot \balpha) \times \{1\} &\to& A^+ \bb A^+.
		\end{eqnarray*}
		The union of these three maps defines a map from~$\cd((D(p_1,\ldots, p_k)\cdot \balpha) \times I)$ to $A^+ \bb A^+$, but this boundary is a $(k-1)$-sphere and so, since $A^+ \bb A^+$ is contractible the $(k-1)$-sphere bounds a $(k)$-disk. We can compatibly extend the map over this disk to create the required homotopy
		$$
		h_k: (D(p_1,\ldots, p_k)\cdot \balpha) \times I \to  A^+ \bb A^+
		$$
		which agrees on the boundary with the three maps above. Now define $h_k$ on any~$k$-cell~$D(p_1,\ldots, p_k)\cdot~\balpha m$ by the following: for $x$ in $D(p_1,\ldots, p_k)\cdot \balpha$ we set~$$h_k(x\cdot m,t)=h_k(x,t)\cdot m.$$ Then by construction $h_k$ is $M^+$ equivariant and, since both $f_k$ and $g_k$ are $M^+$ equivariant,~$h_k$ satisfies $h_k(x,0)=f_k$ and~$h_k(x,1)=~g_k$.
	\end{proof}
\end{prop}

\begin{defn}\label{def:double homotopy quotient}
	Given a monoid $M$ and two submonoids $N_1$ and $N_2$ we can define the \emph{double homotopy quotient} $N_1 \bb M \ff N_2$ to be the geometric realisation of the bi-semi-simplicial space (recall Definition \ref{def:bisemi}) defined by taking the two simplicial directions arising from the bar constructions $B_\bt(*,N_1,M)$ and $B_\bt(M, N_2,*)$. The $p,q$ level of the associated bi-semi-simplicial space $X_{\bt \bt}$ has simplices
	$$
	X_{p,q}=N_1^p \times M \times N_2^q
	$$
	and face maps inherited from $B_\bt(*,N_1,M)$ in the $p$ direction ($\cd_{p,\bt}$) and $B_\bt(M, N_2,*)$ in the $q$ direction ($\cd_{\bt,q}$). Then $[n_1,\ldots, n_p]m[n'_1,\ldots, n'_q]$ represents an element in the $p,q$ level, where $n_i$ in $N_1$ and $n'_j$ in $N_2$ for $1\leq i\leq p$ and $1\leq j \leq q$. We note that the face maps on the left and right commute, since the only maps which act on the same coordinates are $\cd_{p,\bt}$ in the $p$ direction and $\cd_{\bt,0}$ in the $q$ direction and these commute:
	\begin{eqnarray*}
		\cd_{p,\bt}(\cd_{\bt,0}([n_1,\ldots, n_p]m[n'_1,\ldots, n'_q]))&=&\cd_{p,\bt}([n_1,\ldots, n_p]mn'_1[n'_2,\ldots, n'_q])\\
		&=&[n_1,\ldots, n_{p-1}] n_pmn'_1[n'_2,\ldots, n'_q]\\
		&=&\cd_{\bt,0}(\cd_{p,\bt}([n_1,\ldots, n_p]m[n'_1,\ldots, n'_q])).
	\end{eqnarray*}
\end{defn}

\section{Preliminaries concerning the sequence of Artin monoids} \label{section:artin}
This section introduces notation used throughout the remainder of the proof.
 
We consider the sequence of Artin monoids and inclusions
\begin{equation}\label{eq:monoid seq with A0}
A^+_0 \hookrightarrow A^+_1 \hookrightarrow A^+_2 \hookrightarrow \cdots \hookrightarrow A^+_n \hookrightarrow \cdots	
\end{equation}
\noindent
with Artin monoid systems $(A^+_n,\Sigma_n)$ given by the following diagrams. Here the Artin Monoid~$A_i$ corresponds to the Coxeter group~$W_i$, as defined in \cite{Hepworth}, and so we denote the corresponding Coxeter diagrams~$\D_{W_i}$.

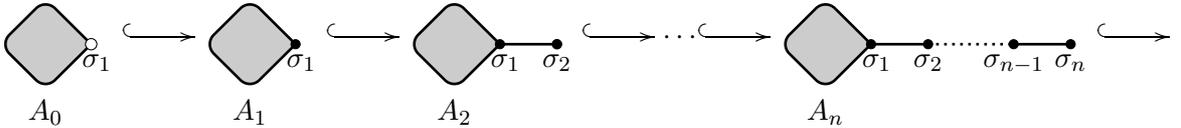
\centerline{	\xymatrix@R=3mm {
		\begin{tikzpicture}[scale=0.15, baseline=0]
		\draw[line width=1, fill = white!80!black, rounded corners=5 pt]
		(5,0) --  (1,4) -- (-3,0) -- (1,-4) -- (5,0);
		\draw (1,-6) node {$A_0$};
		\draw[fill= white] (5,0) circle (0.5);		
		\draw (5.5,0) node[below] {$\sm_1$};
		\end{tikzpicture} \ar@{^{(}->}[r]
		&
		\begin{tikzpicture}[scale=0.15, baseline=0]
		\draw[line width=1, fill = white!80!black, rounded corners=5 pt]
		(5,0) --  (1,4) -- (-3,0) -- (1,-4) -- (5,0);
		\draw (1,-6) node {$A_1$};
		\draw[fill= black] (5,0) circle (0.4);	
		\draw (5.5,0) node[below] {$\sm_1$};
		\end{tikzpicture} \ar@{^{(}->}[r]
		&
		\begin{tikzpicture}[scale=0.15, baseline=0]
		\draw[line width=1, fill = white!80!black, rounded corners=5 pt]
		(5,0) --  (1,4) -- (-3,0) -- (1,-4) -- (5,0);
		\draw[line width=1] (5,0) -- (10,0);
		\draw[fill= black] (10,0) circle (0.4);
		\draw (1,-6) node {$A_2$};
		\draw[fill= black] (5,0) circle (0.4);
		\draw (5.5,0) node[below] {$\sm_1$};
		\draw[fill= black] (10,0) circle (0.4) node[below] {$\sm_2$};
		\end{tikzpicture}  \ar@{^{(}->}[r]
		&
		\cdots  \ar@{^{(}->}[r]
		&
		\begin{tikzpicture}[scale=0.15, baseline=0]
		\draw[line width=1, fill = white!80!black, rounded corners=5 pt]
		(5,0) --  (1,4) -- (-3,0) -- (1,-4) -- (5,0);
		\draw[fill= black] (5,0) circle (0.4);
		\draw (5.5,0) node[below] {$\sm_1$};
		\draw[line width=1] (5,0) -- (10,0);
		\draw[fill= black] (10,0) circle (0.4) node[below] {$\sm_2$};
		\draw[line width=1, dotted] (10,0) -- (17.5,0);
		\draw[fill= black] (17.5,0) circle (0.4)node[below] {$\sm_{n-1}$};
		\draw[line width=1] (17.5,0) -- (22.5,0);
		\draw[fill= black] (22.5,0) circle (0.4)node[below] {$\sm_{n}$};
		\draw (1,-6) node {$A_n$};
		\end{tikzpicture}  \ar@{^{(}->}[r]
		&
}}

\begin{defn}
	Let $(A_0, \Sigma_0)$ be the Artin system corresponding to the Coxeter diagram $\D_{W_1}$, but with the vertex $\sm_1$ and all edges which have vertex $\sm_1$ at one end removed. We depict the diagram as above. Note that~$A_0\hookrightarrow A_1$.
\end{defn}

\begin{rem}
With the generating sets corresponding to the above sequence of diagrams, for all~$p$ every generator and hence every word in the monoid~$A_p^+$ commutes with~$\sm_j$ for~$j\geq p+2$.	
\end{rem}

We now apply the theory developed in Section~\ref{sec:generic artin monoid divisors and cosets 2} to the specific case of a monoid~$A^+_n$ in the sequence of monoids and inclusions \eqref{eq:monoid seq with A0} and the submonoid of~$A^+_n$, given by a previous monoid in the sequence $A^+_p$ where $p<n$. We adopt the following notation for the remainder of this paper. The generating set of~$A^+_n$ will always be given by~$\Sigma_n$, the generating set specified by the diagram $\D_{W_n}$.

\begin{itemize}
	\item Let $\EndMon_p(\alpha)=\EndMon_{A_p}(\alpha)$ and $\EndGen_p (\alpha)=\EndGen_{A_p} (\alpha)$ for $\alpha$ in $A^+_n$, as in Definition \ref{def:end sets}. Then
	\begin{eqnarray*}
		\EndGen_p (\alpha)&=&\{\sm_s \mid s \in S_{A_p^+},  \sm_s \preceq_R \alpha\}\\
		\EndMon_{p} (\alpha)&=&\{\beta \in A_p^+ \mid \beta \preceq_R \alpha\}.
	\end{eqnarray*}
	\item Let $A^+(n;p)$ be the set $A^+(M)$ for $A^+=A^+_n$ and $M=A^+_p$ as in Definition~\ref{def:cosets} (this is the set of elements in $A^+_n$ that do not end in a non trivial element in $A^+_p$).
	\item Let the equivalence class of $\alpha$ in $A^+_n$ under the relation $\approx$ with respect to the submonoid $A^+_p$ (Definition \ref{def:monoid coset}) be denoted $[\alpha]_p$ as opposed to $[\alpha]_{A_p}$. Then $[\alpha]_p$ is the equivalence class of $\alpha$ under $\approx$, the transitive closure of the relation $\sim$ on $A^+_n$ given by
	$$
	\alpha_1 \sim \alpha_2 \iff \alpha_1 \beta_1=\alpha_2 \beta_2 \text{ for some }\beta_1\text{ and }\beta_2\text{ in }A_p^+.
	$$
\end{itemize}

\noindent Then we have from Lemma \ref{lem:barcosetequality} that the equivalence classes under $\approx$ with respect to the submonoid $A^+_p$ are in one to one correspondence with the set $A^+(n;p)$. Recall from Remark \ref{rem:balpha def} that if $\beta$ is the least common multiple of $\EndMon_p(\alpha)$ then we define $\balpha$ in $A^+_n$ to be the element such that $\alpha=\balpha\beta$. Then $A^+(n;p)$ is the set of all such $\balpha$ and for all $\alpha_1$ and $\alpha_2$ in $A^+_n$:
\begin{equation*}
[\alpha_1]_p=[\alpha_2]_p \iff \overline{\alpha_1}=\overline{\alpha_2}.
\end{equation*}

\noindent We also have from Proposition \ref{prop:decomp} the decomposition
$$
A^+_n \cong A^+(n;p) \times A^+_p\text{ for all }p < n.
$$

\section{The semi-simplicial space {$\protect{\Aforheader}$}}\label{sec:introduce tri ssspace}

We now build the semi-simplicial space $\A^n_\bt$ as promised in Section \ref{section:setup}.

\begin{defn}\label{def:complex}
	Define a semi-simplicial space $\C^n_\bt$ by, for $0\leq p \leq(n-1)$, setting levels $\C^n_p$ to be the discrete space of equivalence classes $A^+_n\slash \approx$ where the equivalence relation is taken with respect to the submonoid $A^+_{n-p-1}$, i.e.~$\approx$ is the transitive closure of the relation $\sim$ on $A^+_n$ given by
	$$
	\alpha_1 \sim \alpha_2 \iff \alpha_1 \beta_1=\alpha_2 \beta_2 \text{ for some }\beta_1\text{ and }\beta_2\text{ in }A^+_{n-p-1}.
	$$
	\noindent Face maps are given by
	\begin{eqnarray*}
		\cd^p_k: \C^n_p &\to& \C^n_{p-1} \text{ for } 0\leq k \leq p\\
		\cd^p_k:[\alpha]_{n-p-1} &\mapsto& [\alpha(\sm_{n-p+k}\sm_{n-p+k -1}\ldots \sm_{n-p+1})]_{n-p}.
	\end{eqnarray*}
For example,~$\cd^p_0$ acts on the equivalence class representative by right multiplication by~$e$, and~$\cd^p_p$ acts by right multiplication by~$\sm_n\ldots\sm_{n-p-1}$.
The motivation for this choice of face maps follows Hepworth, as discussed in \cite[Example 35]{Hepworth}.
\end{defn}
\begin{lem}
	The face maps of Definition \ref{def:complex} are well defined.
\begin{proof}
	We want that if $[\alpha]_{n-p-1}=[\eta]_{n-p-1}$ then~$\cd^p_k([\alpha]_{n-p-1})=	\cd^p_k([\eta]_{n-p-1})$. If $[\alpha]_{n-p-1}=[\eta]_{n-p-1}$, then $\balpha=\bar{\eta}$ where the bar is taken with respect to $A_{n-p-1}^+$. Set $\balpha=\gamma$ (recall the definition of $\balpha$ from Remark \ref{rem:balpha def}). It follows that there exist $a$ and $b$ in $A^+_{n-p-1}$ such that~$\alpha=\gamma a$ and $\eta~=~\gamma b$. Then since $a$ and $b$ only contain letters in $A_{n-p-1}^+$ and all of these letters commute with~$(\sm_{n-p+k}\sm_{n-p+k -1}\ldots \sm_{n-p+1})$ it follows that $a$ and $b$ letterwise commute with the face map. Taking equivalence classes with respect to $A^+_{n-p}$ therefore gives
	\begin{eqnarray*}
		& &\lbrack \alpha (\sm_{n-p+k}\sm_{n-p+k -1}\ldots \sm_{n-p+1}) \rbrack_{n-p}\\
		&=& \lbrack(\gamma a )(\sm_{n-p+k}\sm_{n-p+k -1}\ldots \sm_{n-p+1} )\rbrack_{n-p}\\
		&=& \lbrack\gamma (\sm_{n-p+k}\sm_{n-p+k -1}\ldots \sm_{n-p+1}) a\rbrack_{n-p}\\
		&=& \lbrack\gamma (\sm_{n-p+k}\sm_{n-p+k -1}\ldots \sm_{n-p+1}) \rbrack_{n-p}
	\end{eqnarray*}
	and similarly
	\begin{eqnarray*}
		& &\lbrack \eta (\sm_{n-p+k}\sm_{n-p+k -1}\ldots \sm_{n-p+1}) \rbrack_{n-p}\\
		&=& \lbrack\gamma (\sm_{n-p+k}\sm_{n-p+k -1}\ldots \sm_{n-p+1}) \rbrack_{n-p}
	\end{eqnarray*}
	and so the face maps are well defined.
\end{proof}
\end{lem}

\begin{lem}\label{lem:simplicial identities for complex}
	The face maps $\{\cd_k^p\}$ on $\C_\bt^n$ defined in Definition \ref{def:complex} satisfy the simplicial identities, that is, for $0\leq i<j\leq p$:
	$$
	\cd^{p-1}_i\cd^p_j=\cd^{p-1}_{j-1}\cd^p_i:\C^n_p\to \C^n_{p-2}
	$$
	\begin{proof}
		For ease of notation in the proof, we denote $(n-p)$ as $r$. Then the left hand side acts as follows
		
		\centerline{\xymatrix@R=2mm{		
				\C^n_p \ar[r]^-{\cd^p_j}&   \C^n_{p-1} \ar[r]^-{\cd^{p-1}_i}&  \C^n_{p-2}\\
				{[\alpha]_{r-1}} \ar@{|->}[r]^-{\cd^p_j}& {[\alpha (\sm_{r+j}\ldots \sm_{r+1})]_{r}} \ar@{|->}[r]^-{\cd^{p-1}_i}& {[\alpha (\sm_{r+j}\ldots \sm_{r+1})(\sm_{r+i+1}\ldots \sm_{r+2})]_{r+1}}.
		}}
		\ms
		In comparison the right hand side acts as follows
		
		\centerline{\xymatrix@R=2mm{		
				\C^n_p \ar[r]^{\cd^p_i}&   \C^n_{p-1} \ar[r]^{\cd^{p-1}_{j-1}}&  \C^n_{p-2}\\
				{[\alpha]_{r-1}} \ar@{|->}[r]^(.4){\cd^p_i}& {[\alpha (\sm_{r+i}\ldots \sm_{r+1})]_{r}} \ar@{|->}[r]^(.4){\cd^{p-1}_{j-1}}& {[\alpha (\sm_{r+i}\ldots \sm_{r+1})(\sm_{r+j}\ldots \sm_{r+2})]_{r+1}}
		}}
		\ms
		
		Let $x=(\sm_{r+j}\ldots \sm_{r+1})(\sm_{r+i+1}\ldots \sm_{r+2})$ and $y=(\sm_{r+i}\ldots \sm_{r+1})(\sm_{r+j}\ldots \sm_{r+2})$. Note that for $0\leq k <j$ we have
		$$
		(\sm_{r+j}\ldots \sm_{r+1}){\sm_{r+k+1}}=\sm_{r+k}(\sm_{r+j}\ldots \sm_{r+1})
		$$
		\noindent from manipulation of the words using the braiding relations in the monoid. Reiterating this gives us the first equality in the following:
		\begin{eqnarray*}
		x&=&(\sm_{r+j}\ldots\sm_{r+1})(\sm_{r+i+1}\ldots \sm_{r+2})\\
		&=&	(\sm_{r+i}\ldots\sm_{r})(\sm_{r+j}\ldots\sm_{r+1})\sm_{r+2}\\
		&=& (\sm_{r+i}\ldots\sm_{r})(\sm_{r+j}\ldots\sm_{r+3})(\sm_{r+2}\sm_{r+1}\sm_{r+2})\\
		&=&(\sm_{r+i}\ldots\sm_{r})(\sm_{r+j}\ldots\sm_{r+3})(\sm_{r+1}\sm_{r+2}\sm_{r+1})\\
		&=&(\sm_{r+i}\ldots\sm_{r})\sm_{r+1}(\sm_{r+j}\ldots\sm_{r+3})\sm_{r+2}\sm_{r+1}\\
		&=&(\sm_{r+i}\ldots\sm_{r}\sm_{r+1})(\sm_{r+j}\ldots\sm_{r+3}\sm_{r+2})\sm_{r+1}\\
		&=&y\sm_{r+1}
		\end{eqnarray*}
		The result follows since we are taking the equivalence relation with respect to the submonoid $A^+_{r+1}$.
	\end{proof}
\end{lem}

\begin{lem}\label{lem:levelequiv}
	Recall the notation $A^+(n;n-p-1)$, as defined in Section \ref{section:artin}. Then the realisation of the map defined levelwise on the bar construction by
	\begin{eqnarray*}
		B_p(A_n^+,A_{n-p-1}^+, *)&\to& A^+(n;n-p-1)\\
		\alpha[m_1,\ldots, m_p] &\mapsto&\balpha
	\end{eqnarray*}
	where $\alpha \in A^+_n$, $m_i\in A^+_{n-p-1}$ for all i and $\alpha=\balpha \beta$ for $\balpha \in A^+(n;n-p-1)$ and~$\beta \in A_{n-p-1}^+$ is a homotopy equivalence. That is the $p$th level of $\C^n_\bt$ satisfies
	\begin{equation*}
	A^+_n \ff A^+_{n-p-1}\simeq  A^+(n;n-p-1)=\C^n_p.
	\end{equation*}
	\begin{proof}
		This is a direct application of Proposition \ref{prop:homotopy quotient to set} and Proposition \ref{prop:decomp} which gives the decomposition $A^+_n \cong A^+(n;n-p-1) \times A^+_{n-p-1}$.
	\end{proof}
\end{lem}

\begin{defn}\label{def:new ssspace}
	Let $\A^n_\bt$ be the semi-simplicial space with $p$th level the homotopy quotient $\A^n_p = A^+_n \bb \C^n_p$, where the action of $A^+_n$ on $A^+(n;n-p-1)$ is given by
	$$
	a\cdot [\alpha]_{n-p-1}=[a\alpha]_{n-p-1} \text{ for }a, \alpha \in A^+_n
	.$$

	\noindent The face maps are denoted by $\cd^p_k \text{ for } 0\leq k \leq p$
	\begin{eqnarray*}
		\cd^p_k: \A^n_p &\to& \A^n_{p-1} \\
		\cd^p_k: A^+_n \bb \C^n_p &\to&  A^+_n \bb \C^n_{p-1}
	\end{eqnarray*}
	and $\cd^p_k$ acts as the face map $\cd^p_k$ from Definition \ref{def:complex} on the $\C^n_p$ factor of each simplex in the homotopy quotient, and as the identity on the other factors.
	
	\noindent Diagrammatically, $\A^n_\bt$ can be drawn as:
	
	\centerline{	\xymatrix@R=5mm@C=10mm {
			A^+_n \bb  \C^n_{n-1}\ar@<-4ex>[d]\ar@<-2ex>[d] \ar@<2ex>[d]_{\cdots}\ar@<4ex>[d] & \A^n_{n-1}\\
			A^+_n \bb  \C^n_{n-2}\ar@<-4ex>[d]\ar@<-2ex>[d] \ar@<2ex>[d]_{\cdots}\ar@<4ex>[d] & \A^n_{n-2}\\
			A^+_n \bb  \C^n_{n-3}
			\ar@<-4ex>[d]\ar@<-2ex>[d] \ar@<2ex>[d]_{\cdots}\ar@<4ex>[d] & \A^n_{n-3} \\
			\vdots
			\ar@<0ex>[d]\ar@<-2ex>[d] \ar@<2ex>[d]& \vdots \\
			A^+_n \bb  \C^n_{1}  \ar@<-1ex>[d]\ar@<1ex>[d]& \A^n_{1}\\
			A^+_n \bb  \C^n_{0}& \A^n_{0}
	}}
\end{defn}
	
\begin{lem}\label{lem:well defined face maps}
	The factorwise definition of the face maps ${\cd}^p_k$ in Definition \ref{def:new ssspace} gives well defined maps on the homotopy quotients at each level of~$\A_\bt^n$.
	
	\begin{proof}		
	The set of $j$-simplices in $A^+_n \bb \C^n_p$ is identified with $(A^+_n)^j \times \C^n_p$ and a generic element in this set is given by $\lbrack a_1,\ldots, a_j \rbrack \lbrack\alpha\rbrack_{n-p-1}$, where the $a_i$ and $\alpha$ are in $A^+_n$. Then the map ${\cd}^p_k$ acts on this simplex as
	$$
	\cd^p_k(\lbrack a_1,\ldots , a_j\rbrack\lbrack\alpha\rbrack_{n-p-1}) \mapsto {\lbrack a_1,\ldots , a_j\rbrack\lbrack\alpha(\sm_{n-p+k}\sm_{n-p+k -1}\ldots \sm_{n-p+1})\rbrack_{n-p}}
	$$
	and since the multiplication by $(\sm_{n-p+k}\sm_{n-p+k -1}\ldots \sm_{n-p+1})$ is on the right it follows that ${\cd}^p_k$ commutes with all face maps of the bar construction $B_\bt(*,A_n^+, \C^n_p)$ for each $k$. Therefore the definition of ${\cd}^p_k$ on the simplicial level induces a map on the homotopy quotient $A^+_n \bb \C^n_p$.
	\end{proof}
\end{lem}

\begin{lem}\label{lem:simplicial idesntities for new sss}
	The face maps $\cd_k^p$ on $\A_\bt^n$ defined in Definition \ref{def:new ssspace} satisfy the simplicial identities, that is for $0\leq i<j\leq p$:
	$$
	\cd^{p-1}_i\cd^p_j=\cd^{p-1}_{j-1}\cd^p_i.
	$$
	\begin{proof}
		This follows directly from the fact that the simplicial identities are satisfied for $\C^n_\bt$ (Lemma \ref{lem:simplicial identities for complex}), since the face maps for $\A^n_\bt$ are defined via the maps for $\C^n_\bt$.
	\end{proof}
\end{lem}

We now show that there exist homotopy equivalences $\A^n_p \simeq BA^+_{n-p-1}$ for every~$p \geq 0$, as promised in Section \ref{section:setup}.

\begin{lem}\label{lem:homotopy equivalence levels to previous in sequence}
	
	Consider the levelwise maps on $(j,k)$-simplices of $A^+_n \bb A^+_n \ff A^+_{n-p-1}$:
	\begin{eqnarray*}
		(A^+_n \bb A^+_n \ff A^+_{n-p-1})_{(j,k)}&\to& (A^+_n\bb \C_p^n)_j\\
		\lbrack a_1,\ldots, a_j\rbrack \alpha\lbrack a'_1,\ldots, a'_k\rbrack &\mapsto& \lbrack a_1,\ldots, a_j\rbrack\balpha
	\end{eqnarray*}
	\noindent and the projection
	\begin{eqnarray*}
		(A^+_n \bb A^+_n \ff A^+_{n-p-1})_{(j,k)}  &\to& (* \ff A^+_{n-p-1})_k\\
		\lbrack a_1,\ldots, a_j\rbrack \alpha \lbrack a'_1,\ldots, a'_k\rbrack &\mapsto& *\lbrack a'_1,\ldots, a'_k\rbrack
	\end{eqnarray*}
	\noindent
	where $\alpha$ and $a_i \in A^+_n$, $a'_i \in A^+_{n-p-1}$, and $\alpha=\balpha \beta$ for $\balpha \in A^+(n;n-p-1)$ and~$\beta \in A_{n-p-1}^+$. Then these maps are homotopy equivalences
	$$
	A^+_n \bb A^+_n \ff A^+_{n-p-1}\simeq 	\A^n_p
	$$
	\noindent and
	$$
	A^+_n \bb A^+_n \ff A^+_{n-p-1}\simeq BA^+_{n-p-1}
	$$
	respectively. That is, the $p$th level of the space $\A^n_\bt$ satisfies $$
	\A^n_p \simeq A^+_n \bb A^+_n \ff A^+_{n-p-1}\simeq BA^+_{n-p-1}.
	$$
	\begin{proof}
		From Lemma \ref{lem:levelequiv}, $\C^n_p = A^+(n;n-p-1) \simeq A^+_n \ff A^+_{n-p-1}$, and this induces
		$$\A^n_p=A^+_n \bb \C^n_p  \simeq A^+_n \bb A^+_n \ff A^+_{n-p-1}$$ \noindent with the homotopy equivalence given by the required map. We then have the following
		$$
		\A^n_p\simeq A^+_n \bb A^+_n \ff A^+_{n-p-1} = (A^+_n \bb A^+_n) \ff A^+_{n-p-1} \simeq * \ff A^+_{n-p-1}=BA_{n-p-1}^+.
		$$
		\noindent The central equality is due to the fact that the double homotopy quotient is the geometric realisation of a bi-simplicial-set and therefore we can take the realisation in either direction first. The second map in the previous equation is a homotopy equivalence by Lemma \ref{lem:m bb m htpy pt}.
	\end{proof}
\end{lem}

We now define the map from the geometric realisation of $\A^n_\bt$ to the classifying space $BA^+_n$ promised in Section \ref{section:setup}:

$$
\| \A^n_\bt \| \overset{\|\phi_\bt \|}{\longrightarrow} BA^+_n
$$
In Section \ref{section:highcon} we will show that  $\|\phi_\bt\|$ is highly connected.

\begin{lem}
	The geometric realisation~$\| \A^n_\bt \|$ satisfies $\| \A^n_\bt \| \cong A^+_n \bb \| \C^n_\bt \|$.
	\begin{proof}
		The face maps in the bar construction $B_\bt(*, A^+_n, \C^n_p)$ for the homotopy quotient in $\A^n_p=A^+_n\bb \C^n_p$ commute with the face maps in $\C^n_\bt$ (see the proof of Lemma~\ref{lem:well defined face maps}) and therefore with the face maps of $\A^n_\bt$. Therefore the two simplicial directions create a bi-semi-simplicial space and one can realise in either direction first, as in Definition \ref{def:bisemi geom rel}. Realising by taking the homotopy quotients $\A^n_p=A^+_n\bb \C^n_p$ before realising in the $\A^n_\bt$ direction first (which has face maps induced by those of $\C^n_\bt$) gives the left hand side. Realising in the $\C^n_\bt$ direction before taking the homotopy quotient $A^+_n \bb \| \C^n_\bt \|$ gives the right hand side.
	\end{proof}
\end{lem}

Recall that $A^+_n\bb *$ is a model for $BA^+_n$. We therefore define $\|\phi_\bt \|$ as a map from~$A^+_n \bb \| \C^n_\bt \|$ to  $A^+_n\bb *$.

\begin{defn}
	Define $\phi_\bt$ to be the semi-simplicial map from the bar construction~$B_\bt(*,A^+_n,\| \C^n_\bt \|)$ to the bar construction $B_\bt(*,A^+_n,*)$, obtained by collapsing $ \| \C^n_\bt \|$ to a point:	
	\begin{eqnarray*}
		\phi_p: B_p(*,A^+_n,\| \C^n_\bt \|) &\to& B_p(*,A^+_n,*)\\
		\lbrack a_1,\ldots , a_p\rbrack a &\mapsto& \lbrack a_1,\ldots , a_p\rbrack *
	\end{eqnarray*}
	
	\noindent where $a_i$ is in $A^+_n$ for all $i$, and $a$ is in $\| \C^n_\bt \|$.	Then the geometric realisation $\|\phi_\bt\|$ maps the homotopy quotient $A^+_n \bb \| \C^n_\bt \|$ to the homotopy quotient $A^+_n\bb *\simeq BA^+_n$.
\end{defn}

\begin{prop}\label{prop:highcon}
	If $\| \C^n_\bt \|$ is $(k-1)$-connected then the map $\|\phi_\bt\|$ is $k$-connected.
	\begin{proof}
		From \cite[Lemma 2.4]{EberyRandalWilliams} a semi-simplicial map~$f_\bt: X_\bt \to Y_\bt$ satisfies that~$\|f_\bt\|$ is~$k$-connected if~$f_p: X_p \to Y_p$ is~$(k-p)$ connected for all $p \geq 0$. The map $\|\phi_\bt\|$ is defined level-wise as the projection $$\phi_p:(A^+_n)^p \times \| \C^n_\bt \|\to (A^+_n)^p.$$ \noindent Therefore since $\| \C^n_\bt \|$ is $(k-1)$-connected it follows that $\phi_p$ is $k$-connected and in particular it is $(k-p)$-connected for all $p\geq0$. Thus the geometric realisation $\|\phi_\bt\|$ is $k$-connected.
	\end{proof}
\end{prop}

\section{High connectivity}\label{section:highcon}
This section is concerned with the proof of the following theorem.

\begin{thm}\label{thm:highly connected}
	The geometric realisation $\| \C^n_\bt \|$ of the semi-simplicial space $\C^n_\bt$ is~$(n-2)$-connected for all $n$, i.e.~$\pi_i(\| \C^n_\bt \|)=0$ for $0\leq i \leq n-2$.
\end{thm}

Combining this Theorem with Proposition \ref{prop:highcon}, it follows that the map $\|\phi_\bt \|$ is~$(n-1)$-connected as promised in Section\ref{section:setup}. For the remainder of this paper, we will refer to the geometric realisation of the semi-simplicial space as a \emph{complex} (the geometric realisation is, by definition, a cell complex: note that it is not necessarily a simplicial complex).

\subsection{Union of chambers argument}

There is a specific argument, called a \emph{union of chambers argument} that is often used to prove high connectivity of a complex. It is closely related to the notion of \emph{shellability}.

In \cite{Davis}, Davis used a union of chambers argument to prove that the Davis complex $\Sigma_W$ associated to a Coxeter group is contractible. He did this by showing that the Davis complex is an example of a so called \emph{basic construction}. Hepworth's high connectivity results relating to homological stability for Coxeter groups \cite{Hepworth} also used such an argument. In \cite{Paris}, Paris used a union of chambers argument to show that the universal cover of an analogue of the Salvetti complex for certain Artin monoids is contractible, proving the $K(\pi,1)$ conjecture for finite type Artin groups. In this chapter we use a similar union of chambers argument to prove high connectivity. Loosely, the argument consists of breaking the complex up into high dimensional \emph{chambers} and considering how connectivity changes as they are glued together to create the complex. Whilst applying the argument in the case of Artin monoids and the complex we have constructed, numerous technical challenges arise, leading to the proof being split into many separate cases.

To prove high connectivity in our setup we use a union of chambers argument applied to the complex $\|\C^n_\bt\|$. Recall that $\|\C^n_\bt\|$ has dimension~$n-1$. We filter the top dimensional simplices by the natural numbers as follows:

\begin{defn}
	For $k$ in $\mathbb{N}$ we define $\C^n (k)$ as follows:
	\begin{equation*}
	\C^n(k)=\bigcup_{\substack{\alpha \in A^+_n, \\ \lng(\alpha)\leq k}}\lb\alpha\rb_0
	\end{equation*}
	where $\lb\alpha\rb_0$ is the $(n-1)$ simplex in $\| \C^n_\bt \|$ represented by $[\alpha]_0$ in $\C^n_{n-1}$. 
\end{defn}

\begin{rem}
	Note that every simplex in $\| \C^n_\bt \|$ arises as a face of some $\lb\alpha\rb_0$, since smaller simplices are represented by some $\lb\tau\rb_k$ for $k>0$ and this is a face of $\lb\tau\rb_0$. Then $\|\C^n_\bt\|$ is given by $\colim_{k\to \infty}\C^n(k)$.
\end{rem}

The union of chambers argument relies on the following two steps:

\begin{enumerate}[(A)]
	\item \label{point A} If $\lng(\alpha)=k+1$ then $\lb\alpha\rb_0 \cap \C^n(k)$ is a non-empty union of top dimensional faces of $\lb\alpha\rb_0$.
	\item \label{point B} If $\lng(\alpha)=\lng(\beta)=k+1$ and~$\alpha \neq \beta$ then $\lb\alpha\rb_0 \cap \lb\beta\rb_0 \subseteq \C^n(k)$.
\end{enumerate}

\begin{prop} \label{prop: connectivity of C}
	If (\ref{point A}) and (\ref{point B}) hold then $\| \C^n_\bt \|$ is homotopy equivalent to a wedge of~$(n-1)$ spheres, and in particular is $(n-2)$ connected.
	\begin{proof}
		We build up $\| \C^n_\bt \|$ by increasing $k$ in $\C^n(k)$. We start at $\C^n(0)=\lb e \rb_0$, which is contractible. At each step we build up from $\C^n(k)$ to $\C^n({k+1})$ by adding the set of simplices represented by words in $A^+_n$ of length $(k+1)$:
		$$
		\bigcup_{\substack{\alpha \in A^+_n, \\ \lng(\alpha)= k+1}}\lb\alpha\rb_0.
		$$
		Then point (\ref{point A}) says that when $\lb\alpha\rb_0$ is added to $\C^n(k)$, the intersection is a non-empty union of facets of $\lb \alpha \rb_0$. Therefore either the homotopy type doesn't change upon adding the simplex (if not all facets are in the intersection), or the homotopy changes and this change is described by the possible addition of an $(n-1)$ sphere (if all facets are in the intersection).
		Point (\ref{point B}) then says that adding the entirety of the above union to $\C^n(k)$ at the same time only changes the homotopy type in the sense that the individual simplices change it, since each two simplices intersect within $\C^n(k)$. Therefore at each stage we change the homotopy type by at most the addition of several $(n-1)$ spheres and it follows that $\| \C^n_\bt \|$ is $(n-2)$ connected.
	\end{proof}
\end{prop}

\begin{rem}
	The length function gives a partial order on the top dimensional simplices of~$\| \C^n_\bt \|$. By~(\ref{point B}), any linear extension of this partial order to a total order will still satisfy~(\ref{point A}). In this case, the ordering is called a \emph{shelling} (see~\cite{Bjorner}), which we know to be highly connected: giving an alternative proof to the previous proposition.
\end{rem}

The remainder of this section is therefore devoted to the proof of point (\ref{point A}) and point (\ref{point B}). The proof of point (\ref{point A}) is split into several subsections.

\subsection{Proof of point (\ref{point A}): facets of $\lb\alpha\rb_0$}
Recall that the top dimensional faces of a simplex are called \emph{facets}. We start the proof of point (\ref{point A}) with a discussion of the facets of a simplex $\lb\alpha\rb_0$. Consider the face maps 	
\begin{eqnarray*}
	\cd^{n-1}_q: \C^n_{n-1} &\to& \C^n_{n-2} \\
	\cd^{n-1}_q:\lb\alpha\rb_{0} &\mapsto& \lb\alpha\sm_{q+1}\sm_q\ldots \sm_{2}\rb_{1}
\end{eqnarray*}
for $0\leq q \leq n-1$. The map $\cd^{n-1}_0$ is right multiplication by the identity.

\noindent Under these face maps the facets of $\lb\alpha\rb_0$ are given by
\begin{equation*}
\lb\alpha\rb_1, \lb\alpha \sm_2\rb_1,  \lb\alpha \sm_3\sm_2\rb_1,  \lb\alpha \sm_4 \sm_3 \sm_2\rb_1, \cdots ,  \lb\alpha \sm_n \sm_{n-1}\ldots \sm_3\sm_2\rb_1
\end{equation*}

\begin{prop}
	If $\lng(\alpha)=k+1$, at least one of the facets of $\lb\alpha\rb_0$ lies in $\C^n(k)$.
	\begin{proof}
		We must show that at least one facet of $\lb\alpha\rb_0$ is also a facet of some simplex~$\lb\alpha'\rb_0$, where $\lng(\alpha')\leq k$.
		
		Consider $\EndGen_1(\alpha)$. If this is non-empty then there exists~$\eta$ with length at least 1 in $A^+_1$ such that~$\alpha=\alpha'\eta$. It follows that~$\lb\alpha\rb_1=\lb\alpha'\eta\rb_1=\lb\alpha'\rb_1$. Therefore the facet $\lb\alpha\rb_1$ is also a facet of~$\lb\alpha'\rb_0$. Since~$\lng(\eta)\geq 1$ then~$\lng(\alpha')<\lng(\alpha)=k+1$, so~$\lb\alpha'\rb_0$ is in $\C^n(k)$.
		
		Alternatively if~$\EndGen_1(\alpha)=\emptyset$, then~$\lng(\alpha)\geq 1$ implies that~$\EndGen_n(\alpha) \neq \emptyset$. It follows from these two observations that~$\{\sm_2, \ldots \sm_n\}\cap \EndGen_n(\alpha)\neq \emptyset$ i.e.~for some $2\leq j \leq n$,~$\alpha=\alpha'\sm_j$. Applying the face map~$\cd^n_{j-2}$ gives
		\begin{eqnarray*}
			\cd^{n-1}_{j-2}(\lb\alpha\rb_0) &=& \lb\alpha\sm_{j-1}\ldots \sm_{2}\rb_{1}\\
			&=&\lb\alpha' \sm_j \sm_{j-1}\ldots \sm_{2}\rb_{1}\\
			&=& \cd^{n-1}_{j-1}(\lb\alpha'\rb_0)
		\end{eqnarray*}
		and as before $\lng(\alpha')\leq k$. This shows that the facet~$\cd^{n-1}_{j-2}(\lb\alpha\rb_0)$ is also a facet of~$\lb\alpha'\rb_0$ and is therefore in $\C^n(k)$.
	\end{proof}	
\end{prop}

To complete the proof of point (\ref{point A}) we must show that if a lower dimensional face of $\lb\alpha\rb_0$ is contained in $\C^n(k)$ then it is contained in a facet of $\lb\alpha\rb_0$, which is itself contained in~$\C^n(k)$. We first describe a general form for faces of $\lb\alpha\rb_0$.

\subsection{Proof of point (\ref{point A}): low dimensional faces of $\lb\alpha\rb_0$}

\begin{defn} \label{def:facets}
	A face of $\lb\alpha\rb_0$ is obtained by applying a series of face maps to $\lb\alpha\rb_0$. We denote the series of face maps applied by a tuple $(\cd^{n-1}_{i_2}, \cd^{n-2}_{i_3},\ldots, \cd^{n-r+1}_{i_r})$, and we let $a_j := \sm_{i_j-1+j}\ldots \sm_j$. That is, the $(j-1)$st map in the tuple corresponds to right multiplication by $a_j$. We note here that $a_j$ has length $i_j$ and ends with the generator $\sm_j$, unless $i_j=0$ in which case $a_j=e$.
	\begin{eqnarray*}
		\cd^{n-j+1}_{i_j}: \C^n_{n-j+1} &\to& \C^n_{n-j} \\
		: \lb\alpha\rb_{j-2} &\mapsto& \lb\alpha \sm_{i_j-1+j}\ldots \sm_j \rb_{j-1}\\
		&=&\lb\alpha a_j \rb_{j-1}.
	\end{eqnarray*}
	
\end{defn}

\noindent From now on we assume that the first map in a tuple maps from $\C^n_{n-1}$ to $\C^n_{n-2}$, the second map from $\C^n_{n-2}$ to $\C^n_{n-3}$ and so on. We therefore dispense of the superscripts in the $\cd$ notation for the face maps when we write these tuples.

With the above notation, an $(n-p-1)$ subsimplex of $\lb\alpha\rb_0$ occurs when a tuple of face maps $(\cd_{i_2}, \cd_{i_3},\ldots, \cd_{i_{p+1}})$ is applied to $\lb\alpha\rb_0$. The image of these maps is then the subsimplex $\lb\alpha a_2 \ldots a_{p+1}\rb_p$ with $a_j$ defined as in Definition \ref{def:facets} above.

\begin{lem}\label{lem:ijs}
	With the above notation, the tuple of face maps $(\cd_{i_j})_{j=2}^{p+1}$ can be written such that $i_{j+1}\geq i_j$ for all $j$, which translates to $\lng(a_{j+1})\geq \lng(a_j)$.
	\begin{proof}
		This is a direct consequence of Lemma \ref{lem:ijs background}.
	\end{proof}
\end{lem}

\begin{lem} \label{lem:biggerfaces}
	The $(n-p-1)$ subsimplex of $\lb\alpha\rb_0$ given by $(\cd_{i_2}, \cd_{i_3},\ldots, \cd_{i_{p+1}})$, or alternatively $\lb\alpha a_2 \ldots a_{p+1}\rb_p$, is a subsimplex of the following facets of $\lb\alpha\rb_0$:
	\begin{itemize}
		\item $\cd_{i_2}(\lb\alpha\rb_0)=\lb\alpha a_2\rb_1$
		\item $\cd_{i_3+1}(\lb\alpha\rb_0)=\lb\alpha a_3 \sm_2\rb_1$		
		\item $\cd_{i_4+2}(\lb\alpha\rb_0)=\lb\alpha a_4 \sm_3 \sm_2\rb_1$\\
		$\cdots$
		\item $\cd_{i_{p+1}+p-1}(\lb\alpha\rb_0)=\lb\alpha a_{p+1}\sm_p\ldots \sm_2\rb_1$
	\end{itemize}
	In general these facets are given by the face map $$\cd_{i_j+(j-2)}:\lb\alpha\rb_0 \mapsto \lb\alpha a_j \sm_{j-1} \ldots \sm_2\rb_1.$$
	\begin{proof}
		%$p$ face maps in the tuple  $(\cd_{i_2}, \cd_{i_3},\ldots, \cd_{i_{p+1}})$ are applied to $\lb\alpha\rb_0$ to get the $(n-p-1)$ subsimplex  $\lb\alpha a_2 \ldots a_{p+1}\rb_p$. The first of these face maps takes $\lb\alpha\rb_0$ to a face of  $\lb\alpha\rb_0$ containing $\lb\alpha a_2 \ldots a_{p+1}\rb$. Therefore
		It is enough to show that $\cd_{i_j+(j-2)}$ can act as the first face map in the tuple~$(\cd_{i_2}, \cd_{i_3},\ldots, \cd_{i_{p+1}})$ for all~$j$.
		Recall from Lemma \ref{lem:ijs} that in the tuple~$i_{j+1}\geq~i_j$ for all $j$. Using the simplicial identities, the tuple is equivalent to the tuple $(\cd_{i_j+(j-2)},\cd_{i_2}, \cd_{i_3},\ldots, \widehat{\cd_{i_j}}, \ldots,  \cd_{i_{p+1}})$.
	\end{proof}
\end{lem}

For the remainder of this section, let $\alpha$ in $A_n^+$ with $\lng(\alpha)=k+1$. The aim of this section is to show that if the $(n-p-1)$ subsimplex of $\lb\alpha\rb_0$ given by~$(\cd_{i_2}, \cd_{i_3},\ldots, \cd_{i_{p+1}})$ or alternatively $\lb\alpha a_2 \ldots a_{p+1}\rb_p$ is in $\C^n(k)$ then it follows that one of the facets of~$\lb\alpha\rb_0$ from Lemma \ref{lem:biggerfaces} is also in $\C^n(k)$. The proof of point (\ref{point A}) will follow.

\begin{defn} \label{defn:bjs}
	If $\lb\alpha a_2 \ldots a_{p+1}\rb_p$ is in $\C^n(k)$ then it is also a $(n-p-1)$ subsimplex of a simplex $\lb\beta\rb_0$ for some $\beta$ in $A^+_n$ such that $\lng(\beta)\leq k$. The subsimplex is therefore obtained from $\lb\beta\rb_0$ by applying a tuple of face maps, denote these $(\cd_{l_2}, \cd_{l_3},\ldots, \cd_{l_{p+1}})$ and order as in Lemma \ref{lem:ijs} such that $l_{j+1}\geq l_j$ for all $j$. Define $b_j:=\sm_{l_j-1+j}\ldots \sm_j$ and when $l_j=0$ let $b_j=e$. Then $(\cd_{l_2}, \cd_{l_3},\ldots, \cd_{l_{p+1}})$ applied to $\lb\beta\rb_0$ gives the $(n-p-1)$ simplex $\lb\beta b_2\ldots b_{p+1}\rb_p$. By construction $\lb\beta b_2\ldots  b_{p+1}\rb_p=\lb\alpha a_2 \ldots a_{p+1}\rb_p$. We recall here that $\lng(a_j)=i_j$ and $\lng (b_j)=l_j$.
\end{defn}

\begin{lem}\label{lem:ljs minimal}
	Choose $\beta$ and $b_j$ as defined above, such that $\sum_{k=2}^{p+1} l_k$ is minimal, corresponding to  $b_2\ldots b_{p+1}$ being of minimal length. This choice of  $b_2\ldots b_{p+1}$ then corresponds to either:
	$$
	\lb\alpha a_2\ldots a_{p+1}\rb_p=\lb\beta\rb_p \text{ that is, } l_j=0 \, \forall \, j
	$$
	or
	$$
	\lng(\beta)=\lng(\alpha)-1=k.
	$$
	\begin{proof}
		Suppose that $\beta$ and $b_j$ are chosen such that $\sum_{k=2}^{p+1} l_k$ is minimal, and furthermore suppose that $\lng(\beta)<\lng(\alpha)-1$ and $\sum_{k=2}^{p+1} l_k > 0$. Then some $l_k \neq 0$: set $j$ to be minimal such that $l_j \neq 0$. Then $b_j =\sm_{l_j-1+j}\ldots \sm_j \neq e$ and $$\lb\beta b_2\ldots b_{p+1}\rb_p=\lb\beta b_j\ldots b_{p+1}\rb_p=\lb\beta \sm_{l_j-1+j}\ldots \sm_j b_{j+1}\ldots b_{p+1}\rb_p.$$ But this is the tuple of face maps $(\cd_{l_{j}-1}, \cd_{l_{j+1}},\ldots, \cd_{l_{p+1}})$ applied to $\lb\beta \sm_{l_j -1 + j}\rb_0$. Since $\lng(\beta)<\lng(\alpha)-1$ it follows that $\lng(\beta \sm_{l_j -1 + j})\leq \lng(\alpha)-1$ and so $\lb\beta \sm_{l_j -1 + j}\rb_0$ is in~$\C^n(k)$. However the tuple for $\beta \sm_{l_j -1 + j}$ has the sum of its corresponding $l_j$ less than the original tuple for $\beta$. This is a contradiction, as $\beta$ was chosen to have minimal~$\sum_{k=2}^{p+1} l_k$. Therefore either  $\sum_{k=2}^{p+1} l_k=0 $, or~$\lng(\beta)=\lng(\alpha)-1$.
	\end{proof}
\end{lem}

For the remainder of this paper, assume $\beta$ and $b_j$ are chosen such that $\sum_{k=2}^{p+1} l_k$ is minimal, so we have
$$
\lb\beta b_2\ldots  b_{p+1}\rb_p=\lb\alpha a_2 \ldots a_{p+1}\rb_p
$$
for either $\sum_{k=2}^{p+1} l_k=0$ or $\lng(\beta)=\lng(\alpha)-1=k$.
We use the following notation throughout the remainder of this paper.
\begin{defn} \label{def:a b}
	Let $a:=a_2 \ldots a_{p+1}$ and $b:=b_2\ldots  b_{p+1}$. Note that $\sum_{k=2}^{p+1} l_k=0$ corresponds to $b=e$. So we have
	$$
	\lb\alpha a \rb_p=\lb\beta b\rb_p
	$$
	Where either~$b=e$, or~$\lng(\beta)=\lng(\alpha)-1$. We recall that this is equivalent to $\overline{\alpha a}=\overline{\beta b}$ in $A^+(n;p)$. Let $\gamma:=\overline{\alpha a}=\overline{\beta b}$, and define $u$ and $v$ in $A^+_p$ such that
	$$
	\alpha a=\gamma u \text{  and  } \beta b=\gamma v.
	$$
\end{defn}

Recall point \ref{point A}: if $\lng(\alpha)=k+1$ then $\lb\alpha\rb_0 \cap \C^n(k)$ is a non-empty union of top dimensional faces of $\lb\alpha\rb_0$. Recall that we have fixed a face~$\lb\alpha a \rb_p$ of~$\lb\alpha \rb_0$, and we wish to prove this is contained in a facet of~$\lb\alpha\rb_0$ (from Lemma~\ref{lem:biggerfaces}) which is contained in~$\C^n(k)$. We complete the proof of this by splitting into three cases:
\begin{center}
	\begin{enumerate}[(i)]
		\item \label{casei} $\lng(\beta b)<\lng(\alpha a)$
		\item \label{caseii} $\lng(\beta b)=\lng(\alpha a)$
		\item \label{caseiii} $\lng(\beta b)>\lng(\alpha a)$
	\end{enumerate}
\end{center}

\noindent and since multiplication in the Artin monoid corresponds to adding lengths the conditions of these cases correspond to analogous conditions on the lengths of $u$ and $v$.

\begin{rem}
	Note that if $\sum_{k=2}^{p+1} l_k=0$ then $b=e$, and since $\lng(\beta)< \lng(\alpha)$ it follows we are therefore in case (\ref{casei}):~$\lng(\beta b)<\lng(\alpha a)$.
\end{rem}

We prove the three cases one by one in the following subsections. This involves some technical lemmas, and in particular computation of least common multiples of strings of words. We therefore include these technical lemmas in a separate section and refer to them as required.

\subsection{Proof of point (\ref{point A}): preliminary lemmas}

Recall from Definition~\ref{def:facets} that a face of $\lb\alpha\rb_0$ is obtained by applying a series of face maps to~$\lb\alpha\rb_0$. We denote the series of face maps by a tuple~$(\cd^{n-1}_{i_2}, \cd^{n-2}_{i_3},\ldots, \cd^{n-r+1}_{i_r})$, and we let~$a_j = \sm_{i_j-1+j}\ldots \sm_j$ and when $i_j=0$ let $a_j=e$. That is, the $(j-1)$st map in the tuple corresponds to right multiplication by $a_j$. We let~$a=a_{2}\ldots a_{p+1}$. Recall also that if $\lb \alpha a \rb_p$ is in~$\C^n(k)$ then the subsimplex is also obtained from some~$\lb\beta\rb_0$ for~$\lng(\beta)\leq k$, by applying a tuple of face maps~$(\cd_{l_2}, \cd_{l_3},\ldots, \cd_{l_{p+1}})$. Recall~$b_j:=\sm_{l_j-1+j}\ldots \sm_j$ and when $l_j=0$ let $b_j=e$. Let~$b=b_2\ldots b_{p+1}$. By construction~$\lb\beta b\rb_p=\lb\alpha a\rb_p$. Recall from Definition~\ref{def:lcm} that for~$\alpha$ and~$\beta$ two words in~$A^+$, we denote the least common multiple of~$\alpha$ and~$\beta$ (if it exists) by~$\Dl(\alpha, \beta)$.

\begin{lem}\label{lem:sigma manipulation}
	For all~$k> j$, the generators~$\sm_i$ satisfy
	$$
	(\sm_k\ldots\sm_{j+1})\sm_j(\sm_k\ldots\sm_{j+1})=(\sm_{k-1}\sm_k)(\sm_{k-2}\sm_{k-1})\ldots(\sm_{j+1}\sm_{j+2})(\sm_j\sm_{j+1}\sm_j).
	$$
	\begin{proof}
		We proceed by induction on~$k-j$. For the base case let~$k-j=1$ i.e.~$k=j+1$. Then the left hand side of the above equation evaluates to~$\sm_{j+1}\sm_j\sm_{j+1}$ and the right hand side evaluates to~$\sm_j\sm_{j+1}\sm_j$. These are equal by the Artin relations.
		For the inductive hypothesis we assume the Lemma is true for~$k-j<r$, and we prove for~$k-j=r$, i.e.~$k=j+r$. We manipulate the left hand side of the equation, and show equality to the right hand side:
		\begin{eqnarray*}
			&&(\sm_k\ldots\sm_{j+1})\sm_j(\sm_k\ldots\sm_{j+1})\\&=& (\sm_{j+r}\ldots\sm_{j+1})\sm_j({\color{red}\sm_{j+r}}\ldots\sm_{j+1})\\
			&=&(\sm_{j+r}\sm_{j+r-1}{\color{red}\sm_{j+r}}\ldots\sm_{j+1})\sm_j(\sm_{j+r-1}\ldots\sm_{j+1})\\
			&=&(\sm_{j+r}\sm_{j+r-1}\sm_{j+r})(\sm_{j+r-2}\ldots\sm_{j+1})\sm_j(\sm_{j+r-1}\ldots\sm_{j+1})\\
			&=&(\sm_{j+r-1}\sm_{j+r}\sm_{j+r-1})(\sm_{j+r-2}\ldots\sm_{j+1})\sm_j(\sm_{j+r-1}\ldots\sm_{j+1})\\
			&=&(\sm_{j+r-1}\sm_{j+r})(\sm_{j+r-1}\sm_{j+r-2}\ldots\sm_{j+1})\sm_j(\sm_{j+r-1}\ldots\sm_{j+1})\\
			&=&(\sm_{j+r-1}\sm_{j+r})(\sm_{j+r-2}\sm_{j+r-1})\ldots(\sm_{j+1}\sm_{j+2})(\sm_j\sm_{j+1}\sm_j)
		\end{eqnarray*} 
		where the final equality applies the inductive hypothesis.
	\end{proof}
\end{lem}

\begin{lem}\label{lem: lcm a_j+1 s_j}
	With notation as above, $\Dl(a_{j+1}, \sm_j)=a_{j+1}\sm_ja_{j+1}$ (and similarly~$\Dl(b_{j+1}, \sm_j)=b_{j+1}\sm_jb_{j+1}$).
	\begin{proof}
		The proof is the same for both the~$a_j$ and~$b_j$ case, so we prove it for the~$a_j$ case. We must show
		\begin{enumerate}[(a)]
			\item $a_{j+1} \preceq_R a_{j+1}\sm_ja_{j+1}$ and $\sm_j \preceq_R a_{j+1}\sm_ja_{j+1}$.
			\item if $x$ in $A^+_n$ is a common multiple of $a_{j+1}$ and $\sm_j$, then $a_{j+1}\sm_ja_{j+1} \preceq_R x$.
		\end{enumerate}	
		
		Recall~$a_{j+1} := \sm_{i_{j+1}+j}\ldots \sm_{j+1}$. Without loss of generality, we relabel~$j=1$ and~$i_{j+1}+j=k$. Then~$a_{j+1}=\sm_k\ldots \sm_2$ and~$\sm_j=\sm_1$.
		
		To prove (a) note that~$a_{j+1} \preceq_R a_{j+1}\sm_ja_{j+1}$ by observation, and also
		\begin{eqnarray*}
			a_{j+1}\sm_ja_{j+1}&=&(\sm_{k}\ldots \sm_2)\sm_1({\color{black}\sm_{k}}\ldots \sm_{2})\\
			&=&(\sm_{k-1}\sm_{k}\sm_{k-2}\sm_{k-1}\ldots \sm_{2}\sm_{3})(\sm_{1}\sm_{2}\sm_{1})
		\end{eqnarray*}
		by Lemma~\ref{lem:sigma manipulation}, so $\sm_1=\sm_j \preceq_R a_{j+1}\sm_ja_{j+1}$.
		
		To prove (b) we show by induction on $\lng(a_{j+1})$ that any common multiple $x$ must satisfy $a_{j+1}\sm_ja_{j+1} \preceq_R x$. When $\lng(a_{j+1})=1$, $a_{j+1}=\sm_{2}$ and we have
		$$\Dl(\sm_{2},\sm_1)=\sm_{2}\sm_1\sm_{2}=a_{j+1}\sm_ja_{j+1}.$$ For $\lng(a_{j+1})=r-1$ when $r\geq 2$, assume that $\Dl(a_{j+1}, \sm_j)=a_{j+1}\sm_ja_{j+1}$ and prove for $\lng(a_{j+1})=r$. Assume $x$ satisfies $a_{j+1}\preceq_R x$ and $\sm_j \preceq_R x$. Since $\lng(a_{j+1})=r$, $a_{j+1}=\sm_{r+1}\ldots \sm_{2}$ and so $\sm_{r+1}\ldots \sm_{2} \preceq_R x$ which in particular gives $\sm_{r}\ldots \sm_{2} \preceq_R x$. By the inductive hypothesis it follows that $$\Dl (\sm_{r}\ldots \sm_{2}, \sm_1)=(\sm_{r}\ldots \sm_{2})\sm_1(\sm_{r}\ldots \sm_{2}).$$ and this is in $\EndMon_n(x)$ by Lemma \ref{lem:deltaendsets}. Let $x=x'(\sm_{r}\ldots \sm_{2})\sm_1(\sm_{r}\ldots \sm_{2})$. Then since $\sm_{r+1}\ldots \sm_{2} \preceq_R x$, by cancellation of $\sm_{r}\ldots \sm_{2}$ it follows that $$\sm_{r+1}\preceq_R x'(\sm_{r}\ldots \sm_{2})\sm_1=x'\sm_{r}(\sm_{r-1}\ldots \sm_{2}\sm_1).$$ Since $\sm_{r+1}$ letterwise commutes with $(\sm_{r-1}\ldots \sm_{2}\sm_1)$, from Lemma \ref{lem:mini} we have $\sm_{r+1}\preceq_R x'\sm_{r}$. From Lemma \ref{lem:deltaendsets} it follows $\Dl (\sm_{r+1},\sm_{r})=\sm_{r}\sm_{r+1}\sm_{r} \preceq_R x'\sm_{r}$. By cancellation of $\sm_{r}$ this gives $x'=x''\sm_{r}\sm_{r+1}$, so
		\begin{eqnarray*}
			x&=&(x')(\sm_{r}\ldots \sm_{2})\sm_1(\sm_{r}\ldots \sm_{2})\\
			&=& (x''\sm_{r}\sm_{r+1})(\sm_{r}\ldots \sm_{2})\sm_1(\sm_{r}\ldots \sm_{2})\\
			&=& x''(\sm_{r}\sm_{r+1}\sm_{r})(\sm_{r-1}\ldots \sm_{2})\sm_1(\sm_{r}\ldots \sm_{2})\\
			&=& x''(\sm_{r+1}\sm_{r}\sm_{r+1})(\sm_{r-1}\ldots \sm_{2})\sm_1(\sm_{r}\ldots \sm_{2})\\
			&=& x''(\sm_{r+1}\sm_{r}{\color{red}\sm_{r+1}}\sm_{r-1}\ldots \sm_{2})\sm_1(\sm_{r}\ldots \sm_{2})\\
			&=& x''(\sm_{r+1}\sm_{r}\sm_{r-1}\ldots \sm_{2})\sm_1({\color{red}\sm_{r+1}}\sm_{r}\ldots \sm_{2})\\
			&=&x'' a_{j+1} \sm_j a_{j+1}
		\end{eqnarray*}
		as required.
	\end{proof}
\end{lem}

\begin{lem}\label{lem:hat}
	Recall from Lemma \ref{lem: lcm a_j+1 s_j} that $\Dl(a_{j+1}, \sm_j)=a_{j+1}\sm_ja_{j+1}$. Then when~$j\geq 2$ this expression satisfies
	$$
	a_{j+1}\sm_ja_{j+1}=\hat{a}_ja_ja_{j+1}\sm_j
	$$
	where $\hat{a}_j=\sm_{i_{j+1}+j-1}\ldots \sm_{i_{j}+j}$ and letterwise commutes with $a_2\ldots a_{j-1}$.
	
	When~$j=1$, the expression satisfies
	$$
	a_2\sm_1a_2=\hat{a}_1\sm_1a_2\sm_1
	$$
	i.e.~the same equality holds, setting~$a_1:=\sm_1$. The analogous statements hold for the~$b_j$.
	\begin{proof}
		Recall $a_{j+1} = \sm_{i_{j+1}+j}\ldots \sm_{j+1}$ and $a_j = \sm_{i_j-1+j}\ldots \sm_j$. Without loss of generality, relabel $j=1$, $i_{j+1}+j=k$, and $i_{j}-1+j=l$. Then $a_{j+1}=\sm_k\ldots \sm_2$ , $\sm_j=\sm_1$, and $a_j=\sm_l\ldots \sm_1$. Note that since $i_{j+1}\geq i_j$ then $k>l$. We wish to show that
		$a_{j+1}\sm_ja_{j+1}=\hat{a}_ja_ja_{j+1}\sm_j$ where $\hat{a}_j=\sm_{k-1}\ldots \sm_{l+1}$.
		We proceed by induction on the length of~$a_{j+1}$. For the base case, when~$\lng(a_{j+1})=1$ this implies that~$a_{j+1}=\sm_2$. Then we have
		\begin{equation*}
			a_{j+1}\sm_ja_{j+1}=\sm_2\sm_1\sm_2=\sm_1\sm_2\sm_1.
		\end{equation*}
		Since~$0\leq i_j\leq i_{j+1}=1$ there are now two options. In the case~$i_j=1$ the right hand side is~$a_ja_{j+1}\sm_j$ and~$\hat{a}_j=e$. In the case~$i_j=0$ then~$a_j=e$ and the right hand side is $\hat{a}a_{j+1}\sm_j$ with~$\hat{a}_j=\sm_1=\sm_{i_j+j}$.
		
		For the inductive hypothesis we assume true for~$\lng(a_{j+1})\leq{r-1}$ and prove for $\lng(a_{j+1})={r-1}$, i.e.~$k=r$.
		Recall from Lemma~\ref{lem:sigma manipulation} that
		\begin{eqnarray*}
			a_{j+1}\sm_ja_{j+1}&=&(\sm_{k}\ldots \sm_2)\sm_1(\sm_{k}\ldots \sm_{2})\\
			&=&(\sm_{r}\ldots \sm_2)\sm_1(\sm_{r}\ldots \sm_{2})\\
			&=&(\sm_{r-1}\sm_{r})(\sm_{r-2}\sm_{r-1})\ldots (\sm_{2}\sm_{3})(\sm_{1}\sm_{2}\sm_{1})\\
			&=&(\sm_{r-1}{\sm_{r}})(\sm_{r-2}\ldots\sm_{l+1})a_j(\sm_{r-1}\ldots\sm_2)\sm_j
		\end{eqnarray*}
		where the final equality applies the inductive hypothesis. Then~$\sm_{r}$ commutes with~$(\sm_{r-2}\ldots\sm_{l+1})a_j$ since~$\lng(a_j)\leq \lng(a_{j+1})$. This gives the following:
		\begin{eqnarray*}
			a_{j+1}\sm_ja_{j+1}
			&=&(\sm_{r-1}{\sm_{r}})(\sm_{r-2}\ldots\sm_{l+1})a_j(\sm_{r-1}\ldots\sm_2)\sm_j\\
			&=&(\sm_{r-1}\sm_{r-2}\ldots\sm_{l+1})a_j({\sm_{r}}\sm_{r-1}\ldots\sm_2)\sm_j\\
			&=& \hat{a}_ja_ja_{j+1}\sm_j
		\end{eqnarray*}		
	 Since $i_j\geq i_{j-1}$ it follows that $l-1$ is the maximal index of a generator appearing in $a_{j-1}$ and hence in the string $a_2\ldots a_{j-1}$. Therefore $\hat{a}_j$ letterwise commutes with $a_2\ldots a_{j-1}$ since the indices of the generators in each word pairwise differ by at least two.
		
		Since the~$b_j$ have the same form as the~$a_j$, with difference only in word length, the analogous statements hold for the~$b_j$. 
	\end{proof}
\end{lem}

	Recall the definition of $a_j$ and $b_j$ for $2\leq j \leq p+1$, from Definition \ref{def:facets} and Definition \ref{defn:bjs} respectively.
	
\begin{defn}\label{def:c}
For $2\leq j \leq p+1$ define $c_j$ as follows
	$$
	c_j =
	\left\{
	\begin{array}{ll}
	a_j  & \mbox{if } \lng(a_j)\geq \lng(b_j) \\
	b_j & \mbox{if } \lng(a_j)< \lng(b_j)
	\end{array}
	\right.
	$$
	for $2\leq j \leq p+1$. Define $c:=c_2\ldots c_{p+1}$.
	Let 	
	$$a_j' =
	\left\{
	{\begin{array}{ll}
		e  & \mbox{if } \lng(a_j)\geq \lng(b_j) \\
		\sm_{l_j+j-1}\ldots \sm_{i_j+j} & \mbox{if } \lng(a_j)< \lng(b_j)
	\end{array} }\right.
	$$
	and similarly
	$$
		b_j' =
	\left\{
	{\begin{array}{ll}
	e  & \mbox{if } \lng(b_j)\geq \lng(a_j) \\
	\sm_{i_j+j-1}\ldots \sm_{l_j+j} & \mbox{if } \lng(b_j)< \lng(a_j).
	\end{array} } \right.
	$$
	Define~$a'=a_2'\ldots a_{p+1}'$ and~$b'=b_2'\ldots b_{p+1}'$.
\end{defn}

\begin{lem} \label{lem:c is lcm}
	With $c, a'$ and~$b'$ as in Definition \ref{def:c} and $a$ and $b$ as defined in Definition \ref{def:a b} we have $c=\Dl(a,b)$ and in particular $c=a'a=b'b$.
	\begin{proof}
		We prove that
		\begin{enumerate}[(a)]
			\item $c=a'a=b'b$
			\item if $x$ in $A^+_n$ is a common multiple of $a$ and $b$, then $c \preceq_R x$.
		\end{enumerate}
		To prove (a), we show that $c=a'a$: the proof that $c=b'b$ is symmetric. It follows from the definitions that
		$c_j=a_j'a_j$.
		The smallest generator index in~$a_j'$ is~$(i_{j}+j)$ and the largest generator index in~$a_2\ldots a_{j-1}$ is~$(i_{j-1}+(j-1)-1)$. The elements~$a_j'$ and~$a_2\ldots a_{j-1}$ letterwise commute, since~$i_j \geq i_{j-1}$ so
		$$\mid (i_{j}+j) - (i_{j-1}+(j-1)-1)\mid =\mid (i_{j} - i_{j-1})+2)\mid \geq 2.$$ Let $a'=a_2'\ldots a_{p+1}'$. Then we compute
		\begin{eqnarray*}
			c&=&c_2\ldots c_{p+1}\\
			&=& (a_2'a_2)({\color{red}a_3'}a_3)\ldots (a_{p+1}'a_{p+1})\\
			&=& a_2'{\color{red}a_3'}a_2a_3\ldots( a_{p+1}'a_{p+1})\\
			&=& a_2'a_3'\ldots a_{p+1}'a_2a_3\ldots a_{p+1}\\
			&=& (a_2'a_3'\ldots a_{p+1}')(a_2a_3\ldots a_{p+1})\\
			&=&a'a
		\end{eqnarray*}
		which completes the proof of (a).
		
		To prove (b) assume $x$ is a common multiple of $a$ and $b$.
		
		\noindent{\bf Claim:} If $c_k\ldots c_{p+1}\preceq_R x$ for some $2 \leq k \leq p+1$ then $x=x_kc_k\ldots c_{p+1}$ for some~$x_k$ in $A^+_n$. We claim that $x_k$ satisfies $a_2\ldots a_{k-1} \preceq_R x_k$ and $b_2\ldots b_{k-1} \preceq_R x_k$.
		
		Given the claim, the proof of (b) will follow since~$a=(a_2\ldots a_{p+1})\preceq_R x$ and~$b=(b_2\ldots b_{p+1}) \preceq_R x$ implies that~$c_{p+1}\preceq_R x$, so~$x=x_{p+1}c_{p+1}$. But then~$x_{p+1}$ satisfies~$a_2\ldots a_{p} \preceq_R x_{p+1}$ and~$b_2\ldots b_{p} \preceq_R x_{p+1}$ by the claim for~$k=p+1$. In particular~$c_p\preceq_R x_{p+1}$ and it follows that~$x=x_{p}c_pc_{p+1}$. Continuing in this manner we arrive at~$x=x_2(c_2\ldots c_{p+1})=x_2c$ and so~$c\preceq_R x$. It therefore remains to prove the claim.
		
		Since $c_k\ldots c_{p+1}=(a_k'a_k)\ldots (a_{p+1}'a_{p+1})=(a_k'\ldots a_{p+1}')(a_k \ldots a_{p+1})$ it follows that
		\begin{eqnarray*}
			x&=&x_k(c_k\ldots c_{p+1})\\
			&=&x_k (a_k'\ldots a_{p+1}')(a_k \ldots a_{p+1})\\
			&=&y_k(a_k \ldots a_{p+1}) \text{  for } y_k=x_k(a_k'\ldots a_{p+1}').
		\end{eqnarray*}
		Since $x$ is a common multiple of $a$ and $b$ then we also have $a=(a_2 \ldots a_{p+1}) \preceq_R x$, i.e.~for some $z_k$.
		$$
		x=z_k(a_2 \ldots a_{p+1})
		$$
		Therefore by cancellation of $(a_k \ldots a_{p+1})$,
		$$
		y_k=z_k(a_2 \ldots a_{k-1})
		$$
		By Lemma \ref{lem:LWCends}, $\Dl((a_k'\ldots a_{p+1}'), (a_2\ldots a_{k-1} )) \preceq_R y_k$. Since the two words letterwise commute $\Dl((a_k'\ldots a_{p+1}'), (a_2\ldots a_{k-1}) )=(a_2\ldots a_{k-1})(a_k'\ldots a_{p+1}')$ and so $$y_k=w_k(a_2\ldots a_{k-1})(a_k'\ldots a_{p+1}')$$ for some $w_k$ in $A_n^+$. It follows
		\begin{eqnarray*}
			x&=&x_k(c_k\ldots c_{p+1})\\
			&=&y_k(a_k \ldots a_{p+1})\\
			&=&w_k(a_2\ldots a_{k-1})(a_k'\ldots a_{p+1}')(a_k \ldots a_{p+1})\\
			&=&w_k(a_2\ldots a_{k-1})((a_k'\ldots a_{p+1}')(a_k \ldots a_{p+1}))\\
			&=&w_k(a_2\ldots a_{k-1})(c_k\ldots c_{p+1})
		\end{eqnarray*}
		and by cancellation of $c_k\ldots c_{p+1}$ on the first and final lines of the above equation, $(a_2\ldots a_{k-1})\preceq_R x_k$ as required. The proof for $(b_2\ldots b_{k-1})\preceq_R x_k$ is identical. This completes the proof of the Claim and thus of (b).
	\end{proof}
\end{lem}

	Recall that we have fixed $\alpha$ in $A_n^+$ with $\lng(\alpha)=k+1$, and we have fixed a face $\lb\alpha a_2 \ldots a_{p+1}\rb_p$ in $\C^n(k)$. We want to show that one of the facets of~$\lb\alpha\rb_0$ from Lemma \ref{lem:biggerfaces} is also in $\C^n(k)$. Recall since $\lb\alpha a\rb_p$ is in $\C^n(k)$, there exists $\beta$ in $A^+_n$ such that $\lng(\beta)\leq k$ and
	$$
	\lb\alpha a \rb_p=\lb\alpha a_2 \ldots a_{p+1}\rb_p=\lb\beta b_2\ldots  b_{p+1}\rb_p=\lb\beta b\rb_p
	$$ from Definition~\ref{def:a b}, where $\lng(a_j)=i_j$ and $\lng (b_j)=l_j$. We have assumed $\beta$ and $b_j$ are chosen such that $\sum_{k=2}^{p+1} l_k$ is minimal, so we have either~$b=e$, or~$\lng(\beta)=\lng(\alpha)-1$. Recall $\gamma:=\overline{\alpha a}=\overline{\beta b}$, and that we defined $u$ and $v$ in $A^+_p$ such that
		$$
		\alpha a=\gamma u \text{  and  } \beta b=\gamma v.
		$$

We prove in the next three Lemmas that in the case $\EndGen_p(\alpha a)\neq \emptyset$ we are done.

\begin{lem}\label{lem:hypothesis1}
	If $\EndGen_0(\alpha a) \neq \emptyset$ then the facet $\lb\alpha a_2\rb_1$ containing $\lb\alpha a\rb_p$ is in~$\C^n(k)$.
	\begin{proof}
		Consider $\tau$ in $\EndGen_0(\alpha a)$. Then since the generators $S_0$ of $A^+_0$ commute with $\sm_2, \ldots, \sm_n$ it follows that $\tau$ letterwise commutes (Definition \ref{def:LW}) with~$a$, because $a=a_2\ldots a_{p+1}$ only contains generators in the set of $\{\sm_2, \ldots \sm_n\}$. Since~$\tau$ and~$a$ are both in $\EndMon_n(\alpha a)$ and they letterwise commute, it follows from Lemma~\ref{lem:mini} that $\tau$ is in $\EndMon_n(\alpha)$ i.e.~some $\alpha'$ in $A^+_n$ with~$\lng(\alpha')<\lng(\alpha)$ satisfies~$\alpha=\alpha'\tau$ .
		
		The facet $\lb\alpha a_2\rb_1$ therefore satisfies
		$$
		\lb\alpha a_2\rb_1=\lb\alpha' \tau a_2\rb_1=\lb\alpha' a_2 \tau\rb_1=\lb\alpha' a_2\rb_1.
		$$
		Here the final equality is due to $\overline{\alpha' a_2 \tau}=\overline{\alpha' a_2}$ where the reduction is taken with respect to $A^+_1$ (from Lemma \ref{lemma:barmult}). The penultimate equality is due to the fact $\tau$ and $a_2$ letterwise commute. Since $\lng(\alpha')<\lng(\alpha)$, $\lb\alpha'\rb_0$ is in $\C^n(k)$ and $\lb\alpha' a_2\rb_1$ is a facet of $\lb\alpha'\rb_0$. Therefore $\lb\alpha a_2\rb_1$ is in $\C^n(k)$ and this completes the proof.
	\end{proof}
\end{lem}

The case where $\EndGen_p(\alpha a) \neq \emptyset$ but $\EndGen_0(\alpha a) = \emptyset$ requires the following technical lemma.

\begin{lem}\label{lem:a' LW coms and face map}
	Suppose $a_j \neq e$, then the words $a_j$ and $a_{j+1}$ as in Definition \ref{def:facets} satisfy $a_{j+1}\sm_j=\ba_ja_j$, for some $\ba_j$ in $A^+_n$ with $\lng(\ba_j)\geq 1$. Furthermore $\ba_j$ letterwise commutes with $a_2 \ldots a_{j-1}$. Regardless of whether or not $a_j=e$, $a_{j+1}\sm_j$ corresponds to the face map $\cd^{n-j+1}_{i_{j+1}+1}$. The analogous results hold for the~$b_j$.
	\begin{proof}
		If $a_j\neq e$ then
		\begin{eqnarray*}
			a_{j+1}\sm_j&=&(\sm_{i_{j+1}+j}\ldots \sm_{j+1})\sm_j\\
			&=&(\sm_{i_{j+1}+j}\ldots \sm_{i_j+j})(\sm_{i_j+j-1}\ldots \sm_{j+1})\sm_j \\
			&=&(\sm_{i_{j+1}+j}\ldots \sm_{i_j+j})(\sm_{i_j+j-1}\ldots \sm_{j+1}\sm_j) \\
			&=& (\sm_{i_{j+1}+j}\ldots \sm_{i_j+j})a_j\\
			&=&\ba_ja_j
		\end{eqnarray*}
		so $\ba_j=\sm_{i_{j+1}+j}\ldots \sm_{i_j+j}$, and~$\lng(\ba_j)\geq 1$ since $\lng( a_{j+1})\geq \lng(a_j)\geq 1$. The generators appearing in the word $a_2 \ldots a_{j-1}$ are~$\{\sm_2, \ldots, \sm_{i_{j-1}+(j-1)-1}\}$ and so to prove that~$\ba_j$ letterwise commutes with $a_2 \ldots a_{j-1}$ it is enough to show that the sets~$A=\{\sm_{i_{j}+j},\ldots, \sm_{i_{j+1}+j}\}$ and~$B=\{\sm_2, \ldots, \sm_{i_{j-1}+(j-1)-1}\}$ pairwise commute. The largest index of a generator in $B$ is~$i_{j-1}+(j-1)-1$ and the smallest index of a generator in~$A$ is~$i_{j}+j$ so it is enough to show~
		$$\mid (i_{j}+j) - (i_{j-1}+(j-1)-1)\mid=\mid (i_{j} - i_{j-1})+2\mid \geq 2.$$ 
		This holds since~$i_j\geq i_{j-1}$, and so~$\ba_j$ and~$a_2 \ldots a_{j-1}$ letterwise commute.
		Regardless of whether or not $a_j=e$, $a_{j+1}\sm_j=\ba_ja_j=\sm_{i_{j+1}+j}\ldots \sm_j$ corresponds to the face map~$\cd^{n-j+1}_{i_{j+1}+1}$ defined in Definition \ref{def:facets}.
		Since the~$b_j$ have the same form as the~$a_j$, with difference only in word length, the analogous statements hold for the~$b_j$.
	\end{proof}
\end{lem}

\begin{lem}\label{lem:hypothesis2}
	If $\EndGen_p(\alpha a) \neq \emptyset$ but $\EndGen_0(\alpha a) = \emptyset$ then some $\sm_j$ is in $\EndGen_p(\alpha a)$ for $1\leq j \leq p$. Then the facet $\lb\alpha a_j\sm_{j-1}\ldots \sm_2\rb_1$ containing $\lb\alpha a\rb_p$ is in $\C^n(k)$.
	\begin{proof}
		If $\EndGen_0(\alpha a)=\emptyset$ and $\EndGen_p(\alpha a) \neq \emptyset$ it follows that $$\{\sm_1, \sm_2, \ldots \sm_p\}\cap \EndGen_p(\alpha a)\neq \emptyset,$$ so some $\sm_j$ is in $\EndGen_p(\alpha a)$ for $1\leq j \leq p$. Then $\sm_j$ and $a=a_2\ldots a_{p+1}$ are both in $\EndMon_n(\alpha a)$. In particular $\sm_j$ and $a_{j+2}\ldots a_{p+1}$ are both in $\EndMon_n(\alpha a)$. Since~$\sm_j$ and~$a_{j+2}\ldots a_{p+1}$ letterwise commute we have from Lemma~\ref{lem:mini} that~$\sm_j$ is in $\EndMon_n(\alpha a_2 \ldots a_{j+1})$. Since $a_{j+1}$ is also in $\EndMon_n(\alpha a_2 \ldots a_{j+1})$, from Lemma~\ref{lem:deltaendsets}~$\Dl(a_{j+1}, \sm_j)$ is in $\EndMon_n(\alpha a_2 \ldots a_{j+1})$. Also, from Lemma~\ref{lem: lcm a_j+1 s_j}, $\Dl(a_{j+1}, \sm_j)=~a_{j+1}\sm_ja_{j+1}$. By cancellation of $a_{j+1}$ it follows that $a_{j+1}\sm_j$ is in $\EndMon_n(\alpha a_2 \ldots a_{j})$, so 		
		\begin{equation*}
			\alpha a_2 \ldots a_{j}=\alpha'(a_{j+1}\sm_j)\label{eq:star}\tag{$\dagger$}
		\end{equation*}
		for some $\alpha'$ in $A^+_n$.
				
		Recall Lemma \ref{lem:a' LW coms and face map} and split into two cases: either
		\begin{enumerate}[(a)]
			\item $a_j \neq e$, or
			\item $a_2=\cdots =a_j=e$ since~$\lng(a_i)\leq \lng(a_{i+1}) \forall i$.
		\end{enumerate}
		For case (a) recall from Lemma~\ref{lem:a' LW coms and face map} that $a_{j+1}\sm_j=\ba_ja_j$ and $\ba_j$ letterwise commutes with $a_2\ldots a_{j-1}$. This gives
		\begin{eqnarray*}
			\alpha a_2 \ldots a_{j} &=& \alpha' (a_{j+1}\sm_j) \text{ from Equation }\eqref{eq:star}\\
			&=&\alpha'(\ba_ja_j)\\
			\Rightarrow 	\alpha a_2 \ldots a_{j-1} &=&\alpha'\ba_j \text{ by cancellation of } a_j
		\end{eqnarray*}
		Now $\alpha (a_2 \ldots a_{j-1}) = \alpha'\ba_j$ and $\ba_j$ letterwise commutes with $a_2\ldots a_{j-1}$. By Lemma \ref{lem:mini} it follows that $\ba_j$ is in $\EndMon_n(\alpha)$, that is there exists $\alpha''$ in $A^+_n$ such that~$\alpha=~\alpha '' \ba_j$.
		
		Then the facet $\lb\alpha a_j\sm_{j-1}\ldots \sm_2\rb_1$ satisfies
		\begin{eqnarray*}
			&&\lb\alpha a_j\sm_{j-1}\ldots \sm_2\rb_1\\
			&=&\lb\alpha''\ba_j a_j\sm_{j-1}\ldots \sm_2\rb_1
		\end{eqnarray*}
		and by Lemma \ref{lem:a' LW coms and face map} $\ba_ja_j$ is a face map~$\cd^{n-j+1}_{i_{j+1}+1}$, so $\ba_j a_j\sm_{j-1}\ldots \sm_2$ is also a face map~$\cd^{n-1}_{i_{j+1}+j-1}$, and therefore $\lb\alpha a_j\sm_{j-1}\ldots \sm_2\rb_1$ is also a facet of~$\lb\alpha''\rb_0$. Since $\lng (\ba_j)\geq~1$ by Lemma~\ref{lem:a' LW coms and face map} it follows $\lng(\alpha'')<\lng(\alpha)$ and so $\lb\alpha a_j\sm_{j-1}\ldots \sm_2\rb_1 \in \C^n(k)$.
		
		For case (b), $a_2=\cdots =a_j=e$ implies $a_{j+1}\sm_j$ is in $\EndMon_n(\alpha)$, so~$\alpha=~\alpha'a_{j+1}\sm_j$ for some $\alpha'$ in $A^+_n$ with $\lng(\alpha')<\lng{(\alpha)}$. Then the facet $\lb\alpha a_j\sm_{j-1}\ldots \sm_2\rb_1$ satisfies
		\begin{eqnarray*}
			&&\lb\alpha a_j\sm_{j-1}\ldots \sm_2\rb_1\\
			&=&\lb(\alpha'a_{j+1}\sm_{j}) a_j\sm_{j-1}\ldots \sm_2\rb_1\\			
			&=&\lb\alpha'(a_{j+1}\sm_{j}\sm_{j-1}\ldots \sm_2)\rb_1 \text{ since } a_j=e
		\end{eqnarray*}
		and as before by Lemma \ref{lem:a' LW coms and face map} this is a face of $\lb\alpha'\rb_0$ which is in $\C^n(k)$ as required.
	\end{proof}	
\end{lem}

\begin{prop}\label{prop:hypothesis}
	If $\EndGen_p(\alpha a)\neq \emptyset$ and $\lb\alpha a\rb_p$ is in $\C^n(k)$ then a facet containing $\lb\alpha a\rb_p$ is in $\C^n(k)$.
	
	\begin{proof}
		Putting together Lemmas \ref{lem:hypothesis1} and \ref{lem:hypothesis2} gives the required result.
	\end{proof}
\end{prop}

\subsection{Proof of point (\ref{point A}): case (\ref{casei}):$\lng(\beta b)<\lng(\alpha a)$} \label{sec: case i}

\begin{prop}
	Under the hypotheses of case (\ref{casei}), $\EndGen_p(\alpha a)\neq \emptyset$.
	\begin{proof}
		Recall that for some $u$ and $v$ in $A^+_p$, $
		\alpha a=\gamma u$ and $\beta b=\gamma v$. If $\lng(\beta b)<\lng(\alpha a)$ then it follows $\lng(\gamma v)<\lng(\gamma u)$ and consequently $\lng(v )<\lng(u)$, since multiplication in $A^+_n$ corresponds to addition of lengths. Since the inequality is strict, it follows that $\lng(u) \neq 0$, i.e.~$u \neq e$. Then since $\alpha a=\gamma u$, $u \in \EndMon_p(\alpha a)$ so in particular $\EndGen_p(\alpha a) \neq \emptyset$.
	\end{proof}
\end{prop}

Applying Proposition \ref{prop:hypothesis} concludes the proof of case (\ref{casei}).

\subsection{Proof of point (\ref{point A}): case (\ref{caseii}): $\lng(\beta b)=\lng(\alpha a)$}

Recall that for some $u$ and $v$ in $A^+_p$, and $\gamma$ in $A^+_n$ with $\EndMon_p(\gamma)=\emptyset$, that
$\alpha a=\gamma u$ and $\beta b=\gamma v$. 
\begin{prop} \label{prop: al a = be b = gam}
If we are in case (\ref{caseii}) then we only need to consider when $\alpha a=\beta b = \gamma$.
	\begin{proof}
		Case (\ref{caseii}) states that $\lng(\beta b)=\lng(\alpha a)$. This implies that $\lng(\gamma u)=\lng(\gamma v)$, which in turn implies $\lng(u)=\lng(v)$ by cancellation. If $u\neq e$ then $\alpha a$ satisfies $\EndGen_p(\alpha a)\neq \emptyset$. Then, by Proposition \ref{prop:hypothesis} a facet containing $\lb \alpha a \rb_p$ lies in~$\C^n(k)$. Therefore we can assume $u=e$, which implies $v=e$ since they have the same length. Therefore $\alpha a= \beta b = \gamma$.
	\end{proof}
\end{prop}

Recall the definition of $c_j$ from Definition \ref{def:c}:
$$
c_j =
\left\{
\begin{array}{ll}
a_j  & \mbox{if } \lng(a_j)\geq \lng(b_j) \\
b_j & \mbox{if } \lng(a_j)< \lng(b_j)
\end{array}
\right.
$$
for $2\leq j \leq p+1$. Recall $c=c_2\ldots c_{p+1}$.
Recall that since $\lng(\beta)<\lng(\alpha)$ then in case~(\ref{caseii}): $\lng(\beta b)=\lng (\alpha a)$ that it follows $\lng(b)>\lng(a)$.

\begin{prop}
	With the notation as above, there exists at least one $j$ for which $c_j=b_j \neq a_j$. Consider the maximal $j$ for which $c_j=b_j \neq a_j$. Then the facet $\lb\alpha a_j \sm_{j-1}\ldots \sm_2\rb_1$ of $\lb\alpha\rb_0$ containing $\lb\alpha a\rb_p$ is in $\C^n(k)$.
	\begin{proof}
		Recall $c=a'a=b'b$ where $a'=a_2'\ldots a_{p+1}'$ and $b'=b_2'\ldots b_{p+1}'$  as in Definition~\ref{def:c}. We fist prove the existence of $j$ in the statement. Note since $\lng(\beta)<\lng(\alpha)$ it follows that $b\neq e$ and so from Lemma \ref{lem:ljs minimal} it follows that $\lng(\beta)=\lng(\alpha)-1$ which gives $\lng(b)=\lng(a)+1$. Since $c=a'a=b'b$, this gives $\lng(a')=\lng(b')+1$ and in particular $\lng(a')\geq 1$. It follows that at least one $a'_j\neq e$ i.e.~$c_j=b_j \neq a_j$.

		Recall also that $\alpha a= \beta b =\gamma$ from Proposition \ref{prop: al a = be b = gam}.
		Therefore $a$ and $b$ are in $\EndMon_n(\alpha a)$ and it follows from Lemma \ref{lem:deltaendsets} that $\Dl(a,b)$ is in $\EndMon_n(\alpha a)$. From Lemma \ref{lem:c is lcm} $\Dl(a,b)=c$ so it follows that $c$ is in $\EndMon_n(\alpha a)$ i.e.~for some~$\alpha'$ in $A_n^+$ with $\lng(\alpha')<\lng(\alpha)$
		$$
		\alpha a = \alpha'(c)=\alpha'(a'a).
		$$
		By cancellation of $a$ we have $\alpha=\alpha'a'$.
			
		Consider the maximal $j$ for which $c_j=b_j \neq a_j$. Then $a_{j+1}'=\cdots=a_{p+1}'=e$, i.e.~$a'=a_2'\ldots a_j'$. It follows that the facet $\lb\alpha a_j \sm_{j-1}\ldots \sm_2\rb_1$ satisfies
		\begin{eqnarray*}
			\lb(\alpha) a_j \sm_{j-1}\ldots \sm_2\rb_1
			&=& \lb(\alpha'a') a_j \sm_{j-1}\ldots \sm_2\rb_1\\
			&=& \lb(\alpha'a_2'\ldots a_j') a_j \sm_{j-1}\ldots \sm_2\rb_1\\
			&=& \lb\alpha'a_2'\ldots( a_j' a_j) \sm_{j-1}\ldots \sm_2\rb_1\\		
			&=& \lb\alpha'a_2'\ldots( c_j) \sm_{j-1}\ldots \sm_2\rb_1\\
			&=& \lb\alpha'a_2'\ldots a'_{j-1} (b_j) \sm_{j-1}\ldots \sm_2\rb_1.
		\end{eqnarray*}
		Recall~$\lng(b_j)=~l_j$, so post multiplication by $b_j \sm_{j-1}\ldots \sm_2$ corresponds to the face map $\cd^{n-1}_{l_j+j-2}$. Therefore $ \lb\alpha a_j \sm_{j-1}\ldots \sm_2\rb_1$ is a facet of $ \lb\alpha'a_2'\ldots a'_{j-1}\rb_0$ and we have that $\lng(\alpha'a_2'\ldots a'_{j-1})< \lng(\alpha)$ since $\alpha=\alpha'a_2'\ldots a'_{j}$ and $\lng(a_j')\geq 1$. Therefore $\lb\alpha a_j \sm_{j-1}\ldots \sm_2\rb_1$ is in $\C^n(k)$.
	\end{proof}
\end{prop}

\subsection{Proof of point (\ref{point A}): case (\ref{caseiii}): $\lng(\beta b)>\lng(\alpha a)$}

Recall that for some $u$ and~$v$ in $A^+_p$, and $\gamma$ in $A^+_n$ with $\EndMon_p(\gamma)=\emptyset$, that
$\alpha a=\gamma u$ and $\beta b=\gamma v$.

\begin{prop} \label{prop: case iii}
	 If we are in case (\ref{caseiii}) then $b \neq e$. Furthermore, we only need to consider the case when $\gamma=\alpha a$ so $\beta b = \gamma v = \alpha a v$. In this case it follows $\EndGen_p(\beta b)\neq \emptyset$.
	\begin{proof}
		Case (\ref{caseiii}) states that $\lng(\beta b)>\lng(\alpha a)$, and note that this can only happen when $b\neq e$ since $\lng(\beta )<\lng(\alpha)$. Recall this implies $\lng(\beta)=\lng(\alpha)-1$ from Lemma~\ref{lem:ljs minimal}. If $u\neq e$ then $\alpha a$ satisfies $\EndGen_p(\alpha a)\neq \emptyset$. Then by Proposition~\ref{prop:hypothesis}, a facet containing $\lb \alpha a \rb_p$ lies in $\C^n(k)$. Therefore we can assume $u=e$. Then $\alpha a = \gamma$ and it follows that $\beta b=\gamma v=\alpha a v$. Since $\lng(\beta b)>\lng(\alpha a)$ it follows $\lng(v)\geq 1$ and therefore $\EndGen_p(\beta b)\neq \emptyset$.
	\end{proof}
\end{prop}

We now prove a technical Lemma required for the rest of this section.

\begin{lem}\label{lem:contracdicting b}
	If there exists~$\beta'\in A_n^+$ such that~$\lng(\beta')=\lng(\beta)-1=\lng(\alpha)-2$ and
	$$
	\lb \beta b \rb_p= \lb\beta' b \rb_p
	$$
	then this contradicts our choice of~$b$.
	\begin{proof}
		Write~$b=\sm'b'$ i.e~$\sm'$ is the leftmost generator of the word $b$. Then~$\lng(\beta'\sm')=\lng(\alpha)-1$ and so~$$\lb \beta'\sm'b'\rb_p=\lb \beta b\rb_p$$
		where~$\lng(b')<\lng (b)$ and~$\lng(\beta'\sm')=\lng(\beta)=\lng(\alpha)-1$. This contradicts our choice of~$b$: we chose $b$ such that $\sum_{k=2}^{p+1}l_k$ was minimal, as in Lemma \ref{lem:ljs minimal} and therefore no such~$b'$ can exist.
	\end{proof}
\end{lem}

\begin{prop} \label{prop: case (iii) prop 1}
	$\EndGen_0(\beta b)= \emptyset$.
	\begin{proof}
		Suppose $\EndGen_0(\beta b)\neq \emptyset$ for a contraction. Let $\tau$ in $\EndGen_0(\beta b)$.
		Then since $\tau$ letterwise commutes with $b_2\ldots b_{p+1}$ it follows that $\tau$ is in $\EndGen_0(\beta)$ from Lemma \ref{lem:mini}. Then $\beta = \beta' \tau$ for some $\beta'$ in $A^+_n$ with $\lng(\beta')=\lng(\beta)-1=\lng(\alpha)-2$. It follows
		\begin{eqnarray*}
			\lb\beta b\rb_p &=& \lb(\beta' \tau) b\rb_p\\
			&=&\lb\beta'  \tau b\rb_p \\
			&=&\lb\beta' b \tau \rb_p \\
			&=& \lb\beta' b \rb_p
		\end{eqnarray*}
		and by Lemma \ref{lem:contracdicting b} this is a contradiction.
		\end{proof}
\end{prop}

\begin{prop}\label{prop: case (iii) prop 2}
	The generator $\sm_1$ is not in $\EndGen_p(\beta b)$.
	\begin{proof}
		Suppose $\sm_1$ is in $\EndGen_p(\beta b)\neq \emptyset$ and work for a contradiction. Since~$\sm_1$ letterwise commutes with $b_3\ldots b_{p+1}$ it follows that $\sm_1$ is in $\EndGen_p(\beta b_2)$ by Lemma \ref{lem:mini}. From Lemma \ref{lem: lcm a_j+1 s_j}, $\Dl(\sm_1, b_2)=b_2\sm_1b_2$ and by Lemma \ref{lem:LWCends} this is in $\EndMon_n(\beta b_2)$, giving by cancellation of $b_2$ that $b_2 \sm_1$ is in $\EndMon_n(\beta)$. So $\beta=\beta'b_2\sm_1$ for some $\beta'$ in $A^+_n$. Then
		\begin{eqnarray*}
			\lb(\beta )(b)\rb_p &=& \lb(\beta' b_2 \sm_1) (b)\rb_p\\
			&=& \lb(\beta' b_2 \sm_1 )(b_2\ldots b_{p+1})\rb_p
		\end{eqnarray*}
		and by Lemma \ref{lem:hat}, $b_2\sm_1 b_2$ can be written as $\hat{b}_1\sm_1 b_2 \sm_1$. So we have
		\begin{eqnarray*}
			\lb(\beta) (b)\rb_p &=& \lb(\beta' b_2 \sm_1) (b_2\ldots b_{p+1})\rb_p\\
			&=& \lb\beta' (b_2 \sm_1 b_2) (b_3\ldots b_{p+1})\rb_p\\
			&=&\lb\beta'(\hat{b}_1\sm_1 b_2 {\color{red}\sm_1}) (b_3\ldots b_{p+1})\rb_p\\
			&=&\lb\beta'(\hat{b}_1\sm_1 b_2) (b_3\ldots b_{p+1}){\color{red}\sm_1}\rb_p\\			
			&=&\lb\beta'(\hat{b}_1\sm_1) (b_2b_3\ldots b_{p+1})\sm_1\rb_p\\
			&=&\lb\beta'\hat{b}_1\sm_1( b )\sm_1\rb_p\\
			&=&	\lb\beta'\hat{b}_1\sm_1 b\rb_p
		\end{eqnarray*}
		with $\lng(\beta'\hat{b}\sm_1)=\lng(\beta)-1$, since~$\beta=\beta'b_2\sm_1$ and~$\hat{b}_1=\sm_{i_2}\ldots \sm_2$ is a subword of~$b_2=~\sm_{i_2+1}\ldots \sm_2$ satisfying~$\lng(\hat{b})=\lng(b)-1$. By Lemma \ref{lem:contracdicting b} this is a contradiction.
	\end{proof}
\end{prop}

\begin{prop}\label{prop: case (iii) prop 3}
	The generator~$\sm_j$ is not in $\EndGen_p(\beta b)$ for~$2\leq j\leq p$.
	\begin{proof}
	Suppose $\sm_j$ is in $\EndGen_p(\beta b)=\EndGen_p(\beta (b_2\ldots b_{p+1}))$ for some $2\leq j \leq p$ and work for a contradiction. Since $\sm_j$ letterwise commutes with $b_{j+2}\ldots b_{p+1}$ it follows from Lemma \ref{lem:mini} that $\sm_j$ is in $\EndGen_p(\beta b_2\ldots b_{j+1})$. From Lemma \ref{lem: lcm a_j+1 s_j}, $\Dl(\sm_j, b_{j+1})=b_{j+1}\sm_jb_{j+1}$ and by Lemma \ref{lem:LWCends} this is in $\EndMon_n(\beta b_2\ldots b_{j+1})$, giving by cancellation of $b_{j+1}$ that $b_{j+1} \sm_j$ is in $\EndMon_n(\beta b_2\ldots b_{j})$. We first handle the two cases where $b_j=e$, namely the case where~$b_{j+1}=e$ and the case where~$b_{j+1}\neq e$.
	
	When~$b_{j+1}=e$, it follows from the conditions on the~$l_i$ that~$b_j=b_{j-1}=\cdots=b_2=e$, so $\sm_j$ is in $\EndMon_n(\beta)$ i.e. there exists~$\beta' \in A^+_n$ such that~$\beta=\beta'\sm_j$ (in particular~$\lng(\beta')=\lng(\beta)-1$. In this case
	\begin{eqnarray*}
	\lb \beta b \rb_p&=&\lb \beta'\sm_j b_{j+2}\ldots b_{p+1} \rb_p\\
	&=&\lb \beta'b_{j+2}\ldots b_{p+1} \sm_j \rb_p\\
	&=& \lb \beta'b_{j+2}\ldots b_{p+1} \rb_p\\
	&=&\lb \beta' b \rb_p
	\end{eqnarray*}
	and by Lemma \ref{lem:contracdicting b} this is a contradiction.
	
	When~$b_j=e$ but~$b_{j+1}\neq e$, it follows that~$b_{j-1}=\cdots=b_2=e$, so $b_{j+1}\sm_j$ is in $\EndMon_n(\beta)$ i.e. there exists~$\beta' \in A^+_n$ such that~$\beta=\beta'b_{j+1}\sm_j$ and therefore
	\begin{eqnarray*}
		\lb \beta b \rb_p&=&\lb \beta'(b_{j+1}\sm_j b_{j+1})\ldots b_{p+1} \rb_p\\
		&=&\lb \beta'(\hat{b}_jb_jb_{j+1}\sm_{j})b_{j+2}\ldots b_{p+1} \rb_p \text{ by Lemma }\ref{lem:hat}\\
		&=&\lb \beta'\hat{b}_jb_{j+1}\sm_{j}b_{j+2}\ldots b_{p+1} \rb_p\text{ since }b_j=e \\
		&=&\lb \beta'\hat{b}_jb_{j+1}b_{j+2}\ldots b_{p+1} \sm_j \rb_p \\
		&=&\lb \beta'\hat{b}_jb \rb_p.
	\end{eqnarray*}
	We note that in this case, since~$b_j=e$,~$\hat{b}_j=\sm_{l_{j+1}+j-1}\ldots\sm_j$ and this has the same length as~$b_{j+1}=\sm_{l_{j+1}+j}\ldots \sm_{j+1}$. Therefore~$\lng(\beta'\hat{b}_j)=\lng(\beta'{b}_{j+1})=\lng(\beta)-1$ and by Lemma~\ref{lem:contracdicting b} this is again a contradiction.
	
	Now assume~$b_j\neq e$. By Lemma \ref{lem:a' LW coms and face map}, $b_{j+1} \sm_j=\bar{b}_jb_j$ with~$\lng(\bar{b}_j)\geq 1$, and so by cancellation of $b_j$, $\bar{b}_j$ is in $\EndMon_n(\beta b_2\ldots b_{j-1})$. From Lemma \ref{lem:mini}, since $\bar{b}_j$ letterwise commutes with $b_2\ldots b_{j-1}$ we have $\bar{b}_j$ is in $\EndMon_n(\beta)$ so $\beta=\beta'\bar{b}_j$ for some $\beta'$ in $A^+_n$. Then it follows that
		\begin{eqnarray*}
			\lb(\beta) b\rb_p &=& \lb(\beta'\bar{b}_j) (b)\rb_p\\
			&=& \lb(\beta'\bar{b}_j) (b_2\ldots b_{p+1})\rb_p\\
			&=&\lb(\beta'\bar{b}_j)(b_2\ldots b_{j-1})b_j(b_{j+1}\ldots b_{p+1})\rb_p\\
			&=&\lb(\beta')(b_2\ldots b_{j-1})(\bar{b}_jb_j)(b_{j+1}\ldots b_{p+1})\rb_p\\
			&=&	\lb\beta'(b_2\ldots b_{j-1})(b_{j+1}\sm_j)(b_{j+1}\ldots b_{p+1})\rb_p \text{    since } \bar{b}_jb_j=b_{j+1}\sm_j\\		
			&=&	\lb\beta'(b_2\ldots b_{j-1})(b_{j+1}\sm_jb_{j+1})(b_{j+2}\ldots b_{p+1})\rb_p \\		
			&=&	\lb\beta'(b_2\ldots b_{j-1})(\hat{b}_jb_jb_{j+1}\sm_j)(b_{j+2}\ldots b_{p+1})\rb_p \text{ by Lemma }\ref{lem:hat}\\
			&=&	\lb\beta'(b_2\ldots b_{j-1})(\hat{b}_j)(b_jb_{j+1})(\sm_j)(b_{j+2}\ldots b_{p+1})\rb_p\\
			&=&	\lb\beta'\hat{b}_j(b_2\ldots b_{j-1})(b_jb_{j+1})\sm_j(b_{j+2}\ldots b_{p+1})\rb_p\\
			&=&	\lb\beta'\hat{b}_j(b_2\ldots b_{j-1}b_jb_{j+1}b_{j+2}\ldots b_{p+1})\sm_j\rb_p\\
			&=&	\lb\beta'\hat{b}_j(b)\sm_j\rb_p\\
			&=& \lb\beta'\hat{b}_jb\rb_p.
		\end{eqnarray*}
		Since $\beta b=\beta'\hat{b}_jb\sm_j$, it follows from the additive property of the length function that~$\lng(\beta'\hat{b}_j)=\lng(\beta)-1$ and so by Lemma~\ref{lem:contracdicting b} this is a contradiction.
	\end{proof}
\end{prop}

By Propositions~\ref{prop: case (iii) prop 1},~\ref{prop: case (iii) prop 2} and \ref{prop: case (iii) prop 3} it follows that~$\EndGen_p(\beta b)=\emptyset$. This contradicts the statement of Proposition~\ref{prop: case iii} and therefore concludes the proof of case (\ref{caseiii}) and hence the proof of point (\ref{point A}).

\subsection{Proof of point (\ref{point B})}

Recall point \ref{point B}: If~$\lng(\alpha)=\lng(\beta)=k+1$ and~$\alpha \neq \beta$ then~$\lb\alpha\rb_0 \cap \lb\beta\rb_0 \subseteq \C^n(k)$.

\begin{prop}
	Suppose $\alpha \neq \beta$ in $A^+_n$. If $\lng(\alpha)=\lng(\beta)=k+1$ then it follows that~$\lb\alpha\rb_0 \cap \lb\beta\rb_0 \subseteq \C^n(k)$.
	\begin{proof}
		Suppose $\lb\alpha\rb_0 \cap \lb\beta\rb_0 \neq \emptyset$. Then  for some~$1 \leq p \leq n-1$ there exists $a$ and $b$ as in Definition~\ref{def:a b} such that $\lb\alpha a\rb_p=\lb\beta b\rb_p$. It follows that there exists~$\gamma$ in $A^+_n$ and $u, v$ in $A^+_p$ such that
		$$
		\alpha a = \gamma u \text{  and  } \beta b = \gamma v.
		$$
		Suppose that $u\neq e$. Then by Proposition~\ref{prop:hypothesis} it follows that a facet of $\lb\alpha\rb_0$ containing $\lb\alpha a\rb_p$ is in $\C^n(k)$. Hence $\lb \alpha a \rb_p=\lb\beta b\rb_p$ itself is in $\C^n(k)$. Similarly if~$v\neq~e$ then a facet of $\lb \beta \rb_0$ containing $\lb\beta b\rb_p=\lb\alpha a\rb_p$ is in $\C^n(k)$, and hence~$\lb\beta b\rb_p=~\lb\alpha a\rb_p$ itself is in $\C^n(k)$. So we are left with the case that $u=v=e$, giving
		$$
		\alpha a = \gamma  =\beta b
		$$
		and since $\lng(\alpha)=\lng(\beta)$ it follows that $\lng(a)=\lng(b)$. Since $\alpha \neq \beta$ it follows $a \neq b$.
		Recall the definition of $c$,~$a'$ and~$b'$ from Definition \ref{def:c}. From Lemma~\ref{lem:c is lcm}~$c=~\Dl(a,b)$ and~$c=a'a=b'b$. Since~$\lng(a)=\lng(b)$ then $\lng(a')=\lng(b')$. Suppose $a'=e$, then~$\lng(a')=~\lng(b')$ gives $b'=e$ and hence $c=a=b$. But $a\neq b$ so it follows that $a'\neq e$ and in particular $\lng (a')\geq1$.
		
		From Lemma \ref{lem:LWCends}, since $a$ and $b$ are in $\EndMon_n(\alpha a)$ it follows that $\Dl(a,b)=c$ is in $\EndMon_n(\alpha a)$, so $\alpha a=\alpha' c = \alpha' (a'a)$ for some $\alpha'$ in $A^+_n$. By cancellation of $a$ we have $\alpha=\alpha'a'$ and $\lng(\alpha')<\lng(\alpha)$. Then
		\begin{eqnarray*}
			\lb\alpha a\rb_p&=&\lb(\alpha'a')a\rb_p\\
			&=&\lb\alpha' c\rb_p
		\end{eqnarray*}
		and $\lb\alpha' c\rb_p$ is in $\C^n(k)$ since $c$ represents a series of face maps originating at $\lb\alpha'\rb_0$, with each face map given by the map corresponding to right multiplication by $c_j$, which is either the face map corresponding to $a_j$ or $b_j$.
	\end{proof}
\end{prop}

This completes the proof of point (\ref{point B}), and hence by Proposition \ref{prop: connectivity of C} it follows that $\| \C^n_\bt \|$ is $(n-2)$ connected.

\section{Proof of Theorem {\ref{thm:A}}}\label{sec:HS proof}

This section proves the required results on the differentials of the spectral sequence introduced in Section \ref{section:setup}, before putting together the results of the previous sections and running the spectral sequence argument to complete the proof of Theorem \ref{thm:A}.

\subsection{Results on face and stabilisation maps}

Recall the definition of the face maps of $\A^n_\bt$ from Definition \ref{def:new ssspace}:
$$
\cd^p_k: \A^n_p \to \A^n_{p-1} \text{ for } 0\leq k \leq p
$$
and given by
\begin{eqnarray*}
	\cd^p_k: \A^n_p &\to& \A^n_{p-1} \\
	\cd^p_k: A^+_n \bb \C^n_p &\to&  A^+_n \bb\C^n_{p-1}
\end{eqnarray*}
where $\cd^p_k$ is induced by the face maps of $\C_\bt^n$, which we recall are a composite of right multiplication of the representative for the equivalence class in~$\C^n_p=A^+(n;n-p-1)$ by~$(\sm_{n-p+k}\sm_{n-p+k -1}\ldots \sm_{n-p+1})$, before the inclusion to the equivalence class in~$\C^n_{p-1}$.

Recall from Lemma \ref{lem:levelequiv} that for each $0\leq p \leq n-1$ there is a homotopy equivalence
$$ A^+_n \ff A^+_{n-p-1} \simeq A^+(n;n-p-~1)= \C^n_p, $$	
\noindent given by the map defined levelwise on the bar construction by
\begin{eqnarray*}
	B_k(A_n^+,A_{n-p-1}^+, *)&\to& A^+(n;n-p-1)\\
	\alpha[m_1,\ldots, m_k] &\mapsto&\balpha
\end{eqnarray*}
where $\alpha \in A^+_n$, $m_i\in A^+_{n-p-1}$ for all $i$ and $\alpha=\balpha \beta$ for $\balpha \in A^+(n;n-p-~1)$ and~$\beta \in A_{n-p-1}^+$.

\begin{defn}\label{def:face maps on double htpy quotient}
	Define the map
	$$d^p_k:A^+_n \bb A^+_n \ff A^+_{n-p-1} \to A^+_n \bb A^+_n \ff A^+_{n-p}$$
	\noindent as the composition of two maps $\iota_p\circ\bar{d}^p_k$. The first map
	$$
	\bar{d}^p_k:
	A^+_n \bb A^+_n \ff A^+_{n-p-1} \to  A^+_n \bb A^+_n \ff A^+_{n-p-1}
	$$
	is given by right multiplication of the central term in the double homotopy quotient by $(\sm_{n-p+k}\sm_{n-p+k-1}\ldots \sm_{n-p+1})$.
	
	The set of $(j,k)$-simplices in $A^+_n \bb A^+_n \ff A^+_{n-p-1}$ is identified with the product~$(A^+_n)^j \times A^+_n \times (A^+_{n-p-1})^k$ and a generic element is given by~$\lbrack a_1,\ldots, a_j \rbrack a \lbrack a'_1,\ldots, a'_k \rbrack$ where $a_i$ and $a$ are in $A^+_n$ and $a'_i$ are in $A^+_{n-p-1}$. The map $\bar{d}^p_k$ acts on this simplex as
	$$
	\bar{d}^p_k([a_1,\ldots, a_j]a[a'_1,\ldots, a'_k])=[a_1,\ldots, a_j]a(\sm_{n-p+k}\sm_{n-p+k -1}\ldots \sm_{n-p+1})[a'_1,\ldots, a'_k]
	$$
	The second map $\iota_p$ is the map
	$$
	\iota_p:A^+_n \bb A^+_n \ff A^+_{n-p-1} \to  A^+_n \bb A^+_n \ff A^+_{n-p}
	$$
	induced by the inclusion $A^+_{n-p-1} \hookrightarrow A^+_{n-p}$. Note that~$\bar{d}^p_0$ is the identity map, and therefore $d_0^p=\iota_p$.
\end{defn}

\begin{lem}
The map $\bar{d}^p_k$ in Definition \ref{def:face maps on double htpy quotient} gives a well defined map on the double homotopy quotient~$A^+_n \bb A^+_n \ff A^+_{n-p-1}$.
	
	\begin{proof}
		Since $(\sm_{n-p+k}\sm_{n-p+k -1}\ldots \sm_{n-p+1})$ letterwise commutes with every word in~$A^+_{n-p-1}$, it follows that $\bar{d}^p_k$ commutes with all face maps of the bi-semi-simplicial space $A^+_n \bb A^+_n \ff A^+_{n-p-1}$. Therefore the map on the central term of each simplex gives a map on the whole bi-semi-simplicial space, and hence its geometric realisation: the double homotopy quotient $A^+_n \bb A^+_n \ff A^+_{n-p-1}$.
	\end{proof}
\end{lem}

\begin{lem}\label{lem:face maps commuting diagram}
		The diagram
		\ms
		
		\centerline{\xymatrix@R=10mm@C=30mm{
				A^+_n \bb A^+_n \ff A^+_{n-p-1}\ar[d]_{d_k^p}\ar[r]^-{\simeq} & A^+_n \bb \C^n_p\ar[d]^{\cd_k^p}\\
				A^+_n \bb A^+_n \ff A^+_{n-p}\ar[r]^-{\simeq} & A^+_n \bb \C^n_{p-1}
		}}
		\ms
		\noindent commutes for all $p\geq0$.
		
		\begin{proof}
			Recall from Lemma \ref{lem:homotopy equivalence levels to previous in sequence} that the horizontal homotopy equivalence is given by the levelwise maps on $(j,k)$-simplices of $A^+_n \bb A^+_n \ff A^+_{n-p-1}$:
			\begin{eqnarray*}
				(A^+_n \bb A^+_n \ff A^+_{n-p-1})_{(j,k)}&\to& (A^+_n\bb \C_p^n)_j\\
				\lbrack a_1,\ldots, a_j\rbrack \alpha\lbrack a'_1,\ldots, a'_k\rbrack &\mapsto& \lbrack a_1,\ldots, a_j\rbrack\balpha
			\end{eqnarray*}
			\noindent where $\alpha$ and $a_i$ are in $A^+_n$, the $a'_i$ are in $A^+_{n-p-1}$, and $\alpha=\balpha \beta$ for $\balpha$ in $A^+(n;n-p-1)$ and $\beta$ in $A_{n-p-1}^+$. Diagram chasing using the definition of $d_k^p$ in Definition~\ref{def:face maps on double htpy quotient} gives that levelwise these maps commute, and so taking homotopy quotients and the corresponding maps induced by these levelwise maps yields the required result.
		\end{proof}
\end{lem}

\begin{lem}\label{cor:facemaps}
	The face maps $\cd^p_k$ of $\A^n_\bt$ are all homotopic to the zeroth face map~$\cd^p_0$.
	\begin{proof}
		The map $\bar{d}^p_k$ restricted to $A^+_n \bb A^+_n$ is $A^+_{n-p-1}$-equivariant, and the same holds for the identity map $\id_{A^+_n \bb A^+_n}$. Applying Proposition \ref{prop: equivariant htpy between maps in double homotopy quotient} to these two maps therefore gives an $A^+_{n-p-1}$-equivariant homotopy between them. It follows that they induce homotopic maps $\bar{d}^p_k$ and $\id_{A^+_n \bb A^+_n \ff A^+_{n-p-1}}$ on $A^+_n \bb A^+_n \ff A^+_{n-p-1}$. Applying the inclusion $\iota_p$ to both maps and the homotopy between them yields a homotopy $h_k$ from $d^p_k$ to $\iota_p$. However $\iota_p$ is precisely the map $d^p_0$, and thus $h_k$ is a homotopy from $d^p_k$ to $d^p_0$ for all $k$. Then the image of $h_k$ under the homotopy equivalence in Lemma \ref{lem:face maps commuting diagram} yields a homotopy from $\cd^p_k$ to the zeroth face map $\cd_0^k$, as required.
		\centerline{\xymatrix@R=10mm@C=30mm{
				A^+_n \bb A^+_n \ff A^+_{n-p-1}\ar@/_.6pc/[d]_{d_k^p}\ar@/^.6pc/[d]_{\simeq}^{d_0^p}\ar[r]^-{\simeq} & A^+_n \bb \C^n_p\ar@/_.6pc/[d]_{\cd_k^p}\ar@/^.6pc/[d]_{\simeq}^{\cd_0^p}\\
				A^+_n \bb A^+_n \ff A^+_{n-p}\ar[r]^-{\simeq} & A^+_n \bb \C^n_{p-1}
		}}
\end{proof}
\end{lem}

\begin{lem}\label{lem:facestab}
	The following diagram commutes up to homotopy.
	
	\centerline{\xymatrix@R=10mm@C=30mm{
			\ar[d]^{s_*} BA_{n-p-1}\ar@{<-}[r]^-{\simeq}_-{f^p}&A^+_n \bb A^+_n \ff A^+_{n-p-1}\ar@/_.6pc/[d]_{d_k^p}\ar@/^.6pc/[d]_{\simeq}^{d_0^p}\ar[r]^-{\simeq} & A^+_n \bb \C^n_p\ar@/_.6pc/[d]_{\cd_k^p}\ar@/^.6pc/[d]_{\simeq}^{\cd_0^p}\\
			BA_{n-p}\ar@{<-}[r]^-{\simeq}_-{f^{p-1}}&A^+_n \bb A^+_n \ff A^+_{n-p}\ar[r]^-{\simeq} & A^+_n \bb \C^n_{p-1}
	}}
	
	\noindent i.e.~Under the homotopy equivalence $\A^n_p \simeq BA^+_{n-p-1}$ of Lemma \ref{lem:homotopy equivalence levels to previous in sequence}, the zeroth face map $\cd^p_0: \A^n_p \to \A^n_{p-1}$  is homotopy equivalent to the map $s_*:BA^+_{n-p-1} \to BA^+_{n-p}$ induced by the stabilisation map $s:A^+_{n-p-1} \hookrightarrow A^+_{n-p}$.
		\begin{proof}
		From the proof of Lemma \ref{cor:facemaps}, the right hand square commutes up to homotopy.
		
		From Lemma~\ref{lem:homotopy equivalence levels to previous in sequence}, the map from the centre to the left is given on the $(j,k)$-simplices of the geometric realisation by
		\begin{eqnarray*}
			f^p_{(j,k)}:(A^+_n \bb A^+_n \ff A^+_{n-p-1})_{(j,k)}  &\to& (* \ff A^+_{n-p-1})_k\\
			\lbrack a_1,\ldots, a_j\rbrack a \lbrack a'_1,\ldots, a'_k\rbrack &\mapsto& *\lbrack a'_1,\ldots, a'_k\rbrack
		\end{eqnarray*}
		where $a$ and $a_i$ are in $A^+_n$ and $a'_i$ is in $A^+_{n-p-1}$.
		The map $d^p_0$ is the map
		$$
		d_0^p:A^+_n \bb A^+_n \ff A^+_{n-p-1} \to  A^+_n \bb A^+_n \ff A^+_{n-p}
		$$
		induced by the inclusion $A^+_{n-p-1} \hookrightarrow A^+_{n-p}$. Restricting this map to $(j,k)$-simplices of the double homotopy quotient gives
		\begin{eqnarray*}
			(d^p_0)_{(j,k)}:(A^+_n \bb A^+_n \ff A^+_{n-p-1})_{(j,k)} &\to&  (A^+_n \bb A^+_n \ff A^+_{n-p})_{(j,k)}\\
			\lbrack a_1,\ldots, a_j\rbrack a \lbrack a'_1,\ldots, a'_k\rbrack &\mapsto& \lbrack a_1,\ldots, a_j\rbrack a \lbrack a'_1,\ldots, a'_k\rbrack
		\end{eqnarray*}
		where $a$ and $a_i$ are in $A^+_n$ and the $a'_i$ are in $A^+_{n-p-1}$, hence $a'_i$ is in $A^+_{n-p}$. Applying this map before the homotopy equivalence to the classifying space gives
		
		\centerline{	\xymatrix@R=15mm@C=20mm {	
				(A^+_n \bb A^+_n \ff A^+_{n-p-1})_{(j,k)} \ar[d]_{f^p_{(j,k)}} \ar[r]^{(d^p_0)_{(j,k)}}& (A^+_n \bb A^+_n \ff A^+_{n-p})_{(j,k)} \ar[d]^{f^{p-1}_{(j,k)}}\\			
				(* \ff A^+_{n-p-1})_k \ar@{=}[d]	& (* \ff A^+_{n-p})_k	\ar@{=}[d]\\
			(	BA^+_{n-p-1})_k \ar@{.>}[r]& (BA^+_{n-p})_k
		}}
		
		\noindent and on a $(j,k)$ simplex this map is given by
		
		\centerline{	\xymatrix@R=15mm@C=20mm {	
				\lbrack a_1,\ldots, a_j\rbrack a \lbrack a'_1,\ldots, a'_k\rbrack \ar@{|->}[d]_{f^p_{(j,k)}} \ar@{|->}[r]^{(d^p_0)_{(j,k)}} & 	\lbrack a_1,\ldots, a_j\rbrack a \lbrack a'_1,\ldots, a'_k\rbrack \ar@{|->}[d]^{f^{p-1}_{(j,k)}} \\			
				{*}\lbrack a'_1,\ldots, a'_k\rbrack\ar@{|.>}[r]	& 	{*}\lbrack a'_1,\ldots, a'_k\rbrack.	\\
		}}
		\noindent We note that the dotted map is precisely the map which defines the natural inclusion $BA^+_{n-p-1}\to BA^+_{n-p}$ under the identification of $* \ff A^+_{r}$ with $BA^+_r$ for all $r$. The natural inclusion is in turn induced by the stabilisation map $A^+_r \overset{s}{\hookrightarrow} A^+_{r+1}$ and so we denote it $s_*$. This gives that the left hand square commutes up to homotopy.
	\end{proof}
\end{lem}

\subsection{Spectral sequence argument}

\setcounter{figure}{0}
In this section we run a first quadrant spectral sequence for filtration of $\|\A^n_\bt \|$, e.g.~see Randal-Williams \cite[2 (sSS)]{RandalWilliams}. Recall the points we proved regarding $\|\A^n_\bt \|$:

\begin{enumerate}
	\item there exist homotopy equivalences $\A^n_p \simeq BA^+_{n-p-1}$ for $p \geq 0$ \label{point 2 END}
	\item there is an $(n-1)$ connected map $\|\phi_\bt \|$ from the geometric realisation of~$\A^n_\bt$ to the classifying space $BA^+_n$:
	$$
	\| \A^n_\bt \| \overset{\|\phi_\bt \|}{\longrightarrow} BA^+_n
	$$
	i.e.~$\|\phi_\bt \|$ induces an isomorphism on homotopy groups $\pi_r$ for $0\leq r \leq (n-2)$, and a surjection for $r=(n-1)$. \label{point 4 END}
\end{enumerate}

The first quadrant spectral sequence of the simplicial filtration of $\|\A_\bt^n\|$ satisfies
$$
E^1_{k,l}=H_l(\A^n_k) \Rightarrow H_{k+l}(\|\A^n_\bt\| ).
$$
By point (\ref{point 2 END}) the left hand side is given by $E^1_{k,l}=H_l(\A^n_k)=H_l(BA^+_{n-k-1})$. The first page of the spectral sequence is depicted in Figure \ref{fig:ss1}. By point (\ref{point 4 END}) $\|\phi_\bt\|$ induces an isomorphism
$$
H_{k+l}(\|\A^n_\bt\|)\cong H_{k+l}(BA_n^+) \text{ when  } (k+l)< n-1$$
and a surjection
$$
H_{k+l}(\|\A^n_\bt\|)\twoheadrightarrow H_{k+l}(BA_n^+) \text{ when  } (k+l)= n-1.
$$

\begin{figure}
	\begin{center}
		\begin{tikzpicture}
		\matrix (m) [matrix of math nodes,
		nodes in empty cells,nodes={minimum width=5ex,
			minimum height=5ex,outer sep=-5pt},
		column sep=1ex,row sep=1ex]{
			\vdots	&  \vdots			& \vdots 			& \vdots 		& \\
			3 	&  H_3(BA^+_{n-1})	& \xleftarrow{d^1} H_3(BA^+_{n-2})  &  \xleftarrow{d^1} H_3(BA^+_{n-3})& \cdots \\
			2 	&  H_2(BA^+_{n-1})	& \xleftarrow{d^1} H_2(BA^+_{n-2})  &  \xleftarrow{d^1} H_2(BA^+_{n-3})& \cdots \\
			1   &  H_1(BA^+_{n-1})	& \xleftarrow{d^1}H_1(BA^+_{n-2})  & \xleftarrow{d^1} H_1(BA^+_{n-3})& \cdots \\
			0   &  H_0(BA^+_{n-1})  & \xleftarrow{d^1}H_0(BA^+_{n-2})  & \xleftarrow{d^1} H_0(BA^+_{n-3})& \cdots \\
			\quad\strut &   0  			&  1  &  2  & \cdots\strut \\};
		\draw[thick] (m-1-1.east) -- (m-6-1.east) ;
		\draw[thick] (m-6-1.north) -- (m-6-5.north) ;
		\end{tikzpicture}
		\caption{The $E^1$ page of the spectral sequence, with the groups identified.}
		\label{fig:ss1}
	\end{center}
\end{figure}

\noindent The differential $d^1$ is given by an alternating sum of face maps in $\A^n_\bt$. By Corollary~\ref{cor:facemaps} the face maps are all homotopic to each other and by Lemma~\ref{lem:facestab} they are all homotopic to the stabilisation map $s_*$, via~$\A^n_p \simeq BA^+_{n-p-1}$. Therefore the alternating sum of face maps in the differential $d^1$ will cancel out to give the zero map when there are an even number of terms, and will give the stabilisation map when there are an odd number of terms, i.e.
\begin{eqnarray*}
	d^1: E^1_{\text{even}, l} \to E^1_{\text{odd}, l} && \text{ odd number of terms, so equals the stabilisation map } s_*\\
	d^1: E^1_{\text{odd}, l} \to E^1_{\text{even}, l} && \text{ even number of terms, so equals the zero map }  0.
\end{eqnarray*}

\noindent This gives the $E^1$ page shown in Figure \ref{fig:ss2}.

\begin{figure}
	\begin{center}
		\begin{tikzpicture}
		\matrix (m) [matrix of math nodes,
		nodes in empty cells,nodes={minimum width=5ex,
			minimum height=5ex,outer sep=-5pt},
		column sep=1ex,row sep=1ex]{
			\vdots	&  \vdots			& \vdots 			& \vdots 		& \\
			3 	&  H_3(BA^+_{n-1})	& \xleftarrow{0} H_3(BA^+_{n-2})  &  \xleftarrow{s_*} H_3(BA^+_{n-3})& \xleftarrow{0} \cdots \\
			2 	&  H_2(BA^+_{n-1})	& \xleftarrow{0} H_2(BA^+_{n-2})  &  \xleftarrow{s_*} H_2(BA^+_{n-3})& \xleftarrow{0} \cdots \\
			1   &  H_1(BA^+_{n-1})	& \xleftarrow{0}H_1(BA^+_{n-2})  & \xleftarrow{s_*} H_1(BA^+_{n-3})& \xleftarrow{0} \cdots \\
			0   &  H_0(BA^+_{n-1})  & \xleftarrow{0}H_0(BA^+_{n-2})  & \xleftarrow{s_*} H_0(BA^+_{n-3})&\xleftarrow{0} \cdots \\
			\quad\strut &   0  			&  1  &  2  & \cdots \strut \\};
		\draw[thick] (m-1-1.east) -- (m-6-1.east) ;
		\draw[thick] (m-6-1.north) -- (m-6-5.north) ;
		\end{tikzpicture}
		\caption{The $E^1$ page of the spectral sequence, with groups and~$d^1$ differentials identified.}
		\label{fig:ss2}
	\end{center}
\end{figure}

We proceed by induction on~$n$, for the sequence of monoids~$A^+_n$, and assume that homological stability holds for previous groups in the sequence.\newline {\bf Inductive Hypothesis:} The map induced on homology by the stabilisation map
$$
H_i(BA_{k-1}^+) \overset{s_*}{\longrightarrow} H_i(BA^+_k)
$$
is an isomorphism for $k>2i$ and is a surjection for $k=2i$ whenever $k<n$.

Here we note that Theorem \ref{thm:A} holds for the base case $n=1$, since we have to check $H_0(BA^+_0) \to H_0(BA^+_1)$ is a surjection, which is true since $BA_n^+$ is connected for all $n$.

\begin{figure}
	\subfigure[$n$ odd]{
				\begin{tikzpicture}[scale=0.5]
		\draw[<->] (-1,9)--(-1,-1)--(9.5,-1);			
		\draw (-2,0) node {$0$};					
		\draw (0,-2) node {$0$};		
		\draw (-2,1) node {$1$};					
		\draw (1,-2) node {$1$};
		\draw (-2,2) node {$2$};					
		\draw (2,-2) node {$2$};	
		\foreach \x in {(-2,4), (0,4), (1,4), (2,4)}	
		\draw \x node {$\vdots$};
		\foreach \x in {(4,-2), (3,0), (3,1), (3,2)}					
		\draw \x node {$\cdots$};	
		\draw(5.7,4) node[rotate=40]  {$\vdots$};			
		\draw(6.6,3.7) node[rotate=40]  {${\color{red}\vdots}$};				
		\draw (6,-2) node[rotate=90]  {$\scriptstyle{n-3}$};			
		\draw (7,-2) node[rotate=90]  {$\scriptstyle{n-2}$};						
		\draw (9,-2) node[rotate=90]  {$\scriptstyle{n}$};							
		\draw(8,-2) node[rotate=90]  {$\scriptstyle{n-1}$};
		\draw (-2.1,6) node {$\scriptstyle{\left \lfloor{\frac{n}{2}}\right \rfloor -2}$};
		\draw (-2.1,7) node {$\scriptstyle{\left \lfloor{\frac{n}{2}}\right \rfloor -1}$};
		\draw (-2,8) node {$\scriptstyle{\left \lfloor{\frac{n}{2}}\right \rfloor }$};
		\foreach \x in {(0,0),(0,1),(0,2), (6,2), (7,1), (7,2), (8,1), (8,2), (0,6), (0,7), (0,8), (1,8), (2,8), (5,2), (5,6), (3,7), (4,7), (6,6), (0,5)}
		\draw [fill] \x circle [radius=.1];		
		\foreach \x in {(1,0), (2,0), (1,1), (1,2), (2,1), (2,2),(7,0), (1,6), (2,6), (3,6), (1,7), (1,5), (2,5), (3, 5),(4,5), (5,5), (6,0), (6,1), (5,1), (5,0), (9,0),(8,0), (9,1),(4,2), (4,1), (4,0), (2,7), (4,6),(6,5)}
		\draw \x node {$\scriptstyle{0}$};
		\draw[color=red, thick] (4.5,2.5) -- +(0,-1)--+(2,-1)--+(2,-2)--+(4,-2)--+(4,-3);		
		\draw[color=red, thick] (2.5,7.5) -- +(0,-1)--+(2,-1)--+(2,-2)--+(4,-2)--+(4,-3);	
		\draw[color=blue, thick] (-0.5,8.5) -- (0.5,8.5)--(0.5,-0.5)--(-0.5,-0.5)--cycle;
		\draw[color=red, thick] (2.5,7.5) -- (0.5, 7.5);
		\end{tikzpicture}}
	\subfigure[$n$ even]{
		\begin{tikzpicture}[scale=0.5]
		\draw[<->] (-1,9)--(-1,-1)--(9.5,-1);			
		\draw (-2,0) node {$0$};					
		\draw (0,-2) node {$0$};		
		\draw (-2,1) node {$1$};					
		\draw (1,-2) node {$1$};
		\draw (-2,2) node {$2$};					
		\draw (2,-2) node {$2$};	
		\foreach \x in {(-2,4), (0,4), (1,4), (2,4)}	
		\draw \x node {$\vdots$};
		\foreach \x in {(4,-2), (3,0), (3,1), (3,2)}					
		\draw \x node {$\cdots$};	
		\draw(5,4) node[rotate=40]  {$\vdots$};			
		\draw(6,3.7) node[rotate=40]  {${\color{red}\vdots}$};				
		\draw (6,-2) node[rotate=90]  {$\scriptstyle{n-3}$};			
		\draw (7,-2) node[rotate=90]  {$\scriptstyle{n-2}$};						
		\draw (9,-2) node[rotate=90]  {$\scriptstyle{n}$};							
		\draw(8,-2) node[rotate=90]  {$\scriptstyle{n-1}$};
		\draw (-2.1,6) node {$\scriptstyle{\left \lfloor{\frac{n}{2}}\right \rfloor -2}$};
		\draw (-2.1,7) node {$\scriptstyle{\left \lfloor{\frac{n}{2}}\right \rfloor -1}$};
		\draw (-2,8) node {$\scriptstyle{\left \lfloor{\frac{n}{2}}\right \rfloor }$};
		\foreach \x in {(0,0),(0,1),(0,2), (6,2), (7,1), (7,2), (8,0), (8,1), (8,2), (0,6), (0,7),  (2,7), (0,8), (1,8), (2,8), (5,2), (5,6), (4,6), (3,7), (4,7), (6,6), (0,5),(6,5)}
		\draw [fill] \x circle [radius=.1];		
		\foreach \x in {(1,0), (2,0), (1,1), (1,2), (2,1), (2,2),(7,0), (1,6), (2,6), (3,6), (1,7), (1,5), (2,5), (3, 5),(4,5), (5,5), (6,0), (6,1), (5,1), (5,0), (9,0), (9,1),(4,2), (4,1), (4,0)}
		\draw \x node {$\scriptstyle{0}$};
		\draw[color=red, thick] (4.5,2.5) -- +(0,-1)--+(2,-1)--+(2,-2)--+(3,-2)--+(3,-3);		
		\draw[color=red, thick] (1.5,7.5) -- +(0,-1)--+(2,-1)--+(2,-2)--+(4,-2)--+(4,-3);	
		\draw[color=blue, thick] (-0.5,7.5) -- (0.5,7.5)--(0.5,-0.5)--(-0.5,-0.5)--cycle;
		\draw[color=red, thick] (1.5,7.5) -- (0.5,7.5);
		\end{tikzpicture}}
	\caption{The $E^2$ page of the spectral sequence, under the inductive hypothesis. To the left of the red line, all groups are zero except at positions~$E^2_{0,l}$ for~$2l<n$ - these are highlighted in blue.}
	\label{fig:spectral sequence E2 page}
\end{figure}
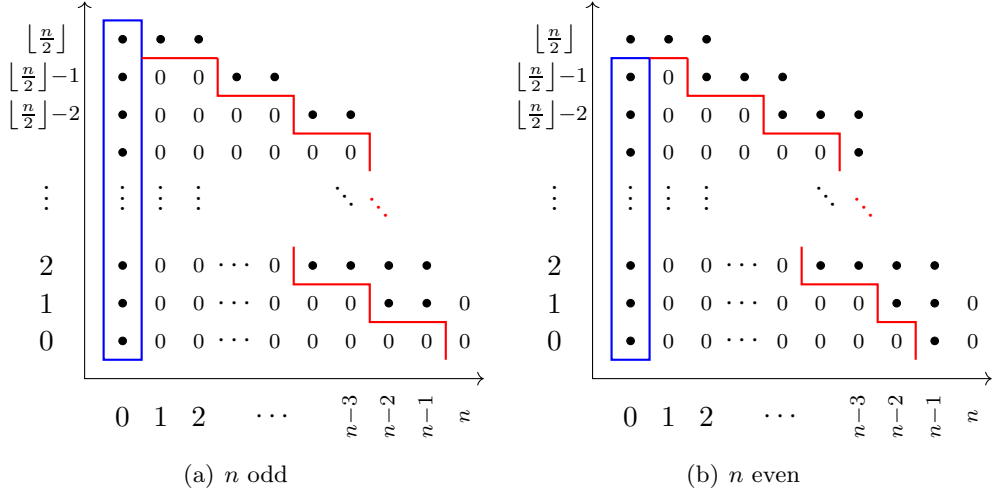

\begin{lem}\label{lem:ss}
	Under the inductive hypothesis, the $E_{0,l}$ terms stabilise on the $E^1$ page for $2l< n$, i.e.
	$$
	E^1_{0,l}=E^\infty_{0,l} \text{ when }2l< n.
	$$
	In particular the $d^1$ differential does not alter these groups, and all possible sources of differentials mapping to $E_{0,l}$ for $2l< n$ are trivial from the $E^2$ page.
	\begin{proof}
		The $d^1$ differentials are given by either the zero map or the stabilisation map as shown in Figure \ref{fig:ss2}. The $d^1$ differentials
		$$
		d^1:E^1_{1,l} \to E^1_{0,l}
		$$
		are given by the zero map, and the  $E^1_{-1,l}$ terms are zero, since this is a first quadrant spectral sequence. Therefore the  $E^2_{0,l}$ terms are equal to the $E^1_{0,l}$ terms.
		
		To show that the sources of all other differentials to $E_{0,l}$ for $2l< n$ are zero, we invoke the inductive hypothesis. This implies that on the~$E^1$ page in the interior of the triangle of height $\lfloor \frac{n}{2} \rfloor$ and base $n$, the stabilisation maps, or $d^1$ differentials satisfy the inductive hypothesis. The resulting maps are shown in Figure~\ref{fig:spectral sequence E1 page}, for the cases~$n$ odd and~$n$ even. Since the $d^1$ differentials going from the odd to the even columns are zero it follows that many groups in the interior of the triangle are zero on the $E^2$ page. This is shown in detail in Figure~\ref{fig:spectral sequence E2 page} for the cases~$n$ odd and~$n$ even. These groups include all the sources of differentials to $E_{0,l}$ for~$2l< n$, hence $E^2_{0,l}=E^{\infty}_{0,l}$ for~$2l< n$.
	\end{proof}
\end{lem}
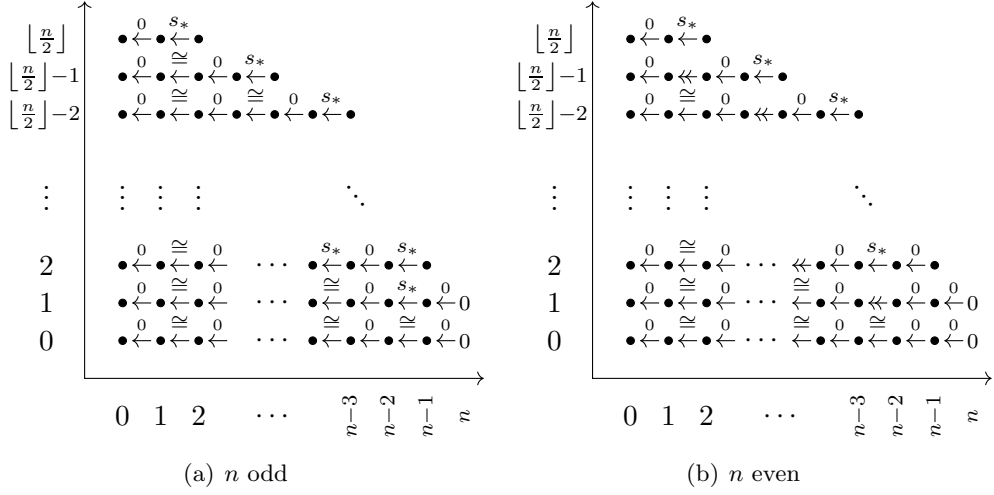
\begin{figure}
	\subfigure[$n$ odd]{
		\begin{tikzpicture}[scale=0.5]
			\draw[<->] (-1,9)--(-1,-1)--(9.5,-1);			
			\draw (-2,0) node {$0$};					
			\draw (0,-2) node {$0$};		
			\draw (-2,1) node {$1$};					
			\draw (1,-2) node {$1$};
			\draw (-2,2) node {$2$};					
			\draw (2,-2) node {$2$};	
			\foreach \x in {(-2,4), (0,4), (1,4), (2,4)}	
			\draw \x node {$\vdots$};
			\foreach \x in {(4,-2), (4,0), (4,1), (4,2)}					
			\draw \x node {$\cdots$};	
			\draw(6,4) node[rotate=40]  {$\vdots$};				
			\draw (6,-2) node[rotate=90]  {$\scriptstyle{n-3}$};			
			\draw (7,-2) node[rotate=90]  {$\scriptstyle{n-2}$};						
			\draw (9,-2) node[rotate=90]  {$\scriptstyle{n}$};							
			\draw(8,-2) node[rotate=90]  {$\scriptstyle{n-1}$};
			\draw (-2.1,6) node {$\scriptstyle{\left \lfloor{\frac{n}{2}}\right \rfloor -2}$};
			\draw (-2.1,7) node {$\scriptstyle{\left \lfloor{\frac{n}{2}}\right \rfloor -1}$};
			\draw (-2,8) node {$\scriptstyle{\left \lfloor{\frac{n}{2}}\right \rfloor }$};
			\foreach \x in {(0,0),(0,1), (0,2), (1,0), (1,1), (1,2), (2,0), (2,1), (2,2), (6,0), (6,1), (6,2), (7,0), (7,1), (7,2), (8,0), (8,1), (8,2), (0,6), (1,6), (2,6), (0,7), (1,7), (2,7), (0,8), (1,8), (2,8), (5,1), (5,0), (5,2), (5,6), (4,6), (3,6), (3,7), (4,7), (6,6)}
			\draw [fill] \x circle [radius=.1];
			\foreach \x in {(9,0), (9,1)}
			\draw \x node {$\scriptstyle{0}$};
			\foreach \x in {(5.75,1), (7.75,0), (5.75,0), (1.75,0), (1.75, 1), (1.75, 2), (1.75, 7), (1.75,6),(3.75, 6)}
			\draw[->] \x -- node[above] {$\scriptstyle{\cong}$} +(-0.5,0);
			\foreach \x in {(6.75,1), (6.75,0), (6.75, 2), (8.75,0), (8.75,1), (0.75,0), (0.75,1), (0.75,2), (0.75, 6), (0.75, 7), (0.75,8), (2.75, 0), (2.75,1), (2.75,2), (2.75,6),(2.75,7), (4.75,6)}
			\draw[->] \x -- node[above] {$\scriptscriptstyle{0}$} +(-0.5,0);
			\foreach \x in {(7.75,1), (7.75,2), (5.75,2),(1.75,8), (3.75,7), (5.75,6)}
			\draw[->] \x -- node[above] {$\scriptstyle{s_*}$} +(-0.5,0);
		\end{tikzpicture}}
	 \subfigure[$n$ even]{
		\begin{tikzpicture}[scale=0.5]
	\draw[<->] (-1,9)--(-1,-1)--(9.5,-1);			
	\draw (-2,0) node {$0$};					
	\draw (0,-2) node {$0$};		
	\draw (-2,1) node {$1$};					
	\draw (1,-2) node {$1$};
	\draw (-2,2) node {$2$};					
	\draw (2,-2) node {$2$};	
	\foreach \x in {(-2,4), (0,4), (1,4), (2,4)}	
	\draw \x node {$\vdots$};
	\foreach \x in {(4,-2), (3.5,0), (3.5,1), (3.5,2)}					
	\draw \x node {$\cdots$};	
	\draw(6,4) node[rotate=40]  {$\vdots$};				
	\draw (6,-2) node[rotate=90]  {$\scriptstyle{n-3}$};			
	\draw (7,-2) node[rotate=90]  {$\scriptstyle{n-2}$};						
	\draw (9,-2) node[rotate=90]  {$\scriptstyle{n}$};							
	\draw(8,-2) node[rotate=90]  {$\scriptstyle{n-1}$};
	\draw (-2.1,6) node {$\scriptstyle{\left \lfloor{\frac{n}{2}}\right \rfloor -2}$};
	\draw (-2.1,7) node {$\scriptstyle{\left \lfloor{\frac{n}{2}}\right \rfloor -1}$};
	\draw (-2,8) node {$\scriptstyle{\left \lfloor{\frac{n}{2}}\right \rfloor }$};
	\foreach \x in {(0,0),(0,1), (0,2), (1,0), (1,1), (1,2), (2,0), (2,1), (2,2), (6,0), (6,1), (6,2), (7,0), (7,1), (7,2), (8,0), (8,1), (8,2), (0,6), (1,6), (2,6), (0,7), (1,7), (2,7), (0,8), (1,8), (2,8), (5,1), (5,0), (5,2), (5,6), (4,6), (3,6), (3,7), (4,7), (6,6)}
	\draw [fill] \x circle [radius=.1];
	\foreach \x in {(9,0), (9,1)}
	\draw \x node {$\scriptstyle{0}$};
	\foreach \x in { (6.75,0),(1.75,0), (1.75, 1), (1.75, 2), (1.75,6), (4.75,1), (4.75,0)}
	\draw[->] \x -- node[above] {$\scriptstyle{\cong}$} +(-0.5,0);
	\foreach \x in {(7.75,1), (7.75,0), (7.75, 2), (8.75,0), (8.75,1), (0.75,0), (0.75,1), (0.75,2), (0.75, 6), (0.75, 7), (0.75,8), (2.75, 0), (2.75,1), (2.75,2), (2.75,6),(2.75,7), (4.75,6), (5.75,2), (5.75, 1), (5.75,0)}
	\draw[->] \x -- node[above] {$\scriptscriptstyle{0}$} +(-0.5,0);
	\foreach \x in {(1.75,8), (3.75,7), (5.75,6), (6.75,2)}
	\draw[->] \x -- node[above] {$\scriptstyle{s_*}$} +(-0.5,0);
	\foreach \x in {(6.75,1), (4.75,2),(1.75, 7), (3.75,6)}
	\draw[->>] \x -- +(-0.5,0);
	\end{tikzpicture}}
\caption{The $E^1$ page of the spectral sequence, with possible non-zero groups represented as circles and the inductive hypothesis applied to the~$d^1$ differentials.}
\label{fig:spectral sequence E1 page}
\end{figure}

We are now in a position to prove Theorem \ref{thm:A}.

\begin{thm}
	The sequence of monoids $A_n^+$ satisfies homological stability, that is
	$$
	H_i(BA_{n-1}^+)\cong H_i(BA_n^+)
	$$
	when $2i< n$, and the map $H_i(BA_{n-1}^+)\to H_i(BA_n^+)$ is surjective when $2i=n$.
	\begin{proof}
		From Lemma \ref{lem:ss}, the spectral sequence satisfies
		$$
		E^\infty_{0,i}=E^1_{0,i}=H_i(BA^+_{n-1})
		$$
		when $2i< n$. From Proposition \ref{prop:highcon} and Theorem \ref{thm:highly connected}
		$$
		H_i(\|\A_\bt^n \|)\cong H_i(BA_n^+)
		$$
		when $i\leq n-2$, and the map $H_i(\|\A_\bt^n \|)\to H_i(BA_n^+)$ is onto when $i=n-1$. The spectral sequence abuts to $H_{k+l}(\|\A^n_\bt\| )$ and from Figure~\ref{fig:spectral sequence E2 page} the only non zero groups on the diagonal $E^{\infty}_{k,l}$ when $k+l=i$ and~$2i<n$ are the groups $E^{\infty}_{0,i}$. Putting these results together yields
		$$
		H_i(BA^+_{n-1})=E^\infty_{0,i}=H_{i+0}(\|\A_\bt^n \|)\cong H_i(BA_n^+)
		$$		
		when both $i< \frac{n}{2}$ and $i\leq n-2$ are satisfied. When $n\geq 2$, $i< \frac{n}{2}$ implies $i\leq n-2$ and the case $n=1$ was the base case of the inductive hypothesis. Therefore an isomorphism is induced when $2i< n$.
		
		\noindent When $i\leq n-1$ and $i< \frac{n}{2}$ it follows that
		$$
		H_i(BA^+_{n-1})=E^\infty_{0,i}=H_{i+0}(\|\A_\bt^n \|)\twoheadrightarrow H_i(BA_n^+)
		$$		
		and for $n\geq 2$, $i< \frac{n}{2}$ implies $i\leq n-1$. Again the case $n=1$ was the base case of the inductive hypothesis. This gives the required range for the surjection, and hence completes the proof.
	\end{proof}
\end{thm}

\bibliography{mybib}{}
\bibliographystyle{alpha}

\end{document}